\documentclass[12pt]{book}
\usepackage{makeidx}
\usepackage{theorem}
\usepackage{amsmath,amssymb,amscd,latexsym,eufrak}

\newcommand{\bvskip}{\vspace{7mm}}
\newcommand{\mvskip}{\vspace{5mm}}
\newcommand{\svskip}{\vspace{3mm}}

\newcommand{\C}{{\Bbb C}}
\newcommand{\R}{{\Bbb R}}

\newcommand{\Z}{{\Bbb Z}}
\newcommand{\BP}{{\Bbb P}}
\newcommand{\Q}{{\Bbb Q}}
\newcommand{\Supp}{{\rm Supp}\:}

\newcommand{\Bk}{{\rm Bk}\:}
\newcommand{\ML}{{\rm ML}\:}

\newcommand{\red}{{\rm red}\:}

\newcommand{\QED}{{\unskip\nobreak\hfil\penalty50\quad\null\nobreak\hfil
{Q.E.D.}\parfillskip0pt\finalhyphendemerits0\par\medskip}}
\newcommand{\Proof}{\noindent{\bf Proof.}\quad}
\newcommand{\A}{{\Bbb A}}
\newcommand{\Spec}{{\rm Spec}\:}
\newcommand{\Proj}{{\rm Proj}\:}
\newcommand{\lto}{\longrightarrow}
\newcommand{\dlto}{\:\cdots\!\!\to}

\newcommand{\lkd}{\ol{\kappa}}

\newcommand{\codim}{{\rm codim}\:}
\newcommand{\height}{{\rm ht}\:}

\newcommand{\trdeg}{{\rm tr.deg}\:}
\newcommand{\rank}{{\rm rank}\:}
\newcommand{\Aut}{{\rm Aut}\:}
\newcommand{\Pic}{{\rm Pic}\:}

\newcommand{\Ker}{{\rm Ker}\:}

\newcommand{\id}{{\rm id}\:}

\newcommand{\propersupset}{\raisebox{-.7ex}{$\,\,\buildrel\displaystyle
{\supset}\over{\scriptstyle\neq}\,\,$}}

\newcommand{\SL}{{\mathcal L}}
\newcommand{\SM}{{\mathcal M}}

\newcommand{\SD}{{\mathcal D}}

\newcommand{\SO}{{\mathcal O}}
\newcommand{\SU}{{\mathcal U}}

\newcommand{\ga}{\EuFrak{a}}

\newcommand{\gp}{\EuFrak{p}}
\newcommand{\gq}{\EuFrak{q}}

\newcommand{\GP}{\EuFrak{P}}

\newcommand{\st}[1]{\stackrel{{#1}}{\longrightarrow}}

\newcommand{\dps}[1]{\displaystyle{#1}}

\newcommand{\wt}{\widetilde}
\newcommand{\ol}{\overline}

\newcommand{\wh}{\widehat}

\newcommand{\is}[2]{({#1}\cdot{#2})}
\newcommand{\sis}[1]{({#1}^2)}
\newtheorem{thm}{Theorem}[section]
\newtheorem{lem}[thm]{Lemma}
\newtheorem{prop}[thm]{Proposition}
\newtheorem{cor}[thm]{Corollary}
\newtheorem{conj}[thm]{Conjecture}
\newtheorem{remark}[thm]{{\sc Remark}}

\newtheorem{example}[thm]{{\sc Example}}
\newtheorem{defn}[thm]{{\sc Definition}}

\makeindex

\usepackage[linkbordercolor={0 0 1}]{hyperref}

\begin{document}
$\mbox{}$
\vspace{1cm}

\begin{center}
{\Huge Lectures\\ 
on \\ 
Geometry and Topology \\ 
of Polynomials \\ 
{\LARGE--Surrounding the Jacobian Conjecture--}\\
\vspace{4cm}

Masayoshi Miyanishi}
\end{center}


\newpage
${}$
\newpage
\pagenumbering{roman}

\begin{center}
{\LARGE Introduction}
\end{center}
\bvskip

A polynomial ring over a field is very simple to define but yet very difficult to handle in 
the case with two or more variables. Although polynomial rings modulo ideals give rise to affine algebras, 
which are fundamental in algebraic geometry, complexity of polynomial computations bars trials to looking 
into the properties of polynomials. For example, one would like to know the relations between the form 
of a polynomial standardized after some operations and geometric or topological properties of the 
hypersurface defined by the given polynomial. Some scattered results are obtained through various 
investigations around this direction, e.g., polynomial mappings defined by the given polynomial 
and other topics. 

Recent developments of affine algebraic geometry, especially the theory of open algebraic surfaces, 
provide means to systematically explore geometric and topological properties of polynomials in two 
variables. Nevertheless, there is one unsurmountable problem remained even in the case of two variables, 
that is the Jacobian conjecture. It has been unsolved for more than sixty years since O.H. Keller 
[Ganze Cremona-Transformationen, Monats. Math. Physik {\bf 47} (1939), 299--306] proved erroneously 
that a polynomial endomorphism $\varphi$ of the complex affine plane $\C^2$ (which we denote also by 
$\A^2$) with everywhere non-vanishing Jacobian determinant is an automorphism. 

Difficulty of the conjecture lies probably in the point that the affine plane $\C^2$ is so simple and 
there are no clues to tackle the problem geometrically. The following four approaches or extra 
conditions to be put additionally will explain what we need or expect to have as clues in dimension two. 
In order to explain these approaches, we note first that the condition on the Jacobian determinant to be nowhere vanishing is equivalent to the condition that $\varphi$ is \'etale (or unramified if there are 
singular points on varieties concerned). 
\svskip

\noindent
(1)\ \ Let $G$ be a finite subgroup of $\Aut(\C^2)$. Suppose that the endomorphism $\varphi$ commutes with 
the $G$-action. Then $\varphi$ induces an endomorphism $\ol{\varphi}$ of the quotient surface $\C^2/G$ 
such that $\ol{\varphi}$ is unramified. Since $G$ is regarded as a linear subgroup, the quotient 
surface $\C^2/G$ has a nice structure of a Platonic $\C^*$-fiber space if the unique singular point 
removed, which is, by definition, a smooth algebraic surface with a $\C^*$-fibration over the projective 
line $\BP^1$ and with three multiple fibers whose multiplicity sequence is a Platonic triplet. If one can show that $\ol{\varphi}$ (or its restriction onto $\C^2/G-\{\mbox{a singular point}\}$) is an 
automorphism, then $\varphi$ itself is an automorphism. This is the so-called {\em equivariant 
version} of the Jacobian conjecture in dimension two.
\svskip

\noindent
(2)\ \ Let $C$ be a curve on $\C^2$. Suppose that $\varphi$ satisfies the condition $\varphi^{-1}(C) 
\subseteq C$. Then $\varphi$ induces an \'etale endomorphism of $\C^2\setminus C$. It is shown in 
\cite{Aoki, KM} that if $C$ is irreducible the conjecture is affirmative with this extra condition.
But the case where $C$ is reducible is not worked out.
\svskip

\noindent
(3)\ \ Let $X$ be a smooth affine surface and let $\varphi$ is an \'etale endomorphism of $X$. One 
can ask as an analogy of the Jacobian conjecture if $\varphi$ is an automorphism. But this is not quite 
so as the $n$-th power mapping of the algebraic torus group $G_m$ is an \'etale finite endomorphism of 
degree $n$. So we ask instead if $\varphi$ is a finite morphism, and we call this {\em the generalized 
Jacobian conjecture}. It will be reviewed later that the generalized Jacobian conjecture is again false 
in general as there are many counterexamples. But one can expect that the conjecture holds for surfaces 
which are quite similar to $\C^2$. One of the candidate surfaces is the complement $\BP^2\setminus C$ of 
a plane curve $C$ defined by $X_0X_1^{d-1}=X_2^d$ for $d \ge 2$. This surface is a kind of 
$\Q$-homology plane and has a structure which we generalize as an {\em affine pseudo-plane}. Namely 
it has an $\A^1$-fibration over the affine line $\A^1$ which has a unique multiple fiber. We shall 
discuss in the subsection 1.4 some results on affine pseudo-planes concerning the cancellation problem. 
See also \cite{tD1, KM3}. It is strongly desired to verify if the generalized Jacobian conjecture holds 
for affine pseudo-planes.
\svskip

\noindent
(4)\ \ Suppose that there exists a non-isomorphic, \'etale endomorphism $\varphi$ of $\C^2$. We are 
interested in algebraic surfaces $Y$ which might factorize $\varphi$ in the sense that 
$\varphi : \C^2 \st{\varphi_1} Y \st{\varphi_2} \C^2$. Note that $\varphi_1$ and $\varphi_2$ are 
\'etale and we may assume that they are almost surjective in the sense that the image includes all 
codimension one points. In general, given a morphism $f : X \to Y$ of smooth affine varieties, we say 
that $X$ is an {\em affine pseudo-covering} of $Y$ if $f$ is almost surjective and \'etale. So, we 
are interested in what kind of surfaces are affine pseudo-coverings of $\C^2$ or what kind of surfaces 
have $\C^2$ as affine pseudo-coverings. These points of view are taken up in the subsection 2.4.
\svskip

The author started to write the present lecture notes for the lecture(s) to be delivered at the 
Mathematics Department of l'Universit\'e de Bordeaux I and, in fact, talked on the beginning parts 
1.1 and 1.2 in January, 2003. The most of the first section was lectured at the Department of Mathematics, Osaka University during the first semester of the year 2003. There are some influences in the notes 
by lecturing in the classes. For example, the subsections 1.1 and 1.2 as well as the other parts are 
elementary, well-known or short of full explanations. But we left these parts in the notes as they are 
because the notes grew as the author teaches in the classes and thinks the Jacobian problem once and 
again. So, the readers might as well skip 1.1, 1.2 and 1.3. The author gradually began to have the idea of 
putting together various results of the author and his collaborators which were published in journals 
and proceedings \cite{Aoki, GM2, GM, MM, ME, KM}. But these notes are not simple reproductions of 
the past results. There are some new results contained in the subsections 1.4 and 2.5. Meanwhile, most 
results in the subsections 2.2, 2.3 and 2.4 are the reproductions from \cite{KM}. We added the 
subsection 2.1 for the readers convenience to have a quick overview on the present (not necessarily 
updated) status of the Jacobian conjecture. 

The author would like to express his sincere gratitude to Professors Phillipe and Pierrette Cassou-Nogues 
of l'Universit\'e de Bordeaux I for their kind hospitality during the stay of the author in Bordeaux and 
to Professor R.V. Grujar for mathematical discussions on the e-mails. Furthermore, the author would like 
to thank Professor A. Fujiki who advised the author to publish these notes as one volume in the 
lecture note series of the Department of Mathematics, Osaka University.
\mvskip

\begin{flushright}
November, 2003
\svskip

M. Miyanishi
\end{flushright}
\newpage
\pagenumbering{arabic}

\chapter{Characterization of polynomial rings}
\bvskip

{\large
\begin{enumerate}
\item[\S 1.]
General properties of polynomial rings
\item[\S 2.]
Squeezing-out one variable
\item[\S 3.]
Fibrations of the affine plane by polynomials
\item[\S 4.]
Plane-like affine surfaces
\end{enumerate}}

\newpage
\section{General properties of polynomial rings}

\mvskip

Let $k$ be a field. A polynomial ring $k[x_1, \ldots, x_n]$ in $n$ variables over $k$ is a 
commutative $k$-algebra consisting of finite sums of monomials 
$a_{d_1\cdots d_n}x_1^{d_1}\cdots x_n^{d_n}$, where 
$a_{d_1\cdots d_n} \in k$ and monomials are multiplied by the rule
\[
x_1^{d_1}\cdots x_n^{d_n}\times x_1^{e_1}\cdots x_n^{e_n}=x_1^{d_1+e_1}\cdots x_n^{d_n+e_n}.
\]
Hence two polynomials $f=\sum_da_{d_1\cdots d_n}x_1^{d_1}\cdots x_n^{d_n}$ and $g =
\sum_db_{d_1\cdots d_n}x_1^{d_1}\cdots x_n^{d_n}$ are added and multiplied as follows:
\[
f+g=\sum_d(a_{d_1\cdots d_n}+b_{d_1\cdots d_n})x_1^{d_1}\cdots x_n^{d_n},
\]
\[
fg=\sum_c\sum_{c=d+e}a_{d_1\cdots d_n}b_{e_1\cdots e_n}x_1^{c_1}\cdots x_n^{c_n}.
\]

For a ring $R$, we denote by $R^*$ the set of invertible elements of $R$. Note that $R^*$ is 
a group under the multiplication of $R$. For an element $a \in R^*$, the group inverse is the 
inverse $a^{-1}$ of $a$ in $R$. For an integral domain $R$, an element $a$ is {\em irreducible}
\index{element!irreducible} if $a=bc$ with $b,c \in R$ implies either $b \in R^*$ or $c \in R^*$. 
$R$ is said to be {\em decomposable}\index{ring!decomposable} if for any nonzero element $a \in R$, 
there is at least one irreducible decomposition $a=ub_1\cdots b_n$, where $u \in R^*$ and 
$b_1, \ldots, b_n$ are irreducible elements.

\begin{lem}\label{Lemma 1.1.1}
If $R$ is a noetherian domain then $R$ is decomposable. 
\end{lem}
\Proof
Suppose that an element $a$ has no irreducible decomposition. Then $a$ is not irreducible. 
Hence $a=a_1b_1$ with non-invertible elements $a_1, b_1 \in R$. Then either $a_1$ or $b_1$ has no irreducible decomposition, say $a_1$. Repeating this argument, we find an infinite series of 
non-invertible elements $a, a_1, a_2, \cdots $ such that $a_i=a_{i+1}b_{i+1}$ and $a_i$ has no 
irreducible decomposition. This implies that there is an ascending chain of ideals 
$aR \subset a_1R \subset \cdots \subset a_nR \subset\cdots$. This contradicts the hypothesis that 
$R$ is noetherian.
\QED
 
Finally, $R$ is said to be {\em factorial}\index{affine $k$-domain!factorial} or a 
{\em unique factorization domain}\index{unique factorization domain} (UFD) if $R$ is 
decomposable and if whenever there exist two expressions $a=ub_1\cdots b_n=vc_1\cdots c_m$ with 
$u, v \in R^*$ then $n=m$ and $b_{\sigma(i)}=u_ic_i \ (1 \le i \le n)$ for a suitable permutation 
$\sigma$ of $\{1,2, \ldots, n\}$ and $u_i \in R^*$. 

Let $R$ be an integral domain. Let $a, b \in R$. Then $a$ {\em divides} $b$ if $b=ac$ 
for some $c \in R$. An element $p \in R$ is called a {\em prime}\index{element!prime} element if 
$p$ divides $ab$ with $a, b \in R$ then $p$ divides either $a$ or $b$. This definition is equivalent 
to saying that the ideal $aR$ is a prime ideal, i.e., the ring $R/(a)$ is an integral domain. 
A prime element $p$ is irreducible, for if $p=ab$ then $a=pu$ or $b=pv$ with $u, v \in R$, 
and $ub=1$ or $av=1$, i.e., either $u$ or $v$ is invertible. 
\svskip

\noindent
{\bf Exercise 1.1.1.}\ \ Let $R$ be an integral domain. Show that if $R$ is factorial then every 
irreducible element is a prime element. Suppose further that $R$ is decomposable. If every irreducible 
element is a prime element, then $R$ is factorial.
\svskip

Suppose that $R$ is noetherian and hence decomposable. Given two elements $a, b$ of $R$, we say that 
$c \in R$ is a {common divisor} of $a, b$ if $c$ divides $a$ and $b$, $c \mid a$ and $c\mid b$ by 
notation. An element $d$ of $R$ is a greatest common divisor of $a, b$ if $d$ is a common divisor 
of $a, b$ and any common divisor $e$ of $a, b$ divides $d$. If $R$ is factorial, such an element $d$ 
exists and is uniquely determined up to  multiplication of an invertible element. So, we denote $d$ 
by $\gcd(a,b)$. For elements $a_1, \ldots, a_m$ of $R$, we can also define the greatest common divisor 
$\gcd(a_1,\ldots,a_m)$ if $R$ is factorial.
\svskip

\noindent
{\bf Exercise 1.1.2.}\ \ Suppose $R$ is noetherian and factorial. For elements $a_1, \ldots, \\
a_m$ of $R$, write
\[
a_i=u_ip_1^{\alpha^{(i)}_1}p_2^{\alpha^{(i)}_2}\cdots p_n^{\alpha^{(i)}_n},
\]
where $u_i \in R^*, \alpha^{(i)}_j \ge 0 \ (1 \le j \le n)$ and $p_1, p_2, \cdots, p_n$ are mutually 
non-divisible prime elements. Show then that $\gcd(a_1,\ldots, a_m)$ is given by
\[
p_1^{\beta_1}p_2^{\beta_2}\cdots p_n^{\beta_n},
\]
where $\beta_j=\min (\alpha^{(j)}_1, \ldots, \alpha^{(j)}_m)$.

\begin{thm}\label{Theorem 1.1.2}
The polynomial ring $k[x_1,\ldots,x_n]$ has the following properties.
\begin{enumerate}
\item[{\rm (1)}]
$k[x_1,\cdots,x_n]^*=k^*$.
\item[{\rm (2)}]
$k[x_1,\ldots,x_n]$ is factorial.
\end{enumerate}
\end{thm}
\Proof
The proof consists of showing that (1) a polynomial ring $k[x]$ is factorial and (2) $R[x]$ is 
noetherian and factorial provided so is $R$. The assertion (1) is Exercise 1.1.3 below. We shall 
prove the assertion (2). Let $K$ be the quotient field of $R$. Note that $R[x] \subset K[x]$. Let 
$f(x) =a_0x^d+a_1x^{d-1}+\cdots+a_{d-1}x+a_d$ be a polynomial in $R[x]$. We say that $f$ is 
{\em primitive}\index{polynomial!primitive } if the coefficients $a_0, a_1, \ldots, a_d$ 
has no common divisor. If 
$f(x)=a_0x^n+a_1x^{n-1}+\cdots+a_n$ and $g(x)=b_0x^m+b_1x^{m-1}+\cdots+b_m$ are primitive polynomials 
with $a_0b_0\ne 0$, then the product $f(x)g(x)=\sum_{i=0}^{n+m}(\sum_{j=0}^ia_jb_{i-j})x^{n+m-i}$ is 
primitive as well. In fact, if $f(x)g(x)$ is not primitive, then there is a prime element $p$ of $R$ 
such that $p$ divides all the coefficients of $f(x)g(x)$. Suppose that $p \mid a_j$ for $0 \le j < \ell_1, 
p \not\! | \; a_{\ell_1}$ and $p \mid b_j$ for $0 \le j < \ell_2$ and $p \not\! | \; b_{\ell_2}$. 
Then $p$ does not divide the coefficient $a_{\ell_1}b_{\ell_2}+
(\sum_{i=0}^{\ell_1+\ell_2}a_ib_{\ell_1+\ell_2-i}-a_{\ell_1}b_{\ell_2})$ 
of the term of degree $\ell_1+\ell_2$. This is a contradiction. 

Given a polynomial $h(x)$ of $K[x]$, we find elements $c, d \in R$ and a primitive polynomial 
$f(x)\in R[x]$ such that $\gcd(c,d)=1$ and $ch(x)=df(x)$. Note that $c=1$ if $h(x) \in R[x]$. 
Furthermore, if $h(x)$ is an irreducible element of $R[x]$, then $d=1$. Hence $h(x)$ is primitive. 
Let $h(x)$ be an irreducible element of $R[x]$ again. If $h(x)$ splits as $h(x)=f'(x)g'(x)$ in 
$K[x]$ with $\deg f'(x) >0$ and $\deg g'(x) > 0$, then $ch(x)=df(x)g(x)$, where $c,d \in R$ with 
$\gcd(c,d)=1$ and $f(x), g(x) \in R[x]$ which are primitive and differ from $f'(x), g'(x)$ by elements 
of $K$. Since $f(x)g(x)$ is primitive, it follows that $c=1$. Hence $h(x)$ is not irreducible. This 
implies that an irreducible element of $R[x]$ is irreducible in $K[x]$ as well. 
We shall show that if $f(x)$ is an irreducible element of $R[x]$ then $f(x)$ is a prime element. 
In fact, if $f(x)$ is constant, i.e., $f(x)=a \in R$, then $a$ is irreducible in $R$ and hence prime. 
So, the ring $R[x]/(f(x))$ is isomorphic to $(R/(a))[x]$, which is an integral domain. So, $f(x)$ is 
a prime element in $R[x]$. Assume that $f(x) \not\in R$. Suppose that $g(x)h(x)=f(x)q(x)$, where 
$g(x), h(x), q(x) \in R[x]$. Since $f(x)$ is irreducible, $f(x)$ divides either $g(x)$ or $h(x)$ 
in $K[x]$. Suppose $f(x) \mid g(x)$ in $K[x]$. Write $g(x)=f(x)q'(x)$ with $q'(x) \in K[x]$ and 
write $cq'(x)=dq(x)$, where $c, d \in R$ with $\gcd(c,d)=1$ and $q(x) \in R[x]$ is primitive. Then 
$cg(x)=df(x)q(x)$, where $f(x)q(x)$ is primitive. Hence $c=1$, and $g(x)=df(x)q(x)$, i.e., 
$f(x) \mid g(x)$. So, $f(x)$ is a prime element of $R[x]$. Since $R[x]$ is noetherian, it follows that 
$R[x]$ is factorial.
\QED

\noindent
{\bf Exercise 1.1.3.}\ \ Show that $k[x]$ is a principal ideal domain (PID) and that a PID is 
factorial.
\svskip

\noindent
{\bf Exercise 1.1.4.}\ \ Let $R$ be an integral domain. Show that if a polynomial ring $R[x]$ is 
factorial so is $R$.

\section{Squeezing-out one variable}

Let $R$ be an integral domain which we assume to be finitely generated over a field $k$ of 
characteristic zero. We are interested in a question asking which conditions make $R$ a polynomial 
ring over $k$. One approach is to squeeze out one variable from $R$ via an algebraic group action. 
Let $\delta$ be a $k$-linear endomorphism of $R$. The we call $\delta$ a {\em $k$-derivation}
\index{$k$-derivation} if $\delta(c)=0$ for $\forall c \in k$ and $\delta(ab)= a\delta(b)+b\delta(a)$ 
for $a, b \in R$. Then we can show that 
\[
\delta^n(ab)=\delta^n(a)b+n\delta^{n-1}(a)\delta(b)+\cdots+\left(\begin{array}{c}n\\ i \end{array}
\right)\delta^{n-i}(a)\delta^i(b)+\cdots+a\delta^n(b).
\]
Furthermore, we say that 
$\delta$ is {\em locally nilpotent}\index{$k$-derivation!locally nilpotent} if for every $a \in R$, 
$\delta^N(a)=0$ for $N \gg 0$. Define a subset $\Ker\delta$ by $\{a \in R \mid \delta(a)=0\}$. 
Then $\Ker\delta$ is a $k$-subalgebra of $R$. It is straightforward to prove the following result 
by making use of the above formula.

\begin{lem}\label{Lemma 1.2.1}
Let $t$ be a variable and define a $k$-linear mapping $\varphi : R \to R[t]$ by 
\[
\varphi(a) =\sum_{n=0}^\infty \frac{1}{n!}\delta^n(a)t^n.
\]
Then $\varphi$ is a $k$-algebra homomorphism.
\end{lem}

In terms of schemes, the homomorphism $\varphi$ gives rise to a morphism $\sigma : G_a\times X \to X$, 
where $G_a=\Spec k[t]$ endowed with a group law and called the {\em additive group scheme}
\index{additive group scheme} and $X=\Spec R$. The group law on $G_a$ is given in such a way that 
for a $k$-algebra $A$ and for elements 
$g, g' \in G_a(A):= {\rm Hom}_{k-alg}(k[t],A)=A$, the product of $g, g'$ under the group law is the 
addition of $g,g'$ in $A$. The morphism $\sigma$ defines a group action on $X$ in the sense that 
for a $k$-algebra $A$ and for $g, g' \in G_a(A)$ and $x \in X(A):={\rm Hom}_{Sch}(\Spec A, X)= 
{\rm Hom}_{k-alg}(R,A)$, we have $\sigma(g,\sigma(g',x))=\sigma(g\cdot g',x)$ and $\sigma(e,x)=x$, 
where $e$ is the identity of the group $G_a(A)$ that is $0$ in $A$. If we identify $g \in G_a(A)$ 
with an element of $A$ and $x \in X(A)$ with a $k$-homomorphism $R \to A$, the element $\sigma(g,x) 
\in X(A)$ corresponds to the $k$-homomorphism $R \to A$ defined by 
$a \mapsto \sum_{n=0}^\infty \frac{1}{n!}x(\delta^n(a))g^n$.

\begin{lem}\label{Lemma 1.2.2}
Let $\delta$ be a non-trivial, locally nilpotent derivation on $R$ and let $R_0=\Ker\delta$. 
Suppose there exists an element $u$ of $R$ such that $\delta(u)=1$. Then $R=R_0[u]$ and $u$ is 
considered to be a variable over $R_0$.
\end{lem}
\Proof
For $a \in R$, define the {\em $\delta$-length $\ell(a)$}\index{$\delta$-length} by 
$\ell(a)=\min\{n \mid \delta^{n+1}(a)=0\}$. It is equal to the degree of a $t$-polynomial 
$\varphi(a)$. We proceed by induction on $\ell(a)$ to show that $a \in R_0[u]$. If $\ell(a)=0$, 
then $a \in R_0$. Let $n=\ell(a)$. Then $\delta^n(a)\ne 0$ but $\delta^{n+1}(a)=0$. Hence 
$\delta^n(a) \in R_0$. Set $a_1=a-\delta^n(a)u^n/n!$. Then $\delta^n(a_1)=0$. Hence $\ell(a_1)
< n$. So, $a_1 \in R_0[u]$ by the hypothesis of induction. Then $a=a_1+(\delta^n(a)/n!)u^n \in 
R_0[u]$. So, $R=R_0[u]$. We have to show that $u$ is a variable over $R_0$, i.e., $u$ is algebraically 
independent over $R_0$. Suppose to the contrary that there exists an algebraic relation
\[
c_0u^n+c_1u^{n-1}+\cdots+c_n=0,
\]
where $c_0, \ldots, c_n \in R_0$ with $c_0c_n \ne 0$. We assume that $n$ is the smallest among all 
such relations. Apply $\delta$ to the left hand side of the equation to get 
\[
nc_0u^{n-1}+(n-1)c_1u^{n-2}+ \cdots+c_{n-1}=0.
\]
Since $c_0\ne 0$, this is a non-trivial relation. Hence this contradicts the choice of $n$. 
\QED

Let $\delta$ be a non-trivial locally nilpotent derivation on $R$. Then there exists an element 
$z$ such that $\delta(z)\ne 0$. Let $n=\ell(z)$. Then $n \ge 1$. So, set $x=\delta^{n-1}(z)$ and 
$a=\delta^n(z)$. Then $\delta(x)=a$ and $a \in R_0$. Denote by $R[a^{-1}]$ the ring of fractions 
of the form $b/a^m$, where $b \in R$ and $m \ge 0$. It is then easy to see that $\delta$ extends 
to a locally nilpotent derivation (denoted by the same letter $\delta$) defined by $\delta(b/a^m)
=\delta(b)/a^m$. Since $\delta(x/a)=1$, Lemma \ref{Lemma 1.2.2} implies that 
$R[a^{-1}]=R_0[a^{-1}][x/a]=R[a^{-1}][x]$. Thus we have shown the following result.

\begin{lem}\label{Lemma 1.2.3}
Let $\delta$ be a non-trivial locally nilpotent derivation on $R$. Then there exist elements 
$a, x$ such that $a \in R_0, \delta(x)=a$ and $R[a^{-1}]=R_0[a^{-1}][x]$.
\end{lem}

An {\em affine $k$-domain}\index{affine $k$-domain} is by definition an integral domain which is 
finitely generated over a field $k$. Let $K$ be the quotient field of $R$. We define 
{\em transcendence degree}\index{transcendence degree} of 
$K$ over $K$ (denoted by $\trdeg_kK$) as follows and then define {\em dimension}\index{dimension} of $R$ 
over $k$ (denoted by $\dim R$) by $\trdeg_kK$. Note that $K$ is finitely generated over $k$ as a field. 
Namely there exists a system of elements $\{\xi_1, \ldots,\xi_n\}$ such that every element of $K$ 
is written as a fraction $f(\xi_1,\ldots,\xi_n)/g(\xi_1,\ldots,\xi_n)$ of two polynomials 
$f(x_1,\ldots,x_n), g(x_1, \ldots,x_n)$ with coefficients in $k$ and with $x_1,\ldots,x_n$ replaced 
by $\xi_1, \ldots,\xi_n$. We apply the following process of taking in and throwing away the 
elements $\xi_1,\ldots,\xi_n$. Consider $\xi_1$. If $\xi_1$ is algebraically independent over $k$, 
then take it in. Otherwise throw it away. Consider next $\xi_2$ over $k(\xi_1)$. If $\xi_2$ is 
algebraically independent over $k(\xi_1)$ then take it is and throw it away otherwise. Continue 
this process. At the $i$-th step, consider $\xi_i$ over $k(\xi_1,\ldots,\xi_{i-1})$ and take it in 
if it is algebraically independent over $k(\xi_1,\ldots,\xi_{i-1})$. Otherwise throw it away. After 
considering all elements $\xi_1,\ldots,\xi_n$, let $\eta_1,\ldots,\eta_d$ be the elements we took in. 
Then $K$ is algebraic over $k(\eta_1,\ldots,\eta_d)$. The above process depends on the choice of 
generators $\{\xi_1, \ldots,\xi_n\}$ and the order of arranging these $n$ elements. But the number 
$d$ of the elements we take in is independent of these choices. So, it is an invariant which is 
proper to $K$. This number $d$ is called {\em transcendence degree}\index{transcendence degree} of $K$ 
over $k$ and denoted by $\trdeg_kK$.

\begin{thm}\label{Theorem 1.2.4}
Let $R$ be an affine domain of $\dim R=1$. If $R$ is given a non-trivial locally nilpotent 
derivation $\delta$, then $R$ is a polynomial ring in one variable over a field $K_0$ which is 
a finite algebraic extension of $k$.
\end{thm}
\Proof
By Lemma \ref{Lemma 1.2.3}, $R[a^{-1}]=R_0[a^{-1}][x]$, where $R_0[a^{-1}]$ is finitely generated over $k$. 
Since $x$ is algebraically independent over $R_0[a^{-1}]$ and since $\dim R=1$, it follows that 
every element of $R_0[a^{-1}]$ is algebraic over $k$. Hence $R_0[a^{-1}]$ is a finite algebraic 
extension of $k$ and coincides with the quotient field $K_0$ of $R_0$.
\QED

\noindent
{\bf Remark.}\ \ It is known \cite{M1} that an affine domain $R$ of $\dim R=1$ defined over 
an algebraically closed field $k$ is a polynomial ring over $k$ if and only if $R$ is factorial 
and $R^*=k^*$.
\svskip

Let $R$ be an affine domain and let $R_0$ be a $k$-subalgebra of $R$ such that $R[a^{-1}]=
R_0[a^{-1}][x]$ with $a \in R_0$ and with $x$ algebraically independent over $R_0$. Define a locally 
nilpotent derivation $\partial$ on $R[a^{-1}]$ by setting $\partial(x)=1$ and $\partial(b)=0$ 
for $\forall b \in R_0$. Write $R=k[u_1,\ldots,u_n]$ with a system of generators $\{u_1,\ldots,u_n\}$
of $R$. Then $a^{m_i}\partial(u_i)\in R$ for $m_i \ge 0$. Let $m=\max\{m_1,\ldots,m_n\}$ and let 
$\delta=a^m\partial$. Then it is easy to see that $\delta$ is a locally nilpotent derivation on $R$ 
with $\delta\mid_{R_0}=0$. We shall show that $\Ker\delta=R_0$ provided $R_0$ satisfies the property 
that $a^rb \in R_0$ with $r > 0$ and $b\in R_0$ implies $b \in R_0$. In fact, we have 
$R_0 \subset \Ker\delta$. Suppose $\delta(b)=0$ for $b \in R$. Then $b \in R_0[a^{-1}]$. Hence 
$a^rb \in R_0$. By the assumed property, we have $b \in R_0$. Hence $\Ker\delta=R_0$.

\begin{thm}\label{Theorem 1.2.5}
Let $R$ be a factorial affine domain defined over an algebraically closed field $k$ of characteristic 
zero. Suppose that there is a non-trivial locally nilpotent derivation $\delta$ on $R$ and that 
$R^*=k^*$. Then $R$ is isomorphic to a polynomial ring in two variables over $k$.
\end{thm}
\Proof
(I)\ \ Let $R_0=\Ker\delta$. Let $a \ne 0$ be an element of $R_0$. Suppose $a=bc$ with $b,c \in R$. 
By applying $\varphi : R \to R[t]$ associated to $\delta$ (cf. Lemma \ref{Lemma 1.2.1}), we have 
$a=\varphi(b)\varphi(c)$. Since $R$ is an integral domain, both $\varphi(b)$ and $\varphi(c)$ are 
constant polynomials in $t$. Namely, $b, c \in R_0$. A $k$-subalgebra $R_0$ of $R$ having this property 
that $a=bc$ with $b, c \in R$ implies $b, c \in R_0$ is said to be {\em inert}\index{subring!inert}. Then 
it is easy to show that an inert $k$-subalgebra of a factorial domain is factorial as well. 
Since $R_0^*=k^*$, $R_0$ is a polynomial ring $k[f]$ by the remark after Theorem \ref{Theorem 1.2.4}. 
Since $R_0$ is inert, it is easy to show that $f-\alpha$ is irreducible in $R$ for every $\alpha \in k$.

(II)\ \ By Lemma \ref{Lemma 1.2.3}, there exist an element $a=\prod_{i=1}^n(f-\alpha_i)$ of $R_0$ and an 
element $g \in R-R_0$ such that $R[a^{-1}]=R_0[a^{-1}][g]$. If $f-\alpha_i$ divides $g-\beta$ for some 
$1 \le i \le n$ and $\beta \in k$, then replace $g$ by $g'=(g-\beta)/(f-\alpha_i)$. If $g'-\beta'$ 
is divisible by some $f-\alpha_i$ for some $\beta'$, we do the same replacement. Hence we may assume 
that $g-\beta$ is not divisible by $f-\alpha_i$ for $\forall \beta\in k$ and $\forall i$. Consider 
an $R_0$-algebra homomorphism $\pi : k[f,g] \to R$. For any $\alpha \in k$, let $\pi_\alpha : 
k[\ol{g}] \to \ol{R}_\alpha$ be the reduction of $\pi$ modulo $(f-\alpha)$. Since $g-\beta$ is not 
divisible by $f-\alpha$, the homomorphism $\pi_\alpha$ is injective. Note that $\ol{R}_\alpha$ 
is an affine domain of dimension one, for $f$ ia algebraically independent over $k$. Under this 
condition, Zariski's main theorem (cf. \cite{H}) says that $\pi$ is an isomorphism. Namely, $R
\cong k[f,g]$. 
\QED

As a corollary of Theorem \ref{Theorem 1.2.5}, we have the following result of Rentschler \cite{R}.

\begin{cor}\label{Corollary 1.2.6}
Let $k$ be an algebraically closed field of characteristic zero and let $k[x,y]$ be a polynomial ring 
in two variables. Let $\delta$ be a locally nilpotent derivation of $k[x,y]$. Then $\delta$ is 
written in the form $\delta=f(x)\partial/\partial y$ after a suitable change of variables.
\end{cor}
\Proof
By Theorem \ref{Theorem 1.2.5}, there are two polynomials $f, g$ such that $\Ker\delta=k[f], 
k[x,y]=k[f,g]$ and $\delta=a\partial/\partial g$, where $a \in k[f]$. Then make the following change of variables $(x,y) \mapsto (f,g)$. Then $\Ker\delta=k[x]$ and $\delta$ has the form 
$f(x)\partial/\partial y$.
\QED

\section{Fibrations of the affine plane by polynomials}

We shall consider mostly a polynomial ring $k[x,y]$ in two variables, where $k$ is an algebraically 
closed field of characteristic zero. We may consider $k$ to be the complex field $\C$ if necessary. Our 
objective is to look into the geometric properties of a given polynomial $f \in k[x,y]$. One way 
is to consider a totality of the polynomials $\{f-\alpha \mid \alpha \in k\}$. 

Hereafter, we use $a,b, c,\cdots$ instead of $\alpha,\beta,\gamma,\cdots$ to denote the elements of 
$k$. For the said purpose, we consider the {\em affine plane}\index{affine plane} $\A^2$ which is 
the set of pairs $\{(a,b)\mid a,b \in k\}$. A {\em curve}\index{curve} on $\A^2$ defined by an equation 
$f=0$ is a subset $\{(a,b)\mid f(a,b)=0\}$, where $f \in k[x,y]$. We use the notations like $V(f), C(f)$ 
to denote this subset. The residue ring $k[x,y]/(f)$ is called the {\em coordinate ring}
\index{coordinate ring} of the curve defined by $f=0$. Sometimes, we consider a topology called the 
{\em Zariski topology}\index{Zariski topology} on $\A^2$ (and the induced topology on $V(f)$ as well) 
such that the open sets are generated by the complements $D(g)$ in $\A^2$ of the curves $V(g)$ with 
$g \in k[x,y]$. If $g=0$, $V(f)$ is not the curve but the total space $\A^2$, while $V(1)$ is the 
empty set. Furthermore, $D(fg)=D(f)\cap D(g)$ and $D(f)\cup D(g)=\A^2-V(\ga)$, where $\ga$ is the ideal 
$(f,g)$ and $V(\ga)$ is the set of points $P$ such that $h(P)=0$ for all $h \in \ga$. For example, 
$D(x)\cap D(y)=\A^2-\{\mbox{$x$-axis}\}\cup\{\mbox{$y$-axis}\}$ and $D(x)\cup D(y)=\A^2-\{(0,0)\}$. 
So, all the points $(a,b)$ are closed points.

A curve $V(f)$ is called {\em irreducible}\index{irreducible} if $f$ is irreducible. Then the residue 
ring $\Gamma(f):=k[x,y]/(f)$ is an integral domain, and one can consider its quotient field which we 
denote by $Q(f)$ or $k(C)$ if the curve is denoted by $C$. We call $Q(f)$ the {\em function field}
\index{function field} of the curve $V(f)$.
\svskip

\noindent
{\bf Example 1.3.1.}\ \ Let $f=y+\lambda(x)$ with $\lambda(x)\in k[x]$. Then $\Gamma(f)\cong k[y], 
\Gamma(f)^*=k^*$ and $Q(f)\cong k(y)$. Any curve whose coordinate ring is isomorphic to a polynomial 
ring in one variable is said to be isomorphic to the {\em affine line}\index{affine line} $\A^1$. 
Hence $V(f)$ for the above $f$ is isomorphic to $\A^1$. Any irreducible curve $V(f)$ is called 
{\em rational}\index{curve!rational} if $Q(f)$ is generated by one element over the field $k$.

A curve isomorphic to $\A^1$ is classified by the following result of Abhyankar-Moh-Suzuki \cite{AM,S}.

\begin{thm}\label{Theorem 1.3.1}
Let $f \in k[x,y]$ satisfy the condition that $V(f)\cong \A^1$. Then there exists an automorphism 
$\sigma$ of $k[x,y]$ such that $\sigma(f)=x$. Hence the curves $V(f-a)$ are isomorphic to $\A^1$ 
for all $a \in k$.
\end{thm}

We can consider a mapping of the sets $\A^2 \to \A^1$ defined by $(a,b)\mapsto f(a,b)$. Then the 
inverse image (called the {\em fiber}\index{fiber}) of a point $c$ of $\A^1$ is the curve $V(f-c)$. 
We denote this mapping by $\Lambda(f)$. It is called a {\em linear pencil}\index{linear pencil} 
defined by $f$. In the case of Theorem \ref{Theorem 1.3.1}, all the fibers of $\Lambda(f)$ are 
isomorphic to $\A^1$. 
\svskip

\noindent
{\bf Example 1.3.2.}\ \ A rational curve which is not isomorphic to $\A^1$ is exemplified by $V(f)$ 
with $f=xy-1$. Then $\Gamma(f)\cong k[x,x^{-1}]$. Hence $\Gamma(f)^*=k^*\times \Z$, where $\Z\cong 
\{x^n\mid m\in \Z\}$, and $Q(f)\cong k(x)$. So, $V(xy-1)$ is not isomorphic to $\A^1$. A curve 
isomorphic to $V(xy-1)$ is called an {\em algebraic torus}\index{algebraic torus} of dimension $1$ 
and denoted by $\A^1_*$.
\svskip

\noindent
{\bf Example 1.3.3.}\ \ It is known that a curve defined by $f_a=y^2-x^3-a$ with $a \in k^*$ is not 
rational. Let $C$ be a curve defined by $f=y^2-x^3$ and let $\Lambda(f) : \A^2 \to \A^1$ be a 
linear pencil defined by $f=y^2-x^3$. Then the fiber $V(f_a)$ with $a\ne 0$ is not rational, 
while the fiber $V(f_0)$ is rational.
\svskip

\noindent
{\bf Exercise 1.3.1.}\ \ Show that the curve $V(f)$ with $f=y^2-x^3-1$ is not rational. In fact, 
if it were rational, there would exist polynomials $g(t), h(t), \ell(t) \in k[t]$ such that 
$\gcd(g(t),h(t),\ell(t))=1$ and $\ell(t)h(t)^2=g(t)^3+\ell(t)^3$. Show that this leads to a 
contradiction.
\svskip

In general, if the total degree of $f$ becomes bigger, then the curve $V(f)$ becomes far from being 
rational. In 1.2, we observed that a non-trivial locally nilpotent derivation $\delta$ on an affine 
$k$-domain gives a variable $x$ such that $R[a^{-1}]=R_0[a^{-1}][x]$, where $R_0=\Ker \delta$ and 
$a \in R_0$. The subalgebra $R_0$ needs not to be an affine $k$-domain as seen by many counterexamples 
to the Fourteenth Problem of Hilbert \cite{N, Ro, KM, F, DF}. \footnote{Kuroda of T\^ohoku University 
has recently shown that a counterexample exists even in case $\dim R_0=3$.} It is known by Zariski 
\cite{N2} that $R_0$ is finitely generated if $\dim R_0 \le 2$. Suppose that $R_0$ is an affine domain. 
We denote by $\Spec R$ an affine variety (or an affine scheme, in general) with coordinate ring $R$. 
Then the inclusion $R_0 \hookrightarrow R$ defines a morphism $f : X \to B$, where $X=\Spec R$ and 
$B=\Spec R_0$. Let $U$ be an open set of $B$ such that $a \ne 0$ on $U$, i.e., $U=\Spec R_0[a^{-1}]$.
Then $f^{-1}(U)\cong U\times \A^1$. So, the fibers $f^{-1}(P)$ over a point $P \in U$ is isomorphic 
to $\A^1$. We define, in general, an {\em $\A^1$-fibration}\index{$\A^1$-fibration} to be a morphism 
of algebraic varieties $f : X \to B$ such that $f^{-1}(U)\cong U \times \A^1$ for an non-empty 
open set $U$ of $B$. Hence the morphism $f : \Spec R \to \Spec R_0$ is an $\A^1$-fibration, 
which we call the {\em quotient morphism}\index{quotient morphism} under $\delta$ or the associated 
$G_a$-action.

Conversely, if given an $\A^1$-fibration $f : X \to B$ with affine varieties $X=\Spec R$ and 
$B=\Spec R_0$, we find an open set $U=D(a) \subset U$ such that $R[a^{-1}]=R_0[a^{-1}][x]$. 
As explained after Theorem \ref{Theorem 1.2.4}, there exists a locally nilpotent derivation $\delta$ 
on $R$ such that $R_0=\Ker \delta$ and $f$ coincides with the quotient morphism under $\delta$. 
This is the case only if $B$ is an affine variety. If $B$ is not affine, $f$ does not come from 
a $G_a$-action.

For an affine variety $X=\Spec R$, we define the {\em Makar-Limanov invariant}
\index{Makar-Limanov invariant} $\ML(X)$ as the intersection of $\Ker \delta$ in the coordinate ring $R$ 
of $X$, where $\delta$ ranges over all locally nilpotent derivations of $R$. Hence $\ML(X)$ is 
a subalgebra of $R$. If $\dim X=1$ and there is a non-trivial locally nilpotent derivation $\delta$, 
then $\Ker\delta=k$. Hence $\ML(X)=k$. When $\dim X=2$, we have the following result. 

\begin{lem}\label{Lemma 1.3.2}
Let $X=\Spec R$ be an affine surface. Then one of the following cases takes place.
\begin{enumerate}
\item[{\rm (1)}]
$\ML(X)=R$ and there are no locally nilpotent derivations on $R$.
\item[{\rm (2)}]
$\ML(X)=\Ker\delta$ and there is a nontrivial locally nilpotent derivation $\delta$. If $\delta'$ is 
another locally nilpotent derivation on $R$, then $\delta'$ is {\em conjugate}
\index{$k$-derivation!conjugate} to $\delta$, which signifies, by definition, that 
$\Ker\delta'=\Ker\delta'$ and $c\delta=d\delta'$ for some $c, d \in \Ker\delta$.
\item[{\rm (3)}]
$\ML(X)=k$ and there are two {\em algebraically independent}
\index{$k$-derivation!algebraically independent} 
locally nilpotent derivations $\delta$ and $\delta'$, which signifies, by definition, that there 
are elements $\xi \in \Ker\delta$ (resp. $\xi'\in\Ker\delta'$) such that $\xi$ (resp. $\xi'$) is 
algebraically independent over $\Ker\delta'$ (resp. $\Ker\delta$).
\end{enumerate}
\end{lem}
\Proof
(1)\ \ If there are no locally nilpotent derivations on $R$, then $\ML(X)=R$. 

(2)\ \ If $\delta$ is a nontrivial locally nilpotent derivation, then $R[a^{-1}]=\\
R_0[a^{-1}][x]$ for $a \in R_0=\Ker\delta$ and $x \in R$. Let $\delta'$ be another locally nilpotent 
derivation on $R$. If $\delta'$ is conjugate to $\delta$ ( $\delta \sim \delta'$ by notation), then 
it is obvious that $\Ker\delta'=R_0$. Suppose $\dim\ML(X)=1$. Then $R_0 \supseteq \ML(X)$ and $R_0$ is 
algebraic over $\ML(X)$. Since $\Ker\delta' \supseteq \ML(X)$ and $R_0$ is algebraic over $\ML(X)$, it 
follows that $\Ker\delta' \supseteq R_0$. This implies that $\Ker\delta'=R_0=\ML(X)$. Considering 
the derivation $\delta'$, we can write $R[b^{-1}]=R_0[b^{-1}][y]$ with $b \in R_0$ and $y \in R$. 
In particular, $R\otimes_{R_0}K_0=K_0[x]=K_0[y]$, where $K_0=Q(R_0)$. So, $y=\alpha x+\beta$ with 
$\alpha, \beta \in K_0$. Hence we may assume $y=x$ by replacing $b$ with a suitable element of $R_0$. 
Since $\delta=d\partial/\partial x$ and $\delta'=c\partial/\partial x$ with $c, d \in R_0$, we have 
$c\delta=d\delta'$. Namely $\delta \sim \delta'$. 

(3)\ \ Suppose $\dim\ML(X)=0$. Then there exist locally nilpotent derivations $\delta, \delta'$ of $R$
such that $\delta'$ is nontrivial on $\Ker\delta$. Hence there exists an element $\xi \in \Ker\delta$ 
such that $\delta'(\xi) \ne 0$. If one writes $R[b^{-1}]=\Ker\delta'[b^{-1}][y]$ with 
$b \in \Ker\delta'$ and $y \in R$, the element $\xi$ is a polynomial in $y$ of positive degree 
with coefficients in $\Ker\delta'[b^{-1}]$. Since $y$ is algebraically independent over 
$\Ker\delta'$, so is the element $\xi$. Since $\delta$ is nontrivial on $\Ker\delta'$, we find similarly 
an element $\xi' \in \Ker\delta'$ which is algebraically independent over $\Ker\delta$. This implies 
that $\dim(\Ker\delta\cap\Ker\delta')=0$ and $\Ker\delta\cap\Ker\delta'$ is algebraic over $k$. So, 
$\Ker\delta\cap\Ker\delta'=\ML(X)=k$.
\QED

\section{Plane-like affine surfaces}

We say that a morphism $f : X \to B$ of algebraic varieties is an $G$-{\em fibration}\index{fibration} 
if almost all fibers $f^{-1}(Q)$ with $Q \in B$ are isomorphic to one and the same algebraic variety 
$G$. We mostly consider the case where $\dim B=1$ and $G$ is isomorphic to $\A^1, \A^1_*$ and $\BP^1$. 
Any fiber $F$ which is not isomorphic to $G$ is called a {\em singular}\index{singular fiber} fiber of $f$.

By Theorem \ref{Theorem 1.2.5}, the affine plane $\A^2$ is an affine surface such that its 
coordinate ring $R$ is factorial, $R^*=k^*$ and there exists an $\A^1$-fibration $\rho : X \to \A^1$. 
In the following, we consider an affine surface $X = \Spec R$ with an $\A^1$-fibration $\rho : 
X \to C:=\A^1$ such that $R^*=k^*$. In order to look into the properties of $X$, we need to consider 
a smooth projective surface $V$ which contains $X$ as an open set and has a $\BP^1$-fibration 
$p : V \to \ol{C}\cong \BP^1$ extending the $\A^1$-fibration $\rho$. This implies that $p\mid_X= 
\rho$. Such a surface $V$ is called a {\em smooth compactification}\index{smooth compactification} 
or a {\em smooth completion}\index{smooth completion} of $X$.

The existence of $V$ and $p$ is shown as follows. Since $X$ is a closed subvariety of a certain 
affine space $\A^N$, we take the closure $\ol{X}$ of $X$ in the projective space $\BP^N$ to which 
$\A^N$ is naturally embedded as the complement of a hyperplane. Though $\ol{X}$ might have singular 
points, we can resolve singularities of $\ol{X}$ which are centered on $\ol{X}\setminus X$. Thus there 
exists a smooth projective surface $V'$ containing $X$ as an open set. Then $\rho : X \to C$ extends 
to a rational mapping $p' : V' \to \ol{C}$, where $\ol{C}$ is a {\em smooth completion} of $C$, 
which means that $\ol{C}$ is a smooth projective curve containing $C$ as an open set. In our 
case $\ol{C}$ is isomorphic to $\BP^1$. If there exists points of $V'$ where $p'$ is not defined, we 
can eliminate indeterminacy by applying a sequence of blowing-ups $\sigma : V \to V'$ with centers 
outside $X$. Namely, the rational mapping $p:=p'\cdot\sigma : V \to \ol{C}$ is a morphism extending 
$\rho$. Note that there a general fiber of $\rho$ is isomorphic to $\A^1$. Hence a general fiber of 
$p$ is a smooth projective curve containing $\A^1$ as an open set. Hence the general fiber of $p$ is isomorphic to $\BP^1$. Namely $p$ is a $\BP^1$-fibration. 

Here is another reasoning of finding an open embedding $(X,\rho)\hookrightarrow (V,p)$ with $\rho=p|_X$. 
Note that $\rho^{-1}(U)\cong U\times \A^1$ for an open set $U$ of $C$, where $\rho|_U$ is identified 
with the first projection of $U\times \A^1$ onto $U$. Since $\A^1$ is embedded into $\BP^1$ as the 
complement of the point at infinity, $U\times \A^1$ is embedded as an open set into 
$V_0:=\ol{C}\times \BP^1$. Hence there exists a birational mapping $\phi_0 : X \dlto V_0$ such that 
$p_0\cdot\phi_0=\rho$, where $p_0$ denotes the first projection $\ol{C}\times \BP^1 \to \ol{C}$. 
If there are points on $V_0$ (resp. $X$) corresponding to irreducible curves on $X$ (resp. $V_0$) 
under the above birational mapping, we blow up the surfaces $X$ and $V_0$ to obtain $X_1$ and $V_1$ 
so that there are no points on $X_1$ (resp. $V_1$) corresponding to irreducible curves on $V_1$ 
(resp. $X_1$). Then $X_1$ is identified with an open set of $V_1$ by virtue of Zariski's Main Theorem. 
Let $\sigma : X_1 \to X$ and $\tau : V_1 \to V_0$ be the sequences of blowing-ups. Since the 
exceptional curves in $X_1$ arising from $\sigma$ is contained in $V_1$ as well, we may contract these
exceptional curves in $V_1$ to obtain a smooth projective surface $V$. Then $X$ is identified with 
an open set of $V$. It is then clear that there is a $\BP^1$-fibration $p : V \to \ol{C}$ such that 
$p|_X=\rho$.

We shall recall some of the basic properties of a $\BP^1$-fibration on a smooth projective surface.

\begin{lem}\label{Lemma 1.4.1}
Let $p: V \to B$ be a $\BP^1$-fibration on a smooth projective surface $V$ with base a smooth curve $B$. 
Let $F$ be a singular fiber of $p$. Then the following assertions hold.
\begin{enumerate}
\item[{\rm (1)}] 
Every irreducible component of $F$ is isomorphic to $\BP^1$. If two irreducible components meets 
each other, they meet transversally. Furthermore, no three components meet in one point. Hence 
$F$ is a tree of smooth rational curves.
\item[{\rm (2)}] 
$F$ contains a $(-1)$-curve, {\em called a $(-1)$-component}\index{$(-1)$-component}, and 
any $(-1)$-curve in $F$ meets at most two other components of $F$. If a $(-1)$-curve $E$ occurs 
with multiplicity $1$ in the scheme-theoretic fiber $F$ then $F$ contains another $(-1)$-curve.
\item[{\rm (3)}] 
By successively contracting $(-1)$-curves in $F$ and its images we can reduce $F$ to a regular fiber.
\end{enumerate}
\end{lem}
For the proof, see Gizatullin \cite{Giz}.

A smooth projective surface $V$ with a $\BP^1$-fibration is called a {\em minimal rational ruled surface}
\index{minimal rational ruled surface} or simply a {\em Hirzebruch surface}\index{Hirzebruch surface} 
if every fiber is irreducible and the base curve is rational. Minimal rational ruled surfaces are 
classified in the following result.

\begin{lem}\label{Lemma 1.4.2}
The following assertions hold true.
\begin{enumerate}
\item[{\rm (1)}]
Every minimal rational ruled surface is isomorphic to $\Sigma_n:=\Proj(\SO_{\BP^1}\oplus\SO_{\BP^1}(n))$,
where $n \ge 0$. 
\item[{\rm (2)}]
There is a section $M$ {\em (called minimal section)}\index{minimal section} on $\Sigma_n$ such that 
$\sis{M}=-n$. If $n > 0$ then $M$ is unique. For other sections $S$, $\sis{S}\ge n$. Hence 
$S \sim M+a\ell$ with $a \ge n$.
\end{enumerate}
\end{lem}

A smooth projective surface $V$ with a $\BP^1$-fibration is obtained from a Hirzebruch surface 
$\Sigma_n$ by applying several blowing-ups. A singular fiber, which is then a reducible fiber of 
the $\BP^1$-fibration, is described in Lemma \ref{Lemma 1.4.1}. Write such a singular 
fiber as $F=\sum_{i=1}^ra_iF_i$, where the $F_i$ are irreducible components and the $a_i$ are their 
coefficients (hence $a_i > 0$). As explained above, any $\BP^1$-fibration has a cross-section, 
say $S$, which comes from the (a) minimal section of a Hirzebruch surface and for which we have 
$\is{F}{S}=1$. So, one of the $a_i$ must be equal to $1$. Hence any $\BP^1$-fibration has no multiple 
fibers. 

As a consequence of the above observations, we have the following result which describes the singular 
fibers of an $\A^1$-fibration on a smooth affine surface. Let $F=\sum_{i=1}^ra_iF_i$ denote anew a 
fiber of such $\A^1$-fibration. Then $F$ is {\em singular}\index{singular fiber} if either $r \ge 2$, 
or $r=1$ and $a_1 \ge 2$. If $F$ is singular, we call $m:=\gcd(a_1,\ldots,a_r)$ the {\em multiplicity}
\index{multiplicity} of $F$. We have the following result.

\begin{lem}\label{Lemma 1.4.3}
Let $X$ be a smooth affine surface with an $\A^1$-fibration $\rho : X \to C$, where $C \cong \A^1$ or 
$C \cong \BP^1$. Then the following assertions hold.
\begin{enumerate}
\item[{\rm (1)}]
Any fiber of $\rho$ is a disjoint union of the affine lines with suitable multiplicities.
\item[{\rm (2)}]
If $C \cong \A^1$, then $\rank\Pic(X)=\sum_{P\in C}(r_P-1)$, where $r_P$ is the number of 
irreducible components of the fiber $\rho^{-1}(P)$. If $C \cong \BP^1$, then $\rank\Pic(X)
=1+\sum_{P\in C}(r_P-1)$.
\item[{\rm (3)}]
Suppose that $C\cong \A^1$ and $r_P=1$ for every $P \in C$. Then $\Pic(X) \cong \prod_{P\in C}\Z/m_P\Z$,
where $m_P$ is the multiplicity of the fiber $\rho^{-1}(P)$.
\item[{\rm (4)}]
With the same assumptions as in the assertion (3), $X$ is a $\Q$-homology plane. Namely, $X$ is a 
smooth affine surface with $H_i(X;\Q)=(0)$ for every $i > 0$.
\end{enumerate}
\end{lem}
\Proof
(1)\ \ Embed $X$ into a smooth projective surface $V$ with a $\BP^1$-fibration $p : V \to \ol{C}$. 
Let $D:=V-X$, which we consider as a reduced effective divisor. Since $X$ is affine, there is a very 
ample divisor whose support is equal to $D$. In particular, $D$ is connected (cf. \cite{H}). Since 
$D$ contains a (unique) cross-section, say $S$, we know that $D$ consists of $S$ and some connected 
parts of the fibers of $\rho$ which meet $S$. By Lemma \ref{Lemma 1.4.1}, the fibers of $\rho$ are 
trees of rational curves. Since $X$ contains no complete curves, it follows that any fiber of $\rho$ 
with the part in $D$ removed off consists of irreducible components which are located at the tips of 
the tree (the ends of tree branches) and each of which is deprived of one point where it meets the part 
in $D$. This means that any fiber of $\rho$ is a disjoint union of the affine lines. 

\noindent
(2)\ \ Note that $\rank\Pic(\Sigma_n)=2$ and each blowing-up increases the rank of the Picard group 
by one. Then $\rank\Pic(V)=2+\sum_{P\in\ol{C}}(s_P-1)$, where $s_P$ is the number of irreducible 
components of the fiber $p^{-1}(P)$. Since $D=S+\sum_{P\in\ol{C}}(D\cap p^{-1}(P))$ and 
$D\cap p^{-1}(P)$ consists of $(s_P-r_P)$ components, we obtain the asserted result.

\noindent
(3)\ \ If $C \cong \A^1$, one point of $\ol{C}$, say $P_\infty$, is not included in $C$. Let 
$\ol{F}_\infty$ be the fiber of $p$ over $P_\infty$. Then every fiber $p^{-1}(P)$ is linearly 
equivalent to $\ol{F}_\infty$. Hence if $m_PF_P$ is the component of $p^{-1}(P)$ left in $X$ with 
multiplicity $m_P$, then we have a relation $m_PF_P \sim 0$. In fact, there are no other relations 
on the $F_P$. Hence we have the asserted result.

\noindent
(4)\ \ We just indicate the key ingredients of the proof without going too much to details. First of 
all, $V$ is a smooth rational surface, whence $H^1(V,\SO_V)=H^2(V,\SO_V)=(0)$. Then $H^1(V;\Q)=(0)$ 
because the first Betti number $b_1$ equals to $2h^{0,1}$, where $h^{0,1}=\dim H^1(V,\SO_V)=0$. 
From the well-known long exact sequence 
\[
\cdots \to H^1(V,\SO_V) \to H^1(V,\SO_V^*) \to H^2(V;\Z) \to H^2(V,\SO_V) \to \cdots ,
\]
it follows that $\Pic(V)\otimes\Q\cong H^2(V;\Q)$, where $\Pic(V)\cong H^1(V,\SO_V)^*$. Since the 
irreducible components of $D$ are independent in $H^2(V;\Q)$ and $\Pic(X)\otimes\Q=(0)$, it follows 
that $H^2(V;\Q)\cong H^2(D;\Q)$. Now consider the cohomology exact sequence of the pair $(V,D)$ with 
$\Q$-coefficients
\begin{eqnarray*}
\cdots &\to & H^1(V) \to H^1(D) \to H^2(V,D) \\
 &\to & H^2(V) \to H^2(D) \to H^3(V,D) \to H^3(V),
\end{eqnarray*}
whence we obtain $H^1(D;\Q)\cong H^2(V,D;\Q)\cong H_2(X;\Q)$ and $H_1(X;\Q)\cong H^3(V,D;\Q)=(0)$. 
Since $D$ is a tree of smooth rational curves, we know that $H^1(D;\Q)=(0)$. Hence $H_2(X;\Q)=(0)$. 
Since $\dim X=2$ it is clear that $H_i(X;\Q)=(0)$ if $i > 2$.
\QED

\begin{cor}\label{Corollary 1.4.4}
Let $X$ be as in Lemma \ref{Lemma 1.4.3}. Assume that $C \cong \A^1$ and that every fiber of 
$\rho$ is irreducible. Then the Picard group of $X$ is a torsion group. In particular, if there 
is only one multiple fiber $mF$, then $\Pic(X)$ is a cyclic group of order $m$.
\end{cor}

There are many kinds of affine surfaces satisfying the same conditions as in the above corollary. 
We define below just one class of such surfaces which include the surfaces $\BP^2-C$, where $C$ is 
a curve $X_0X_1^{d-1}=X_2^d$ with $d \ge 2$ (cf. Lemma \ref{Lemma 1.4.6} below). The remaining part 
of this subsection will be published elsewhere in our forthcoming paper \cite{KM3}.
 
\begin{defn}\label{Definition 1.4.5}
{\rm A smooth affine surface $X$ is an {\em affine pseudo-plane}\index{affine pseudo-plane} if $X$ 
satisfies the following conditions:
\begin{enumerate}
\item[{\rm (1)}]
$X$ has an $\A^1$-fibration $\rho : X \to A$, where $A \cong \A^1$.
\item[{\rm (2)}]
The $\A^1$-fibration $\rho$ has a unique multiple fiber with multiplicity $d \ge 2$.
\end{enumerate}
We say that $X$ has {\em type $(d,n,r)$}\index{affine pseudo-plane!type $(d,n,r)$} if $X$ 
further satisfies the next condition:
\begin{enumerate}
\item[{\rm (3)}]
$X$ has a smooth completion $(V,D)$ such that the dual graph of $D$ is as given in Figure 1 below, 
where $n \ge 1$ and $r \ge 2$. Furthermore, $\ol{F}$ is the closure of $F$ in $V$ and $S'$ is the 
unique cross-section contained in $D$.

\raisebox{-45mm}{
\begin{picture}(140,45)(30,0)
\unitlength=0.95mm
\put(10,20){\circle{1.8}}
\multiput(10,22)(0,1.5){9}{\circle*{0.2}}
\put(10,35){\circle{1.8}}
\put(11,20){\line(1,0){9}}
\multiput(21,20)(1.5,0){5}{\circle*{0.2}}
\put(30,20){\line(1,0){9}}
\put(40,20){\circle{1.8}}
\put(41,20){\line(1,0){13}}
\put(55,20){\circle{1.8}}
\put(55,21){\line(0,1){18}}
\put(55,40){\circle{1.8}}
\put(56,20){\line(1,0){13}}
\put(70,20){\circle{1.8}}
\put(71,20){\line(1,0){9}}
\multiput(81,20)(1.5,0){5}{\circle*{0.2}}
\put(90,20){\line(1,0){9}}
\put(100,20){\circle{1.8}}
\put(101,20){\line(1,0){13}}
\put(115,20){\circle{1.8}}
\put(116,20){\line(1,0){13}}
\put(130,20){\circle{1.8}}
\put(130,21){\line(0,1){13}}
\put(130,35){\circle{1.8}}
\put(2,20){$-2$}
\put(2,35){$-1$}
\put(8,13){$E_{d+r-1}$}
\put(14,35){$\ol{F}$}
\put(36,25){$-2$}
\put(47,25){$-2$}
\put(36,13){$E_{d+1}$}
\put(53,13){$E_d$}
\put(47,40){$-d$}
\put(59,40){$E_1$}
\put(68,13){$E_{d-1}$}
\put(66,25){$-2$}
\put(98,13){$E_2$}
\put(98,25){$-2$}
\put(113,13){$\ell'_0$}
\put(113,25){$-2$}
\put(128,13){$S'$}
\put(125,35){$0$}
\put(134,35){$\ell'_\infty$}
\put(134,20){$-n$}
\put(68,3){{\rm Figure 1}}
\end{picture}}
\end{enumerate}}
\end{defn}

\noindent
If the curves $\ol{F}=E_{d+r}, E_{d+r-1}, \ldots, E_{d+1}, E_d, E_{d-1}, \ldots, E_2, E_1$ are 
contracted in this order, the resulting surface is a Hirzebruch surface $\Sigma_n$ with the 
minimal cross-section $S'$.
\svskip

Choose a point $P$ on the fiber $\ell'_\infty$. Blow up the point $P$ to obtain a $(-1)$ curve 
$E$. Then the proper transform $L$ of $\ell'_\infty$ is a $(-1)$ curve. Contract $L$ to obtain 
the same figure as before with $\ell'_\infty$ replaced by the image of $E$ and with $\sis{S'}$ 
is $-n+1$ if $P \ne S'\cap\ell'_\infty$ and $-n-1$ if $P=S'\cap\ell'_\infty$. This operation is 
called the {\em elementary transformation}\index{elementary transformation} with center $P$. 
After several elementary transformations, we may and shall assume unless otherwise specified that 
$n=-\sis{S'}=1$. So, we say that an affine pseudo-plane has type $(d,r)$ instead of type $(d,1,r)$.

\begin{lem}\label{Lemma 1.4.6}
Let $X$ be an affine pseudo-plane of type $(d,r)$. Then $X$ is isomorphic to the complement of 
$M_0\cup C_d$ if $r < d$ and $M_1\cup C_d$ if $r \ge d$ in the Hirzebruch surface $\Sigma_n$ 
with $n=|r-d|$, where $M_0$ is the minimal section and where $C_d$ and $M_1$ are specified as follows. 
In the case $r < d$, $C_d$ is an irreducible member of the linear system $|M_0+d\ell_0|$ which meets 
$M_0$ in the point $M_0\cap \ell_0$ with multiplicity $r$. In the case $r \ge d$, $M_1$ is a section 
of $\Sigma_n$ with $\sis{M_1}=n$, and $C_d$ is an irreducible member of the linear system 
$|M_1+d\ell_0|$ which meets $M_1$ in the point $M_1\cap \ell_0$ with multiplicity $r$. In both cases, 
$\ell_0\cap X=\ol{F}\cap X$ for the fiber $\ell_0$ of $\Sigma_n$.
\end{lem}
\Proof
Contract $S', \ell'_0, E_2, \ldots, E_d, E_{d+1}, \ldots, E_{d+r-1}$ in this order. Then the resulting 
surface is the Hirzebruch surface $\Sigma_n$ with $n=|r-d|$ and the image of $\ell'_\infty$ provides 
$C_d$. The image of $E_1$ provides $M_0$ or $M_1$ according as $r-d < 0$ or $r-d \ge 0$, while 
the image of $\ol{F}$ is the fiber $\ell_0$.
\QED

Lemma \ref{Lemma 1.4.6} gives a construction of affine pseudo-planes from the Hirzebruch 
surfaces. Whenever we are conscious of this construction from the Hirzebruch surfaces, we denote the 
affine pseudo-plane $X$ in Lemma \ref{Lemma 1.4.6} by $X(d,r)$.

\begin{lem}\label{Lemma 1.4.7}
Let $X$ be an affine pseudo-plane. Then $X$ is isomorphic to $\BP^2-C$, where $C$ is a curve on 
$\BP^2$ defined by $X_0X_1^{d-1}=X_2^d$ with $d \ge 2$, if and only if $X$ has type $(d,d-1)$. 
\end{lem}
\Proof
It is straightforward to see that $\BP^2-C$ is an affine pseudo-plane of type $(d,d-1)$. 
Conversely, if $X$ is an affine pseudo-plane of type $(d,d-1)$, then contract in the Figure 1 of 
Definition \ref{Definition 1.4.5} the curves $S', \ell'_0, E_2, \ldots, E_d,\\ 
E_{d+1}, \ldots, E_{2d-2}, E_1$ in this order to obtain the plane curve $X_0X_1^{d-1}=X_2^d$ which 
is the image of $\ell'_\infty$.
\QED

Let $X$ be a smooth (complex) affine surface. A smooth projective surface $V$ is a {\em normal 
compactification}\index{normal compactification} of $X$ if $V\setminus X$ consists of smooth 
irreducible curves such that they meet each other transversally and that there are no three or more 
curves meeting in one point. Further, it is a {\em minimal}\index{normal compactification!minimal} 
normal compactification if the contraction of any possible $(-1)$ curve in $V\setminus X$ breaks the 
above normality condition. We consider $V\setminus X$ the {\em boundary divisor}\index{boundary divisor} 
and denote it simply by $D:=V-X$. So, $D$ is a reduced effective divisor. 

Suppose that the divisor $D$ is a tree of rational curves. Extracting the definition from \cite{Mumford} 
and \cite{Ra}, we shall define a group associated to a weighted graph. Let $\Gamma$ be a weighted 
connected tree which consists of vertices with weights and edges, each of edges connecting exactly two 
vertices. The group $\pi_1(\Gamma)$ associated to 
$\Gamma$ is generated by as many generators as the vertices $\{v_n\}$ of $\Gamma$ which are subject 
to the following relations. First, order the vertices in some fixed manner as $v_1,v_2,\ldots$. 
For each vertex $v$ with weight $d$, let $v_{i_1},v_{i_2},\ldots,v_{i_r}$ be all the vertices 
connected to $v$ such that $i_1<i_2,\ldots <i_r$ (then $v$ occurs among these vertices, say $v=v_{i_a}$). 
For any vertex $v$, we consider a relation $v_{i_1}\cdots v_{i_{a-1}}\cdot v_{i_a}^d\cdots v_{i_r}=e$. 
Next, if $v,w$ are any two vertices which are connected by an edge, then we impose the relation 
$v\cdot w=w\cdot v$ (i.e., the generators corresponding to $v,w$ commute). The quotient of the free group 
generated by the vertices of $\Gamma$ modulo these relations is the {\em fundamental group at infinity}
\index{fundamental group at infinity} $\pi_{1,\infty}(X)$. It is known that $\pi_{1,\infty}(X)$ is 
independent of the choice of a minimal normal compactification. This group modulo the commutator 
subgroup is the {\em first homology group at infinity}\index{first homology group at infinity} 
$H_{1,\infty}(X)$.

By definition, a {\em $\Q$-homology plane}\index{$\Q$-homology plane} $Y$ is a smooth affine surface 
such that $H_i(Y;\Q)=(0)$ for $i>0$. It is well-known that the divisor at infinity of a $\Q$-homology 
plane $Y$ in a normal compactification is a (connected) tree of rational curves (see \cite{MS}). 
Hence we can use Mumford-Ramanujam's presentation of $\pi_{1,\infty}(Y)$. An {\em appendix}
\index{appendix} of a weighted graph $\Gamma$ is a pair of vertices $\{u,v\}$ such that $u$ is 
a vertex of $\Gamma$ with weight $0$ which is connected only to the vertex $v$ and $v$ is a vertex 
linked to $u$ and at most one other vertex of $\Gamma$. {\em Removing}\index{appendix!removing} 
an appendix $\{u,v\}$ from $\Gamma$ is an operation of removing vertices $u,v$, the edge connecting 
$u,v$ and the edge which connects $v$ to a third vertex, if it exists.

\begin{lem}\label{Lemma 1.4.8}
Let $\Gamma$ be a weighted connected tree. Then the
following assertions hold.  
\begin{enumerate} 
\item[{\rm (1)}] 
Let $\{u,v\}$ be an appendix of $\Gamma$ such that the weight of $v$ is $d$ and let $\Gamma'$ be 
obtained from $\Gamma$ by removing the appendix $\{u,v\}$. Then $\pi_1(\Gamma)=\pi_1(\Gamma')$.  
\item[{\rm (2)}] 
Suppose that $\Gamma$ is a linear chain. Then $\pi_1(\Gamma)$ is a cyclic group. If further the 
determinant of the intersection form on the free abelian group on the vertices of $\Gamma$ is non-zero 
(in particular, if this form is negative definite) then this group is finite cyclic. In the remaining 
case, $\pi_1(\Gamma)$ is an infinite cyclic group.  
\end{enumerate} 
\end{lem}
\Proof 
(1)\ Let $\{u,v\}$ be an appendix of $\Gamma$ with the weight of $v$ being $d$. If $\Gamma$ consists 
only of vertices $u, v$ and an edge connecting them, it is easy to see that $\pi_1(\Gamma)=e$. Otherwise, 
let $w$ be a third vertex connected to $v$. Let the ordering on the vertices of $\Gamma$ be $u,v,w,\ldots$. 
Since $u^0v=e$ and $uv^dw=e$, we have $v=e$ and $u=w^{-1}$. Hence it follows that 
$\pi_1(\Gamma)=\pi_1(\Gamma')$. 

\noindent
(2)\ Let the ordering on $\Gamma$ be $v_1,v_2,\ldots$. From Mumford-Ramanujam's presentation we see 
easily that the generator corresponding to $v_1$ generates $\pi_1(\Gamma)$.  If the intersection form 
is negative definite, then the order of this group is the absolute value of the determinant of the 
intersection form, hence finite. Suppose that the form has at least one positive eigenvalue. It is 
proved in \cite[Lemma 5]{GS}, that there is a connected linear weighted tree $\Gamma_1$ with the 
following properties: 
\begin{enumerate} 
\item[(a)] $\Gamma_1$ has a connected subtree $T$ with vertices $w_1,w_2,\ldots,w_{2n}$ such that 
weights of $w_1,w_2,\ldots,w_{2n}$ are all $0$.  
\item[(b)] The subgraph of $\Gamma_1$ obtained by removing $T$, say $\Gamma_1-T$, is connected 
(it may be empty) and either consists of a single vertex of weight $0$ or has negative definite 
intersection form.  
\item[(c)]
$\pi_1(\Gamma)\cong\pi_1(\Gamma_1)\cong\pi_1(\Gamma_1-T)$. 
\end{enumerate} 
Now if $\Gamma_1=T$ then by the argument in part (1) we see that $\pi_1(\Gamma)=(e)$. If $\Gamma_1-T$ 
is non-empty, then we see by Mumford-Ramanujam's presentation above that $\pi_1(\Gamma_1-T)$ is a 
cyclic group. It is finite if the intersection form on $\Gamma_1-T$ is negative definite.  
\QED

\begin{lem}\label{Lemma 1.4.9}
Let $X$ be an affine pseudo-plane of type $(d,r)$. Then $\pi_{1,\infty}(X)$ is a group generated by 
$x,y$ with relations $x^r=y^d=(xy)^d$. If $r=1$ then $\pi_{1,\infty}(X)$ is a finite cyclic group 
of order $d^2$.
\end{lem}
\Proof
Let $e_i\ (1 \le i \le d+r-1), \ell'_0, s', \ell'_\infty$ be the generators corresponding to 
$E_i\ (1 \le i \le d+r-1), \ell'_0, S', \ell'_\infty$, respectively. Let $x=e_{d+r-1}, y=\ell'_0$ and 
$z=e_1$. Then we have $e_i=x^{d+r-i}$ for $d \le i \le d+r-1, s'=1, \ell'_\infty=y^{-1}, e_i=y^i$ for 
$2 \le i \le d$, and $e_d=z^d$. Furthermore, $x^r=y^d=z^d$ and $zy^{d-1}x^{r-1}=z^{2d}$. Thence we 
obtain $z=xy$. So, $\pi_{1,\infty}(X)$ has generators and relations as described in the statement. 
If $r=1$ then it is easy to see that $\pi_{1,\infty}(X)$ is an abelian group generated by $y$. So, 
$\pi_{1,\infty}(X)=H_{1,\infty}(X)$. A direct computation shows that $H_{1,\infty}(X)$ is a finite group 
for $r$ general, and isomorphic to $\Z/d^2\Z$ provided $r=1$.
\QED

Our next result is the following:

\begin{thm}\label{Theorem 1.4.10}
Let $X$ be an affine pseudo-plane of type $(d,r)$ with $r \ge 2$. Then $\rho$ is 
a unique $\A^1$-fibration on $X$.
\end{thm}

\Proof
Suppose that there exists another $\A^1$-fibration $\sigma : X \to B$ which is different 
from the fixed $\A^1$-fibration $\rho : X \to A$. Then $B \cong \A^1$ because $\Pic(X)_\Q=(0)$ 
and every fiber of $\sigma$ is isomorphic to $\A^1$ if taken with the reduced structure. 
Let $M$ be a linear pencil on $V$ spanned by the closures of general fibers of $\sigma$, where 
the notations $V, D$, etc. are the same as in Definition 1.4.5. Then a general member of $M$ meets 
the curve $\ell'_\infty$, for otherwise the $\A^1$-fibrations $\rho$ and $\sigma$ coincide with 
each other. Suppose that $M$ has no base points. Then the curve $\ell'_\infty$ is a cross-section 
of $M$ and $S'+\ell'_0+E_1+\cdots+E_{d+r-1}$ supports a reducible fiber of $M$. Then $r=d=1$. 
Since $d \ge 2$ by the hypothesis, this case does not take place. Hence $M$ has a base point, 
say $P$, on $\ell'_\infty$. Let $Q:=\ell'_\infty\cap S'$. We consider two cases separately.
\svskip

\noindent
{\sc Case} $P\ne Q$.\ Then $\ell'_\infty+S'+\ell'_0+E_1+E_2+\cdots+E_{d+r-1}$ will support a reducible 
member $G_0$ of the pencil $M$. Let $s=\is{\ol{F}}{G}$, where $G$ is a general member of $M$. 
By comparing the intersection numbers of $G$ with two fibers of $\rho$, $\ell'_\infty$ and the one 
containing $d\ol{F}$, it follows that $\is{\ell'_\infty}{G}=ds$. Let $\mu$ be the multiplicity 
of $G$ at $P$, where $P$ is a one-place point of $G$. We have $ds \ge \mu$. 
Consider first the case $n=1$. The contraction of $S', \ell'_0, E_2,\ldots,E_{d-1}$ makes $E_d$ a 
$(-1)$ curve meeting three components $\ell'_\infty,E_1,E_{d+1}$, and this is impossible. 
So, suppose $n \ge 2$. We need an argument for this case too, which we make use of in the case 
$P=Q$. The elimination of the base points of $M$ will be achieved by blowing up 
the point $P$ and its infinitely near points. After the elimination of the base points of $M$, 
the proper transform $\wt{M}$ gives rise to a $\BP^1$-fibration, and the proper transform of 
$\ell'_\infty$ is a unique $(-1)$ component. If $gs > \mu$ then the point $P$ and its infinitely near 
point of the first order lying on $\ell'_\infty$ are blown up. Hence the proper transform of $\ell'_\infty$ 
is not a $(-1)$ curve. This implies that $ds = \mu$. Let $E$ be the exceptional curve arising from the 
blowing-up of $P$ and let $M'$ be the proper transform of $M$. Then $E$ is contained in the member $G'_0$ 
of $M'$ corresponding to $G_0$ of $M$. In fact, we have otherwise $ds=1$, which is impossible because 
$d \ge 2$. Now contract $\ell'_\infty$ and take the image of $E$ instead of $\ell'_\infty$. Then we have 
the same dual graph as Figure 1 with $\sis{S'}=-(n-1)$. By repeating this argument, we reach to a 
contradiction. 
\svskip

\noindent
{\sc Case} $P=Q$.\ As above, let $G_0$ be a reducible member of $M$ containing 
$S'+\ell'_0+E_1+E_2+\cdots+E_{d+r-1}$. If $\ell'_\infty$ is not contained in $G_0$, the elimination of 
the base points of $M$, which is achieved by blowing up the point $P=Q$ and its infinitely near points, 
yields a $\BP^1$-fibration in which the fiber corresponding to $G_0$ is a reducible fiber not 
containing any $(-1)$ curve. This is a contradiction. Hence $\ell'_\infty$ is contained in $G_0$. 
So, $G_0$ is supported by $S'+\ell'_0+E_1+E_2+\cdots+E_{d+r-1}+\ell'_\infty$. Now apply the elementary 
transformation with center $P$. Then we obtain the same dual graph as Figure 1, where $\sis{S'}=-(n+1)$ 
and $\ell'_\infty$ is replaced by the image of $E$. After repeating the elementary transformations 
several times, we are reduced to the case where $P \ne Q$. So, we reach to a contradiction in the 
present case as well. 
\QED

On the other hand, an affine pseudo-plane of type $(d,1)$ has two algebraically independent 
$G_a$-actions. It follows from a more general result in \cite{GM} when we note that the boundary 
divisor $D$ is then a linear chain for the normal compactification in Definition \ref{Definition 1.4.5} 
and that $\pi_{1,\infty}(X)$ is a finite cyclic group by Lemma \ref{Lemma 1.4.9}. .

\begin{thm} \label{Theorem 1.4.11}
Let $X$ be a smooth affine surface. Then $\ML(X)$ is trivial if and only if $X$ has a minimal normal compactification $V$ such that the dual graph of $D:=V-X$ is a linear chain of rational curves and 
$\pi_{1,\infty}(X)$ is a finite group. 
\end{thm}

One of the possible applications of Theorem \ref{Theorem 1.4.10} is about determining the automorphism 
groups of affine pseudo-planes of type $(d,n,r)$ with $r \ge 2$. 

\begin{thm}\label{Theorem 1.4.12}
Let $X$ be an affine pseudo-plane of type $(d,r)$ with $r \ge 2$. Let $\alpha$ be an automorphism 
of $X$. Then the following assertions hold true.
\begin{enumerate}
\item[{\rm (1)}]
$\alpha$ preserves the $\A^1$-fibration $\rho : X \to A$. Namely, there exists an automorphism 
$\beta$ of $A$ such that $\rho\cdot\alpha=\beta\cdot\rho$, where $\beta$ fixes the point $P_0$ of 
$A$ with $\rho^*(P)$ the multiple fiber $dF_0$ of $\rho$.
\item[{\rm (2)}]
By assigning $\beta$ to $\alpha$, we have a group homomorphism $\pi : \Aut(X) \to \Aut(A-\{P_0\})
\cong G_m$, which is surjective if $(d,r)=(d,d-1)$.
\item[{\rm (3)}]
Let $G$ be the kernel of $\pi$. Then $G$ contains the additive group $G_a$ as a subgroup whose action 
gives rise to the $\A^1$-fibration $\rho$. Furthermore, $G$ does not contain a torus as a subgroup.
\end{enumerate}
\end{thm}
\Proof
The assertion (1) and the first part of the assertion (3) follow readily from Theorem 
\ref{Theorem 1.4.10}. The assertion (2) follows from Lemma \ref{Lemma 1.4.13} below. The second part 
of the assertion (3) follows from Lemma \ref{Lemma 1.4.14} below.

\begin{lem}\label{Lemma 1.4.13}
Let $C$ be a rational curve on $\BP^2$ defined by $X_0X_1^{d-1}=X_2^d$ with $d > 2$ and let $X= 
\BP^2-C$. Let $x=X_0/X_1, y=X_2/X_1$ and $f=x-y^d$. Then the following assertions hold true.
\begin{enumerate}
\item[{\rm (1)}]
Let $\Lambda$ be a linear pencil on $\BP^2$ generated by the curves $C$ and $d\ell_1$, where $\ell_1$ 
is defined by $X_1=0$. Then $\Lambda$ defines an $\A^1$-fibration $\rho : X \to A$, where $A \cong \A^1$. This is a unique $\A^1$-fibration on $X$. Hence there is a canonical projection $\pi : \Aut(X) \to 
\Aut(\A^1-\{0\})\cong G_m$, where $0$ corresponds to the point of $A$ over which lies the unique multiple 
fiber $d(\ell_1\cap X)$.
\item[{\rm (2)}]
Define a $G_m$-action $\tau : G_m\times X \to X$ on $X$ by $(X_0,X_1,X_2) 
\mapsto (\lambda^dX_0,X_1,\lambda X_2)$, where $\lambda \in k$. Then $\tau$ defines a torus subgroup 
$T$ of $\Aut(X)$ such that $\pi_T : T \to \Aut(X-\{0\})$ is the $d$-th power mapping.
\item[{\rm (3)}]
Let $\delta$ be a locally nilpotent derivation on $\Gamma(X,\SO_X)$ such that $\delta(x)=dy^{d-1}f^{-1}$ 
and $\delta(y)=f^{-1}$. Then $\delta$ defines a $G_a$-action $\sigma : G_a\times X \to X$ associated with 
$\rho$. Furthermore, $\tau_\lambda^{-1}\delta\tau_\lambda=\lambda^{d+1}\delta$ for $\lambda \in k^*$.
\item[{\rm (4)}]
Let $G$ be the kernel of $\pi$. Then $G$ consists of automorphisms $\alpha$ such that 
\begin{eqnarray*}
\alpha(y)&=&cy+f^{-m}(a_0f^r+a_1f^{r-1}+\cdots+a_r) \\
\alpha(x)&=&x-y^d+\alpha(y)^d,
\end{eqnarray*}
where $c^d=1, a_i \in k\ (0 \le i \le r)$ with $a_0a_r\ne 0$ and $m > r \ge 0$. Hence $\Aut(X)$ is not 
an algebraic group.
\end{enumerate}
\end{lem}
\Proof
Let $Y=\BP^2-(C+\ell_1)$. Then $Y$ is an affine surface with $\Gamma(Y,\SO_Y)=\C[x,y,f^{-1}]$, 
where $f=x-y^d$. Furthermore, $Y$ has an $\A^1$-fibration $\{f=\lambda \mid \lambda \in \C^*\}$ 
which extends to a unique $\A^1$-fibration $\rho : X\to A$. In fact, the $\A^1$-fibration on 
$Y$ is obtained by removing a unique multiple fiber $d\ell_1$ from $\rho$. 

Let $\alpha$ be an automorphism of $X$. Since $\rho$ is unique, the fibers of $\rho$ are 
transformed by $\alpha$ to the fibers of $\rho$. Since the multiple fiber $d\ell_1$ is unique, 
$\alpha$ induces an automorphism $\alpha_A$ which fixes the point $\rho(\ell_1)$. Thence  
follows the assertion (1). It is clear that $\tau_\lambda(f)=\lambda^d f$. Since $A-\{0\}
=\Spec k[f,f^{-1}]$, the assertion (2) follows from this remark.

We shall show the assertion (3). It is clear that $\delta$ is a locally nilpotent derivation 
of $k[x,y,f^{-1}]=\Gamma(Y,\SO_Y)$ because $\delta(f)=0$. The automorphism $\sigma_\lambda$ 
is given by $\sigma_\lambda(y)=y+f^{-1}\lambda$ and 
$\sigma_\lambda(x)=x-y^d+\sigma_\lambda(y)^d$, where $\lambda \in k$. Hence the assertion 
(4) shows that $\sigma_\lambda$ is an automorphism of $X$. 

We shall show the assertion (4). Let $\alpha$ be an automorphism in $G$. Let $t=X_1/X_2$ and 
$u=X_0/X_2$. Since $\alpha$ induces the identity on $A$, it acts along the fibers of $\rho$. 
Since $Y=X-d\ell_1=\Spec k[y,f,f^{-1}]$ and $\alpha(f)=f$, $\alpha$ is written as 
\begin{eqnarray*} 
\alpha(y)&=&cf^ny+f^s(a_0f^r+a_1f^{r-1}+\cdots+a_r) \\
\alpha(x)&=&x-y^d+\alpha(y)^d,
\end{eqnarray*}
where $c \in k^*, a_i \in k\ (0 \le i \le r)$ with $a_0a_r\ne 0$ and $n, s, r \in \Z$ with 
$r \ge 0$. Conversely, given such an automorphism $\alpha$ of $Y$ as written above, we shall 
consider when $\alpha$ extends to an automorphism of $X$. Note that $C\cap\ell_1=(1,0,0)$ 
and $\ell_1\cap X$ is contained in the open set $\BP^2-\ell_2$, where $\ell_2=\{X_2=0\}$. Let 
$Z=\BP^2-(C\cup\ell_2)$ be an open set of $X$ containing $\ell_1\cap X$. Since $t=1/y, u=x/y$ and 
$f=(ut^{d-1}-1)/t^d$, it follows that $t^df$ is a nowhere vanishing function on $Z-\ell_1$. We write
\[
\alpha(y)=c^*t^{-dn-1}+t^{-(s+r)d}\left(a_0^*+a_1^*t^d+\cdots+a_r^*t^{rd}\right),
\]
where $c^*, a_i^*\ (0 \le i \le r)$ are nowhere vanishing functions on $Z$. Since $\alpha(t)$ 
is divisible by $t$ and not divisible by $t^2$ along the curve $\ell_1\cap Z$ and since $d > 2$, 
it follows that $n=0$ and $s+r \le 0$. Hence $c^*=c$. Then we can write
\begin{eqnarray*}
\alpha(u)&=&\frac{x-y^d+\alpha(y)^d}{\alpha(y)}\\
&=& \frac{u+(c^d-1)t^{-(d-1)}+dc^{d-1}g^*t^{-d(s+r)-(d-2)}+\cdots+g^dt^{-d^2(s+r)+1}}
{c+t^{1-d(s+r)}g^*},
\end{eqnarray*}
where $g^*=\left(a_0^*+a_1^*t^d+\cdots+a_r^*t^{rd}\right)$. Since $1-d(s+r) \ge 1$ in the denominator, 
the numerator of $\alpha(u)$ must be regular along the curve $\ell_1\cap Z$. Namely, we have $c^d=1$ 
and $s+r < 0$. Set $s=-m$. Then we obtain a formula for $\alpha$ as stated in the assertion (4).
\QED

\begin{lem}\label{Lemma 1.4.14}
Let $X$ be a smooth affine surface with an $\A^1$-fibration $\rho : X \to B$. Suppose that $\rho$ has 
a multiple fiber $mA$. Then there are no non-trivial torus actions $\sigma$ which stabilize the fibers 
of the $\A^1$-fibration $\rho$ in the sense that $\rho\cdot\sigma=\rho$.
\end{lem}
\Proof
Suppose that $\sigma : G_m\times X \to X$ be a non-trivial action which stabilizes the fibers of 
the $\A^1$-fibration $\rho$. Then there exists an equivariant smooth compactification $V$ of $X$ 
such that the boundary divisor $D:=V-X$ has smooth normal crossings (cf. \cite{Sumihiro}) and that the 
fibration $\rho$ extends to a $\BP^1$-fibration $p : V \to \ol{B}$. We denote by $\sigma$ the extended 
torus action on $V$ which stabilizes the fibers of the fibration $p$. Since each fiber and its 
irreducible components of $\rho$ are stable under $\sigma$, we may and shall assume that the multiple 
fiber $mA$ is irreducible, i.e., $A$ is isomorphic to $\A^1$. Let $F$ be the fiber $p^{-1}(\rho(A))$ 
and let $\ol{A}$ be the closure of $A$. Since $\sigma$ induces a trivial action on $\ol{B}$, it follows 
that each irreducible component of $F$ is $\sigma$-stable and the induced torus action on each 
component is non-trivial. In fact, if $\sigma$ is trivial along an irreducible component, say $G$, 
of $F$. Then $\sigma$ is trivial on an open neighborhood of $G$, hence on $V$ itself. This implies 
that we can contract the components of $F$ to a smooth fiber by repeating the equivariant blowing-downs. 
Conversely, in order to obtain the fiber $F$ from a smooth fiber, we always blow up the fixed points 
which must be either the fixed points of the starting smooth fiber or the intersection points of adjacent 
components. But, since $A \cong \A^1$ and $mA$ is a multiple fiber, we must blow up a non-fixed point 
to obtain the component $\ol{A}$. This is a contradiction. 
\QED

Let $p : V \to C$ be a $\BP^1$-fibration from a smooth projective surface $V$ onto a smooth projective 
curve $C$, and let $G$ be a singular fiber of $p$. An irreducible component $B$ of $G$ is called a 
{\em branching component}\index{branching component} if $B$ meets three other components of $G$. 
For example, in Figure 1 of Definition 1.4.5, the component $E_d$ is a branching component. 
In connection with Theorem \ref{Theorem 1.4.12}, we have the following result.

\begin{thm}\label{Theorem 1.4.15}
Let $X$ be an affine pseudo-plane of type $(d,n,r)$ with $n > 0$ and $r \ge 3$, where we may assume 
$n=1$. Let $V$ be a smooth normal compactification as described in Definition \ref{Definition 1.4.5}. 
Choose a point $R_1$ on $E_{d+1}$ different from $E_d\cap E_{d+1}$ and $E_{d+1}\cap E_{d+2}$ and 
blow up $R_1$ to obtain the exceptional curve $G_{d+2}$. Choose a point $R_2$ which is different 
from the point $G_{d+2}\cap E_{d+1}$ and blow it up to obtain the exceptional curve $G_{d+3}$. 
Repeat this kind of blowing-ups $(s-1)$ times to obtain a linear chain $G_{d+2}, G_{d+3}, \ldots, 
G_{d+s-1}, \ol{G}=G_{d+s}$. Call the resulting surface $W$. Let 
$Y=W-(\ell'_\infty+S'+\ell'_0+E_1+\cdots+E_{d+r-1}+G_{d+2}+\cdots+G_{d+s-1}$. Then $\Aut(Y)$ has 
no torus group as a subgroup.
\end{thm}
\Proof
Crucial in the above construction is that the adjacent two components $E_d$ and $E_{d+1}$ are branching 
components. By the same arguments in the proof of Theorem \ref{Theorem 1.4.10}, we can show that the 
$\A^1$-fibration $\rho_Y : Y \to A$ inherited from the $\A^1$-fibration $\rho : X \to A$ is a unique 
$\A^1$-fibration on $Y$. Hence if there is a torus action $\sigma : G_m \times Y \to Y$, we can extend it 
to a torus action on the surface $W$ which stabilizes each component of the boundary divisor $W-Y$. It 
is then clear that every branching component is pointwise fixed. Now consider the point $E_d\cap E_{d+1}$. 
Since two curves $E_d$ and $E_{d+1}$ with distinct tangential directions are pointwise fixed, the action 
is trivial in an open neighborhood of the point $E_d\cap E_{d+1}$. Hence $\sigma$ is trivial. This 
implies that, with the notations of Theorem \ref{Theorem 1.4.12}, $\pi : \Aut(Y) \to \Aut(A-\{0\})$ is 
trivial, i.e., $\Aut(Y)=G$. Meanwhile, $G$ does not contain any torus group by Lemma \ref{Lemma 1.4.14}.
\QED

We can strengthen Theorem \ref{Theorem 1.4.10} to the following effect.

\begin{lem}\label{Lemma 1.4.16}
Let $X$ be a smooth affine surface with an $\A^1$-fibration $\rho : X \to A$, where $A$ is isomorphic 
to $\A^1$. Let $\{m_1F_1, \ldots, m_nF_n\}$ exhaust all multiple fibers of $\rho$, where the 
multiplicity $m$ is $\gcd(\mu_1,\ldots,\mu_r)$ for a fiber $F=\sum_{i=1}^r\mu_iC_i$ which is a disjoint 
sum of the irreducible components $C_i \cong \A^1$. Suppose that $n \ge 2$. Then there are no 
non-constant morphisms from $\A^1$ to $X$ whose image is transverse to the given $\A^1$-fibration 
$\rho$. In particular, there are no $\A^1$-fibrations whose general fibers are transverse to the given 
$\A^1$-fibration $\rho$.
\end{lem}

\Proof
Let $\varphi : B \to X$ be a non-constant morphism whose image is transverse to the given 
$\A^1$-fibration $\rho$, where $B \cong \A^1$. Let $V$ be a smooth compactification of $X$ such that 
$D:=V-X$ is a divisor with simple normal crossings and $\rho$ extends to a $\BP^1$-fibration 
$p : V \to \ol{A}$, where $\ol{A}\cong \BP^1$. Then $D$ contains a cross-section $S$ of $p$ and the other 
components of $D$ are contained in the fibers of $p$. We may assume that the fiber $F_\infty:=
p^{-1}(P_\infty)$ is a smooth fiber, where $\{P_\infty\}=\ol{A}-A$. Let $\ol{B}$ be the closure of the 
image $\varphi(B)$ in $V$. Then $\ol{B}$ meets the fiber $F_\infty$, for otherwise $\varphi(B)$ is 
contained in a fiber of $\rho$. Since $B \cong \A^1$, the composite of the morphisms $\rho_B:= 
\rho\cdot\varphi : B \to A$ is surjective. Choose coordinates $x, y$ of $B, A$, respectively. 
Then $y=f(x)$, where $y$ is identified with $(\rho_B)^*(y)$ and $f(x)$ is a polynomial in $x$. Let 
$\rho(F_i)$ be defined by $y=a_i$ for $i=1,2$. Then $y-a_i=(f_i(x))^{m_i}$ for some non-constant 
polynomial $f_i(x)$. Then we have a polynomial relation $f_1(x)^{m_1}-f_2(x)^{m_2}=a_2-a_1$. 
Let $m=\gcd(m_1,m_2)$. We consider the cases $m > 1$ and $m=1$ separately. 

Suppose that $m > 1$. Set $m_i=mn_i$ for $i=1,2$. Then $g_1(x)^m-g_2(x)^m=a_2-a_1\ne 0$, where 
$g_i(x)=f_i(x)^{n_i}$. Then we have 
\begin{eqnarray*}
\lefteqn{g_1(x)^m-g_2(x)^m}\\
&&=\left(g_1(x)-g_2(x)\right)\left(g_1(x)-\zeta g_2(x)\right)\cdots\left(g_1(x)-\zeta^{m-1}g_2(x)\right)
=a_2-a_1,
\end{eqnarray*}
where $\zeta$ is a primitive $m$-th root of unity. Since this is a polynomial identity, it follows that 
$g_1(x)-\zeta^{i-1}g_2(x)=c_i\ne 0$ for $1 \le i \le m-1$. Since $m \ge 2$, $g_1(x)$ and $g_2(x)$ are 
constants. This is a contradiction. 

Suppose that $m=1$. We may assume that $m_1 < m_2$. Let $C$ be an affine curve in $\A^2$ defined by 
$X^{m_1}-Y^{m_2}=1$ and let $\ol{C}$ be a projective plane curve defined by 
$Z^{m_2-m_1}X^{m_1}-Y^{m_2}=Z^{m_2}$. Then the relation $f_1(x)^{m_1}-f_2(x)^{m_2}=a_2-a_1$ implies 
that there is a non-constant morphism $\varphi : \A^1 \to C$. We shall show that the curve $C$ is irrational. The affine curve $C$ is smooth, and the curve $\ol{C}$ has a singular point $Q$ defined by 
$(X,Y,Z)=(1,0,0)$. The singularity type at $Q$ is a cuspidal singularity of type $z^{\ell_1}-y^n=0$, 
where $\ell_1=m_2-m_1$ and $n=m_2$. By the Euclidean algorithm, define positive integers 
$\ell_2, \ldots, \ell_s$ and $a_1, \ldots, a_s$ by 
$n=a_1\ell_1+\ell_2\ (0 < \ell_2 < \ell_1), \ell_1=a_2\ell_2+\ell_3\ (0 < \ell_3 < \ell_2), \ldots, 
\ell_{s-2}=a_{s-1}\ell_{s-1}+\ell_s\ (1=\ell_s < \ell_{s-1}), \ell_{s-1}=a_s\ell_s$. Then the 
multiplicity sequence at $Q$ is $\left\{\ell_1^{a_1}, \ell_2^{a_2}, \ldots, \ell_{s-1}^{a_{s-1}}\right\}$, where $\ell^a$ signifies that the singular points of multiplicity $\ell$ appear $a$-times consecutively 
along the successive blowing-ups. Computing the genus drop by those singularities, we find that the 
geometric genus $g$ of $C$ is equal to 
\[
\frac{1}{2}(n-1)(n-\ell_1-1).
\]
Since $n > 1$ and since $n=\ell_1+1$ implies $m_1=1$, we know that $g > 0$ and $C$ is irrational. 
Then the existence of the non-constant morphism $\varphi : \A^1 \to C$ is a contradiction. 
Hence there are no $\A^1$-fibrations whose general fibers are transverse to $\rho$.
\QED

The following result is noteworthy in view of Theorem \ref{Theorem 1.4.10} and Lemma \ref{Lemma 1.4.16}.

\begin{lem}\label{Lemma 1.4.17}
Let $X$ be an affine pseudo-plane with an $\A^1$-fibration $\rho : X \to A$. Then there exists a 
non-constant morphism $\varphi : \A^1 \to X$ whose image is transversal to the $\A^1$-fibration $\rho$.
\end{lem}
\Proof
We consider a smooth compactification $V$ as specified in Definition \ref{Definition 1.4.5}. Let 
$p : V \to \ol{A}$ be the $\BP^1$-fibration which extends the $\A^1$-fibration $\rho$, where 
$\ol{A}\cong \BP^1$. Let $P_\infty=\ol{A}-A$ and let $P_0=\rho(dF)$, where $dF$ is the unique multiple 
fiber of $\rho$. Let $\tau : \wt{A} \to A$ be a $d$-ple cyclic covering ramifying over the points 
$P_0$ and $P_\infty$. Let $\wt{P}_0$ be the unique point of $\wt{A}$ lying over $P_0$. Let $\wt{X}$ be 
the normalization of the fiber product $X\times_A\wt{A}$. Then $\wt{X}$ has an $\A^1$-fibration 
$\wt{\rho} : \wt{X} \to \wt{A}$ which has a unique reduced reducible fiber with all other fibers 
reduced and irreducible. Namely, the fiber $\wt{\rho}^{-1}(\wt{P}_0)$ consists of $d$ copies of 
the affine line, each of which is counted with multiplicity one. If one deletes all components but one 
from $\wt{\rho}^{-1}(\wt{P}_0)$, the surface $\wt{X}$ with $(d-1)$ copies of the affine line removed 
is isomorphic to the affine plane. In particular, $\wt{X}$ is simply connected. Let $Y$ be the affine 
plane obtained in this fashion. The remaining one component $\ell$ of $\wt{\rho}^{-1}(\wt{P}_0)$ is 
considered to be a coordinate line in $Y$ (cf. Theorem \ref{Theorem 1.3.1}). Let $B$ be another 
coordinate line in $Y$ which meets $\ell$ transversally in one point. Let $\varphi : B \to X$ be the 
restriction onto $B$ of the covering morphism $\wt{X} \to X$. Then $\varphi$ is a non-constant morphism.
\QED

Some partial cases of affine pseudo-planes were observed in tom Dieck \cite{tD1} as examples of 
affine surfaces without cancellation property. In order to state tom Dieck's result, we shall recall 
and generalize a little bit his construction. Write $\Sigma_n=\Proj(\SO_{\BP^1}(-n)\oplus\SO_{\BP^1})$ 
as the quotient of $(\A^2\setminus 0)\times\BP^1$ under the relation
\[
(z_0,z_1),[w_0,w_1] \sim(\nu z_0,\nu z_1),[\nu^nw_0,w_1]
\]
for $\nu \in G_m=k^*$. The projection $\{(z_0,z_1),[w_0,w_1]\} \mapsto [z_0,z_1]$ induces a $\BP^1$ 
fibration $p_n : \Sigma_n \to \BP^1$. In the above definition by quotient and in what follows, the 
integer $n$ could be negative. If $n \ge 0$, the curve $w_0=0$ (resp. $w_1=0$) is a section $M_1$ of 
$p_n$ with $\sis{M_1}=n$ (resp. the minimal section $M_0$ with $\sis{M_0}=-n$). Meanwhile, if $n < 0$, 
then the curve $w_0=0$ (resp. $w_1=0$) is the minimal section $M_0$ (resp. a section $M_1$ with 
$\sis{M_1}=|n|$) of $\Sigma_{|n|}$. Let $d \ge 2$ and $r=d+n \ge 1$. With the notations of Lemma 
\ref{Lemma 1.4.6}, we assume that the fiber $\ell_0$ is defined by $z_0=0$. Let $w=w_0/w_1$. Then 
$\{z_0/z_1,w/z_1^n\}$ is a system of local coordinates at the point $M_1\cap\ell_0$ (resp. 
$M_0\cap\ell_0$) if $n \ge 0$ (resp. $n < 0$). Let $\Lambda$ be a linear subsystem of $|M_1+d\ell_0|$ 
if $n \ge 0$ (resp. $|M_0+d\ell_0|$ if $n < 0$) consisting of members which meet the curve $M_1$ 
(resp. $M_0$) at the point $M_1\cap\ell_0$ (resp. $M_0\cap\ell_0$) with multiplicity $r$ if $n \ge 0$ 
(resp. $n < 0$). Then any member of $\Lambda$ is defined by an equation
\[
\frac{w}{z_1^n}\left\{a_0+a_1\left(\frac{z_0}{z_1}\right)+\cdots+
a_{d-1}\left(\frac{z_0}{z_1}\right)^{d-1}+
a_d\left(\frac{z_0}{z_1}\right)^d\right\}+a_{d+1}\left(\frac{z_0}{z_1}\right)^r=0
\]
or equivalently by 
$$
w_0\left(a_0z_1^d+a_1z_0z_1^{d-1}+\cdots+a_{d-1}z_0^{d-1}z_1+a_dz_0^d\right)+a_{d+1}z_0^rw_1=0
\eqno{(1)}
$$
for $(a_0,a_1, \ldots, a_{d+1}) \in \BP^{d+1}$. In fact, it is readily computed that $\dim\Lambda=d+1$. 
So, the curve $C_d$ in Lemma \ref{Lemma 1.4.6} is defined by such an equation. We shall verify the 
following result.

\begin{lem}\label{Lemma 1.4.18}
Let $X=X(d,r)$ be an affine pseudo-plane with $n=r-d > 0$. Let $\sigma : G_m\times X \to X$ be a 
non-trivial action of the algebraic torus $G_m$. Write $X=\Sigma\setminus M_1\cup C_d$ as in Lemma 
\ref{Lemma 1.4.6}. Then the following assertions hold true.
\begin{enumerate}
\item[{\rm (1)}]
The action $\sigma$ induces an action $\sigma : G_m\times\Sigma_n \to \Sigma_n$ such that 
${}^{\sigma(\lambda)}M_i \subseteq M_i$ for $i=0,1$, ${}^{\sigma(\lambda)}C_d \subseteq C_d$ and 
${}^{\sigma(\lambda)}\ell \sim \ell$ for $\lambda \in k^*$, where $\ell$ is a fiber of $\Sigma_n$.
\item[{\rm (2)}]
The curve $C_d$ is defined by an equation
\[
z_1^dw_0+az_0^rw_1=0.
\]
\end{enumerate}
\end{lem}
\Proof
(1)\ \ Let $\rho : X \to A \cong \A^1$ be the unique $\A^1$-fibration (cf. Theorem \ref{Theorem 1.4.10}).
Then the fibers of $\rho$ are permuted by the action $\sigma$. Hence $\sigma$ extends to the 
cross-section $S'$ and sends $S'$ into itself. Let $W$ be a $G_m$-equivariant smooth normal 
compactification of $X$ whose existence is guaranteed by \cite{Sumihiro}. We may assume that 
$W\setminus X$ contains the cross-section $S'$. Let $F_0$ and $F_\infty$ be two fibers of the 
$\BP^1$-fibration $p : W \to \BP^1$ whose supports partly or totally lie outside of $X$, where 
$F_0$ contains the multiple fiber of $\rho$. We may assume that all $(-1)$ components of $F_0$ and 
$F_\infty$ are fixed componentwise under the action $\sigma$. Then we may assume that $F_\infty$ is 
irreducible and $F_0$ minus the component $\ol{F}$ contains no $(-1)$ components, where $\ol{F}\cap X$ 
gives rise to the multiple fiber of $\rho$. Then we may assume that $W\setminus X$ has the dual graph 
as in Definition \ref{Definition 1.4.5}. So, the action $\sigma$ induces a $G_m$-action on $\Sigma_n$ 
such that ${}^{\sigma(\lambda)}(M_1) \subseteq M_1, {}^{\sigma(\lambda)}(C_d)\subseteq C_d$ and 
${}^{\sigma(\lambda)}(\ell_0) \subseteq \ell_0$ because $M_1, C_d, \ell_0$ are the images on 
$\Sigma_n$ of the components $E_1, \ell'_\infty,\ol{F}$, respectively. The minimal section $M_0$ is 
stable under the $\sigma$ action because the minimal section is unique on $\Sigma_n$. 

\noindent
(2)\ \ The $G_m$-action $\sigma$ on $\Sigma_n$ is given as follows in terms of the coordinates.
\[
\mu\cdot\left((z_0,z_1),[w_0,w_1] \right)=\left((\mu^\alpha z_0,\mu^\beta z_1), [\mu^\gamma w_0, 
\mu^\delta w_1]\right)
\]
for $\mu \in k^*$. Since $C_d$ is stable under the $\sigma$-action, the defining equation (1) must be 
semi-invariant. Note that $a_0a_{d+1}\ne 0$ because $C_d$ is irreducible. Hence we obtain 
$\alpha r+\delta=\beta d+\gamma$. Suppose that $(a_1, \ldots, a_d)\ne (0, \ldots, 0)$. Then we have 
an additional relation $\alpha i+\beta(d-i)+\gamma=\beta d+\gamma$ for some $1 \le i \le d$. The 
last relation implies $\alpha=\beta$. So, the first relation gives $\gamma=\alpha n+\delta$. Then we 
have 
\begin{eqnarray*}
\mu\cdot\left((z_0,z_1),[w_0,w_1] \right)&=&\left((\mu^\alpha z_0,\mu^\beta z_1), [\mu^\gamma w_0, 
\mu^\delta w_1]\right)\\
&=&\left((\mu^\alpha z_0,\mu^\beta z_1),[\mu^{\alpha n+\delta}w_0,\mu^\delta w_1]\right)\\
&\sim&\left((z_0,z_1),[\mu^\delta w_0,\mu^\delta w_1]\right)\\
&=&\left((z_0,z_1),[w_0,w_1]\right).
\end{eqnarray*}
Hence the $\sigma$-action is trivial. This proves the second assertion.
\QED

After tom Dieck \cite{tD1}, we denote by $V(d,r)$ the affine pseudo-plane $X(d,r)$ obtained by 
using a curve defined by an equation
$$
z_1^dw_0+z_0^rw_1=0 \eqno{(2)}
$$
where $r-d$ could be negative. Then $V(d,r)$ has a $G_m$-action defined by
$$
\mu\cdot\left((z_0,z_1),[w_0,w_1]\right)=\left((\mu z_0,z_1),[w_0,\mu^{-r}w_1]\right) \eqno{(3)}
$$
for $\mu\in k^*$. One can show that any $G_m$-action on $V(d,r)$ is reduced to the one defined 
by (3) for
\begin{eqnarray*}
\mu\cdot\left((z_0,z_1),[w_0,w_1] \right)&=&\left((\mu^\alpha z_0,\mu^\beta z_1), [\mu^\gamma w_0, 
\mu^\delta w_1]\right)\\
&=&\left((\mu^\alpha z_0,\mu^\beta z_1),[\mu^{\alpha r-\beta d+\delta}w_0,\mu^\delta w_1]\right)\\
&=&\left((\mu^\alpha z_0,\mu^\beta z_1),[\mu^{\beta n+(\alpha-\beta)r+\delta}w_0,\mu^\delta w_1]\right)\\
&\sim& \left((\mu^{\alpha-\beta}z_0,z_1),[\mu^{r(\alpha-\beta)}w_0,w_1]\right)\\
&=&\left((\mu^{\alpha-\beta}z_0,z_1),[w_0,\mu^{-r(\alpha-\beta)}w_1]\right).
\end{eqnarray*}

In tom Dieck \cite{tD1}, the following result is shown, where $\A^1(c)$ stands for the $\A^1$ with 
a $G_m$-action of weight $c$ and $\wt{V}(d,r)$ for the universal covering of $V(d,r)$.

\begin{thm}\label{Theorem 1.4.19}
With the above notations, there are $G_m$-equivariant isomorphisms 
\[
V(d,r)\times \A^1(-s) \cong V(d,s)\times \A^1(-r)
\]
and 
\[
\wt{V}(d,r)\times \A^1(-s) \cong \wt{V}(d,s)\times \A^1(-r)
\]
for arbitrary positive integers $d,r$ and $s$.
\end{thm}

This result inspires us a very interesting question.
\svskip

\noindent
{\bf Question.}\ \ {\em Let $X(d,r)$ be an affine pseudo-plane with $d \ge 2$ and $r \ge 1$. Does this 
surface have cancellation property?}
\svskip

Towards answering this question, we first consider the universal covering $\wt{X}(d,r)$ of an 
affine pseudo-plane $X(d,r)$.

\begin{lem} \label{Lemma 1.4.20}
The following assertions hold true.
\begin{enumerate}
\item[{\rm (1)}]
The universal covering $\wt{X}(d,r)$ is isomorphic to an affine hypersurface in $\A^3=
\Spec k[x,y,z]$ defined by an equation
$$
x^rz+(y^d+a_1xy^{d-1}+\cdots+a_{d-1}x^{d-1}y+a_dx^d)=1. \eqno{(4)}
$$
\item[{\rm (2)}]
The projection $(x,y,z) \mapsto x$ induces an $\A^1$-fibration $\wt{\rho} : \wt{X}(d,r) \to \A^1$ 
such that every fiber except for $\wt{\rho}^{-1}(0)$ is smooth and the fiber $\wt{\rho}^{-1}(0)$ 
consists of $d$ copies of $\A^1$ which are reduced.
\item[{\rm (3)}]
There is a $G_a$-action on $\wt{X}(d,r)$ defined by 
\begin{eqnarray*}
\lefteqn{c\cdot(x,y,z)}\\
&& =(x,y+cx^r,z-x^{-r}\{((y+cx^r)^d+a_1x(y+cx^r)^{d-1}+\cdots+a_dx^d)\\
&& \quad -(y^d+a_1xy^{d-1}+\cdots+a_{d-1}x^{d-1}y+a_dx^d)\}),
\end{eqnarray*}
where $c \in G_a=k$.
\item[{\rm (4)}]
Let $\omega$ be a $d$-th root of unity. Then there exist uniquely determined polynomials 
$p_\omega(x), q_\omega(x) \in k[x]$ satisfying the following conditions.
\begin{enumerate}
\item[{\rm (i)}]
$\deg p_\omega(x) \le r-1$.
\item[{\rm (ii)}]
$p_\omega(0)=\omega$.
\item[{\rm (iii)}]
$x^rq_\omega(x)+p_\omega(x)^d+a_1xp_\omega(x)^{d-1}+\cdots+a_{d-1}x^{d-1}p_\omega(x)+a_dx^d=1$.
\item[{\rm (iv)}]
$p_{\lambda\omega}(\lambda x)=\lambda p_\omega(x), q_{\lambda\omega}(\lambda x)=\lambda^{-r}q_\lambda(x)$ 
for any $d$-th root $\lambda$ of unity.
\end{enumerate}
By making use of these polynomials, we define the morphism 
\[
\varphi_\omega : \A^2\cong \A^1\times G_a \to \wt{X}(d,r), \quad (x,c) \mapsto c\cdot(x,p_\omega(x),
q_\omega(x))
\]
which is an open immersion onto an open set $U_\omega$ which is the complement of 
$\coprod_{\gamma\ne\omega}G_a\cdot(0,\gamma,0)$. The inverse morphism is defined by 
\[
(x,y,z) \mapsto \left\{ \begin{array}{ll} (x,\dps{\frac{y-p_\omega(x)}{x^r}}) & \mbox{if $x \ne 0$}\\
                                          \mbox{}\\
                                          (0, \dps{\frac{-z+q_\omega(0)}{d\omega^{d-1}}}) & \mbox{if $x=0$}
                                          \end{array} \right.
\]
\item[{\rm (5)}]
$\wt{X}(d,r)$ is obtained by gluing together the $d$-copies of the affine plane $\A^2$ by the 
transition functions
\[
g_{\lambda\omega}:=\varphi_\lambda^{-1}\circ\varphi_\omega : \A^1_*\times\A^1 \to \A^1_*\times\A^1,
\quad (x,c) \mapsto (x,c+\frac{p_\omega(x)-p_\lambda(x)}{x^r}).
\]
\item[{\rm (6)}]
The Galois group is a cyclic group $H(d):=\Z/d\Z$ of order $d$ and acts as 
\[
\lambda\cdot\varphi_\omega(x,c)=\varphi_{\lambda\omega}(\lambda x,\lambda^{1-r}c).
\]
\end{enumerate}
\end{lem}
\Proof
(1)\ \ The surface $X(d,r)$ is the complement in $\Sigma_n$ of the curves $C_d$ defined by the equation 
(1) above and the curve $w_0=0$ which is $M_1$ if $r-d\ge 0$ (resp. $M_0$ if $r-d < 0$), where $n=|r-d|$. 
Since $w_0\ne 0$, we can normalize to $w_0=1$. We can then normalize 
\[
w_0\left(a_0z_1^d+a_1z_0z_1^{d-1}+\cdots+a_{d-1}z_0^{d-1}z_1+a_dz_0^d\right)+a_{d+1}z_0^rw_1\ne 0
\]
to the relation 
\[
z_0^rw_1+(a_0z_1^d+a_1z_0z_1^{d-1}+\cdots+a_{d-1}z_0^{d-1}z_1+a_dz_0^d)=1,
\]
where $a_0 \ne 0$. This normalization comes from the defining equivalence relation 
\[
(z_0,z_1),[w_0,w_1] \sim(\nu z_0,\nu z_1),[\nu^nw_0,w_1].
\]
We may assume that $a_0=1$. Now let $H(d)=\Z/d\Z$ be the cyclic group of order $d$. Then $H(d)$ is 
the Galois group of the covering $\wt{X}(d,r) \to X(d,r)$ acting as
\[
\lambda\cdot(x,y,z) \mapsto (\lambda x,\lambda y, \lambda^{-n}z)
\]
for $\lambda \in H(d)$. In terms of $\varphi_\omega$'s, it is written as in the assertion (6). 

\noindent
(3)\ \ Let $\delta$ be a derivation on the coordinate ring defined by $\delta(x)=0, \delta(y)=x^r$ and 
$\delta(z)= -(dy^{d-1}+(d-1)a_1xy^{d-2}+\cdots+a_{d-1}x^{d-1})$. Then $\delta$ is locally nilpotent. 
Hence it defines a $G_a$-action on $\wt{X}(d,r)$ by Lemma \ref{Lemma 1.2.1}, which is as specified as 
in the assertion. 

\noindent
(4)\ \ Write $p_\omega(x)=\omega+c_1(\omega)x+\cdots+c_{r-1}(\omega)x^{r-1}$, where the coefficients are 
to be determined by the relation
$$
x^rq_\omega(x)+p_\omega(x)^d+a_1xp_\omega(x)^{d-1}+\cdots+a_{d-1}x^{d-1}p_\omega(x)+a_dx^d=1
\eqno{(5)}
$$
which is obtained from the equation (4) above by substituting $p_\omega(x), q_\omega(x)$ for $y,z$. 
By the condition (i), it is easy to see that $p_\omega(x)$ is uniquely determined. Namely the coefficients 
$c_1(\omega), \ldots, c_{r-1}(\omega)$ are uniquely determined by putting the coefficients of the terms 
$x^i\ (1 \le i \le r-1)$ to be zero in the left-hand side of the equation (5) above. Then $q_\omega(x)$ 
is uniquely determined as well. By multiplying $\lambda^d=1$ to the relation (5), we obtain 
\begin{eqnarray*}
\lefteqn{(\lambda x)^r\lambda^{-r}q_\omega(\lambda^{-1}(\lambda x))}\\
&&+(\lambda p_\omega(\lambda^{-1}(\lambda x)))^d+
a_1(\lambda x)(\lambda p_\omega(\lambda^{-1}(\lambda x)))^{d-1}+\cdots+a_d(\lambda x)^d=1.
\end{eqnarray*}
Replace $\lambda x$ by $x$ in the above relation. Then the uniqueness of the polynomials 
$p_{\lambda\omega}(x), q_{\lambda\omega}(x)$ imply that 
$p_{\lambda\omega}(x)=\lambda p_\omega(\lambda^{-1}x)$ and $q_{\lambda\omega}(x)=
\lambda^{-r}q_\omega(\lambda^{-1} x)$. Now replace $x$ by $\lambda x$. Then we obtain the relation (iv). 
Note that $\varphi_\omega : \A^2 \to U_\omega$ is injective and $U_\omega\cong \A^2$. Hence 
$\varphi_\omega$ is an isomorphism by \cite{Ax}. The other assertions are verified in a straightforward 
manner.
\QED

\begin{example}\label{Example 1.4.21}
The following are easy cases to compute.
\begin{enumerate}
\item[{\rm (1)}]
If $r=3, d=2$, then $p_\omega(x)=\omega-\frac{a_1}{2}x+\frac{(a_1^2-4a_2)\omega}{8}x^2$ 
and $q_\omega(x)=-\left\{2c_1(\omega)c_2(\omega)+a_1c_2(\omega)+c_2(\omega)^2x\right\}$.
\item[{\rm (2)}]
If $r=4, d=2$, then $p_\omega(x)=\omega-\frac{a_1}{2}x+\frac{(a_1^2-4a_2)\omega}{8}x^2, c_3(\omega)=0$
and $q_\omega(x)=-c_2(\omega)^2x$.
\end{enumerate}
\end{example}
\svskip

Concerning the isomorphism classes of the hypersurfaces $\wt{X}(d,r)$, we have the following result.

\begin{lem}\label{Lemma 1.4.22}
Let $\wt{X}_1(d,r)$ and $\wt{X}_2(d,r)$ be the hypersurfaces defined by the equations 
$x^rz+y^d+a_2x^2y^{d-2}+\cdots+a_dx^d=1$ and $x^rz+y^d+b_2x^2y^{d-2}+\cdots+b_dx^d=1$, respectively. 
Let $f : \wt{X}_1(d,r) \to \wt{X}_2(d,r)$ be an isomorphism. Suppose $r >d$. Then the induced homomorphism 
of the coordinate rings $\varphi:=f^* : k[x,y,z] \to k[x,y,z]$ is given as $\varphi(x)=cx, \varphi(y)
=uy+x^rG(x)$ and 
\begin{eqnarray*}
\varphi(z)&=&c^{-r}z-(cx)^{-r}\left\{\left((uy+x^rG(x))^d+b_2c^2x^2(uy+x^rG(x))^{d-2}+\right.\right.\\
&&\left.\left.+\cdots+b_dc^dx^d\right)-\left(y^d+a_2x^2y^{d-2}+\cdots+a_dx^d\right)\right\},
\end{eqnarray*}
where $c, u \in k^*$ with $u^d=1$ and $G(x) \in k[x]$ which satisfy the condition $a_i=c^iu^{d-i}b_i$.
\end{lem}
\Proof
Note that $f$ preserves the unique $\A^1$-fibrations of $\wt{X}_1(d,r)$ and $\wt{X}_2(d,r)$ as well as 
the reduced fibers . Hence $\varphi(x)=cx$ with $c \in k^*$. Then $\varphi$ induces an automorphism of
the polynomial ring $k[x,x^{-1}][y]$. So, one can write $\varphi(y)=ux^ey+F(x,x^{-1})$ with $u \in k^*, 
e \in \Z$ and $F(x,x^{-1}) \in k[x,x^{-1}]$. Since $\varphi$ is a homomorphism of the coordinate rings, 
it follows that $e \ge 0$ and $F(x,x^{-1})=F(x)\in k[x]$. Then we have 
\[
c^rx^r\varphi(z)+(ux^ey+F(x)^d+b_2c^2x^2(ux^ey+F(x)^{d-2}+\cdots+b_dc^dx^d=1.
\]
This relation should coincides with the relation 
\[
x^rz+y^d+b_2x^2y^{d-2}+\cdots+b_dx^d=1.
\]
The comparison of the coefficients of the terms $y^d$ and $xy^{d-1}$ implies that $u^d=1, e=0$ and 
$F(x)=0$ is divisible by $x^r$. Furthermore, with the hypothesis $r > d$, we readily see that 
$a_i=c^ib_i\ (2 \le i \le d)$ and $\varphi(z)$ is given as in the statement.
\QED

Let $\wt{X}(d,r)$ be the affine hypersurface in $\A^3$ defined by the equation (4) in Lemma \ref{Lemma 
1.4.20} which has the transition functions given in the  assertion (5) of the same lemma. Let 
$\wt{X'}(d,s)$ be a similar affine hypersurface in $\A^3$ with the equation 
\[
x^sz+(y^d+a'_1xy^{d-1}+\cdots+a'_{d-1}x^{d-1}y+a'_dx^d)=1
\]
and the transition functions 
\[
g'_{\lambda\omega}:={\varphi'}_\lambda^{-1}\circ{\varphi'}_\omega : \A^1_*\times\A^1 \to \A^1_*\times\A^1,
\quad (x,c) \mapsto (x,c+\frac{p'_\omega(x)-p'_\lambda(x)}{x^s}).
\]
As in \cite{tD1}, we define a $3$-dimensional affine variety $W(d,r,s)$ by gluing together $d$-copies 
of the affine $3$-space $\{\omega\}\times\A^3\ (\omega\in H(d))$ by the following identification
\[
(\omega,x,c_1,c_2) \sim (\lambda, x, c_1+\frac{p_\omega(x)-p_\lambda(x)}{x^r}, 
c_2+\frac{p'_\omega(x)-p'_\lambda(x)}{x^s}), \quad x \ne 0.
\]
The projection $(\omega,x,c_1,c_2) \mapsto (\omega,x,c_1)$ yields a morphism $\pi_1 : W(d,r,s) \to 
\wt{X}(d,r)$ which is a principal $G_a$-bundle over $\wt{X}(d,r)$ with $G_a$-action coming. Similarly, 
the projection $(\omega,x,c_1,c_2) \mapsto (\omega,x,c_2)$ gives rise to a  principal $G_a$-bundle 
$\pi_2 : W(d,r,s) \to \wt{X'}(d,s)$. The Galois group action of $H(d)$ onto $\wt{X}(d,r)$ and 
$\wt{X'}(d,s)$ as specified in the assertion (6) of Lemma \ref{Lemma 1.4.20} is lifted to $W(d,r,s)$ 
as follows so that $\pi_1$ and $\pi_2$ are $H(d)$-equivariant:
\[
\lambda\cdot(\omega,x,c_1,c_2)=(\lambda\omega,\lambda x, \lambda^{1-r}c_1,\lambda^{1-s}c_2)\quad 
\mbox{for $\lambda\in H(d)$}.
\]

Since every principal $G_a$-bundle over an affine variety is trivial (cf. \cite{Serre}), we 
have splittings of $W(d,r,s)$ in two ways.

\begin{thm}\label{Theorem 1.4.23}
Let $\wt{X}(d,r)$ and $\wt{X}'(d,s)$ be as above. Then we have isomorphisms
\[
\wt{X}(d,r)\times \A^1 \cong W(d,r,s) \cong \wt{X}'(d,s)\times \A^1.
\]
\end{thm}

We would like to descend these isomorphisms down to the level of $H(d)$-quotient varieties. For 
this purpose, one looks for $H(d)$-equivariant sections of $\pi_1 : W(d,r,s) \to \wt{X}(d,r)$ and 
$\pi_2 : W(d,r,s) \to \wt{X}'(d,s)$.

\begin{lem}\label{Lemma 1.4.24}
For the principal $G_a$-bundle $\pi_1 : W(d,r,s) \to \wt{X}(d,r)$, an $H(d)$-equivariant section is 
given by a family of polynomials $\sigma_\omega \in k[x,c], \omega\in H(d)$ satisfying the following 
conditions:
\begin{enumerate}
\item[{\rm (1)}]
For all $\omega, \lambda \in H(d)$ and $(x,c) \in \A^1_*\times\A^1$,
\[
\sigma_\omega(x,c)+\frac{p'_\omega(x)-p'_\lambda(x)}{x^s}=\sigma_\lambda(x,c+
\frac{p_\omega(x)-p_\lambda(x)}{x^r}).
\]
\item[{\rm (2)}]
For $\omega,\lambda \in H(d)$,
\[
\lambda^{1-s}\sigma_\omega(x,c)=\sigma_{\lambda\omega}(\lambda x,\lambda^{1-r}c).
\]
\end{enumerate}
\end{lem}

We can use the relation (2) in the above lemma to compute $\sigma_\lambda$ from $\sigma_1$
$$
\sigma_\lambda(x,c)=\lambda^{1-s}\sigma_1(\lambda^{-1}x,\lambda^{r-1}c). \eqno{(6)}
$$

In terms of the function $\sigma_1$, the conditions (1) and (2) are reformulated as in the following 
result. The proof is essentially the same as in \cite{tD1} if one takes into account the relation (4) (iv) 
of Lemma \ref{Lemma 1.4.20}.

\begin{lem}\label{Lemma 1.4.25}
Given a polynomial $\sigma=\sigma_1 \in k[x,c]$, define polynomials $\{\sigma_\lambda\mid \lambda
\in H(d)\}$ by the equation (6) above. Then the conditions (1) and (2) in Lemma \ref{Lemma 1.4.24}
are satisfied if and only if $\sigma$ satisfies
$$
\lambda^{1-s}x^s\sigma(\lambda^{-1}x,\lambda^{r-1}(c+\frac{p_1(x)-p_\lambda(x)}{x^r}))
=x^s\sigma(x,c)+p'_1(x)-p'_\lambda(x) \eqno{(7)}
$$
for all $\lambda\in H(d), (x,c) \in \A^1_*\times\A^1$.
\end{lem}

If the polynomial $\sigma$ in Lemma \ref{Lemma 1.4.25} exists, then the principal $G_a$-bundle 
$\pi_1 : W(d,r,s) \to \wt{X}(d,r)$ splits $H(d)$-equivariantly. Namely, define isomorphisms on 
the $\omega$-charts $\theta_\omega : (\omega, \A^2)\times \A^1 \to (\omega, \A^3)$ by 
$((\omega,x,c_1),c_2) \mapsto (\omega,x,c_1,\sigma_\omega(x,c_1)+c_2)$. Then the $\theta_\omega$ 
glue together to give an $H(d)$-equivariant isomorphism $\theta : \wt{X}(d,r)\times \A^1(1-s) \to 
W(d,r,s)$, where $\A^1$ signifies $\A^1$ considered as an $H(d)$-module (i.e., $H(d)$-vector space of 
dimension $1$) with weight $1-s$. 

In order to find solutions of $\sigma(x,c)$ in $k[x,c]$ is a kind of Diophantine Problem, and there 
are no answers to the existence problem of solution in general cases but there is an elaborate 
treating by tom Dieck \cite{tD1} in the case where $\wt{X}(d,r)=\wt{V}(d,r)$ and $\wt{X}'(d,s)=
\wt{V}(d,s)$. In general, we write $\sigma(x,c)$ as a polynomial in $c$ with coefficients in $k[x]$ 
$$
\sigma(x,c)=f_0(x)+f_1(x)c+\cdots+f_t(x)c^t. \eqno{(8)}
$$
We shall show the existence of solutions in the cases where $r=s$ and either $\wt{X}(d,r)=\wt{V}(d,r)$ 
or $\wt{X}'(d,r)=\wt{V}(d,r)$.

\begin{lem}\label{Lemma 1.4.26}
The following assertions hold true.
\begin{enumerate}
\item[{\rm (1)}]
Suppose that $\wt{X}(d,r)=\wt{V}(d,r)$, i.e., $p_\lambda(x)=\lambda$ for all $\lambda \in H(d)$ and 
that the $p'_\lambda(x)$ satisfy the relation
\[
p'_1(\lambda x)-p'_\lambda(\lambda x)=p'_1(x)-p'_\lambda(x) \quad \mbox{for any $\lambda \in H(d)$.}
\]
Put $\sigma(x,c)=f_1(x)c$, where 
\[
f_1(x)=\frac{1}{1-\lambda}\{p'_1(\lambda x)-p'_\lambda(\lambda x)\}.
\]
Then $\sigma(x,c)$ satisfies the relation (7) in Lemma \ref{Lemma 1.4.25}. If $d=2$ and $r=3$ then the 
$p_\lambda(x)$ satisfy the relation above for the $p'_\lambda(x)$ 
(see Example \ref{Example 1.4.21}).
\item[{\rm (2)}]
Suppose that $\wt{X}'(d,r)=\wt{V}(d,r)$ and that $d=2, r=3$. With the notations of Example 
\ref{Example 1.4.21}, let 
\[
\sigma(x,c)=\left(\frac{a_1^2-4a_2}{8}\right)^3x^3+\left\{1-\left(\frac{a_1^2-4a_2}{8}\right)x^2
+\left(\frac{a_1^2-4a_2}{8}\right)^2x^4\right\}c.
\]
Then $\sigma(x,c)$ satisfies the relation (7) in Lemma \ref{Lemma 1.4.25}.
\end{enumerate}
\end{lem}
\Proof
(1)\ \ Let 
\[
\sigma(x,c)=f_0(x)+f_1(x)c
\]
and put it in the relation (7) of Lemma \ref{Lemma 1.4.25}. The relation (7) reads as
\begin{eqnarray*}
\lefteqn{\lambda^{1-r}x^rf_0(\lambda^{-1}x)+x^rf_1(\lambda^{-1}x)c+(1-\lambda)f_1(\lambda^{-1}x)}\\
&&=x^r(f_0(x)+f_1(x)c)+p'_1(x)-p'_\lambda(x).
\end{eqnarray*}
Then we may put 
\[
\begin{array}{c}
f_0(\lambda^{-1}x)=\lambda^{r-1}f_0(x), \quad f_1(\lambda^{-1}x)=f_1(x)\\
\mbox{}\\
(1-\lambda)f_1(\lambda^{-1}x)=p'_1(x)-p'_\lambda(x).
\end{array}
\]
Now put $f_0(x)=0$ and 
\[
f_1(x)=\frac{1}{1-\lambda}\left\{p'_1(\lambda x)-p_\lambda(\lambda x)\right\},
\]
where $p'_1(\lambda x)-p'_\lambda(\lambda x)$ has a factor $1-\lambda$. Then, under the hypothesis 
in the assertion (1), $\sigma(x,c)$ satisfies the relation (7).
\svskip

\noindent
(2)\ \ The computation is straightforward.
\QED

A similar argument as in Lemma \ref{Lemma 1.4.26} gives the following result.

\begin{prop} \label{Prop 1.4.27}
Suppose that $d < r < 2d$ and that $\wt{X}(d,r)$ is defined by $x^rz+y^d+ax^d=1$ with $a \in k$. 
Then the following assertions hold true.
\begin{enumerate}
\item[{\rm (1)}]
$p_\lambda(x)=\lambda-\frac{a}{d}\lambda x^d$ for $\lambda \in H(d)$.
\item[{\rm (2)}]
Define $\sigma(x,c)$ and $\tau(x,c)$ as follows:
\begin{eqnarray*}
\sigma(x,c)&=& \frac{1}{1-\lambda}\left\{p_1(\lambda x)-p_\lambda(\lambda x)\right\}c \\
\tau(x,c)&=& -\frac{a^2}{d^2}x^{2d-r}+\left(1+\frac{a}{d}x^d\right)c.
\end{eqnarray*}
Then we have:
\begin{eqnarray*}
\lambda^{1-r}x^r\sigma(\lambda^{-1}x,\lambda^{r-1}(c+\frac{1-\lambda}{x^r}))
&=&x^r\sigma(x,c)+p_1(x)-p_\lambda(x)\\
\lambda^{1-r}x^r\tau(\lambda^{-1}x,\lambda^{r-1}(c+\frac{p_1(x)-p_\lambda(x)}{x^r}))
&=&x^r\sigma(x,c)+1-\lambda
\end{eqnarray*}
for all $\lambda \in H(d)$. Hence we have an isomorphism 
\[
\wt{X}(d,r)\times\A^1\cong \wt{V}(d,r)\times \A^1.
\]
\end{enumerate}
\end{prop}
\Proof
The relation for $\sigma(x,c)$ follows from the assertion (1) of Lemma \ref{Lemma 1.4.26}. The 
relation for $\tau(x,c)$ can be verified by a straightforward computation.
\QED

\newpage

\chapter{Jacobian Conjecture}
\bvskip

{\large
\begin{enumerate}
\item[\S 1.]
Generalities
\item[\S 2.]
Generalized Jacobian Conjecture 
\item[\S 3.]
Affirmative results
\item[\S 4.]
Normalization by \'etale endomorphisms
\item[\S 5.]
Affine pseudo-coverings
\end{enumerate}}

\newpage
In this chapter, the ground field is an algebraically closed field $k$ of characteristic zero. We often 
confuse $k$ with the complex field $\C$ if the use of $\C$ is more suitable. We shall gather together 
some of the well-known results on the Jacobian conjecture in 2.1 and then shift to the generalized 
Jacobian conjecture in 2.2, 2.3 and 2.4. In the last section 2.5, we deal with affine pseudo-coverings, 
and the treatment contains new results.

\section{Generalities}

The following is a renowned 
\svskip

\noindent
{\bf Jacobian Conjecture.}\ \ {\em Let $f_1,\ldots,f_n$ be elements of a polynomial ring in $n$ 
variables $\C[x_1,\ldots,x_n]$ defined over the complex field $\C$. Suppose that the Jacobian determinant 
\[
J(F;X)=
\left|\begin{array}{ccc}
\frac{\partial f_1}{\partial x_1} & \cdots & \frac{\partial f_1}{\partial x_n}\\
\frac{\partial f_2}{\partial x_1} & \cdots & \frac{\partial f_2}{\partial x_n}\\
\cdots & \cdots & \cdots \\
\frac{\partial f_n}{\partial x_1} & \cdots & \frac{\partial f_n}{\partial x_n}
\end{array}
\right|
\in \C^*.
\]
Then $\C[x_1,\ldots,x_n]=\C[f_1,\dots,f_n]$.}
\svskip

Here $F$ and $X$ signifies a set of polynomials $(f_1,\ldots,f_n)$ and a set of variables 
$(x_1,\ldots,x_n)$ and $\C^*$ stands for the set of nonzero complex numbers $\C\setminus\{0\}$. 
We also denote the Jacobian determinant by 
\[
J\left(\frac{f_1,\ldots,f_n}{x_1,\ldots,x_n}\right).
\]

In geometric terms, the conjecture can be expressed as follows. Define a polynomial morphism 
(or endomorphism) $F : \A^n \to \A^n$ by 
\[
(x_1,\ldots,x_n)\mapsto (f_1(x_1,\ldots,x_n),\ldots,f_n(x_1,\ldots,x_n)).
\]
We also denote this morphism by $F=(f_1,\ldots,f_n)$. If $F$ is an automorphism of the affine space 
$\A^n$, there exists then another polynomial morphism $G=(g_1,\ldots,g_n) :\A^n \to \A^n$ such that 
$G\circ F=\id$ and $F\circ G=\id$. By a theorem of Ax \cite{Ax, Borel} which asserts that 
{\em an injective algebraic endomorphism $\varphi :X\to X$ of an algebraic variety $X$ is an 
automorphism}, a polynomial mapping $F=(f_1,\ldots,f_n) : \A^n \to \A^n$ is an automorphism if there 
exists a polynomial morphism $G=(g_1,\ldots,g_n) : \A^n \to \A^n$ such that $G\circ F=\id$. Denote by 
$J(G\circ F;F)$ the Jacobian determinant $J(G;X)$ with the set of polynomials $F=(f_1,\ldots,f_n)$ 
substituted for the variable $X=(x_1,\ldots,x_n)$. Partial differentiation of a composite function 
then yields the following relation
\[
J(G\circ F;F)\cdot J(F;X)=1\ .
\]
Hence the Jacobian determinant $J(F;X)$ is an invertible element of $\C[x_1,\ldots,x_n]$. Since 
$\C[x_1,\ldots,x_n]^*=\C^*$, it follows that {\em if a polynomial endomorphism $F: \A^n \to \A^n$ is an automorphism, then the Jacobian determinant $J(F;X)$ is a nonzero constant.}
\svskip

The Jacobian Conjecture asserts that the converse of this result holds true. Namely, it asserts that 
\svskip

\noindent
{\bf Jacobian Conjecture (the second form).}\ \ {\em A polynomial endomorphism  
$F=(f_1,\ldots,f_n) : \A^n \to \A^n$ is an automorphism if the Jacobian determinant $J(F;X)$ is 
a nonzero constant.}
\svskip

In the case $n=1$, the conjecture holds because the condition $\partial f/\partial x\in \C^*$ for a 
polynomial in one variable $x$
\[
f=a_0x^m+a_1x^{m-1}+\cdots+a_m,\quad a_0\ne 0
\]
implies $\deg f=1$. In the case $n\geq 2$, the conjecture remains unsolved.  

Set $y_i=f_i(x_1,\ldots,x_n)\ (i=1,\ldots,n)$, where we assume that $f_i(0,\ldots,0)=0$. With the 
condition $J(F;X)\in \C^*$, one can express $x_1,\ldots,x_n$ as formal power series in $y_1,\ldots,y_n$. 
In fact, we can determine the coefficients of the terms of a formal power series $\varphi_i$ by the 
{\em method of undetermined coefficients}\index{method of undetermined coefficients}
\[
x_i=\varphi_i(y_1,\ldots,y_n)=\sum_{\alpha=0}^\infty c_\alpha y^\alpha
\]
after substituting $f_i(x_1,\ldots,x_n)$ for $y_i$, where $y^\alpha=y_1^{\alpha_1}\cdots y_n^{\alpha_n}$ 
for $\alpha=(\alpha_1,\ldots,\alpha_n)$. By the inverse mapping theorem, the formal power series 
$\varphi_i(y_1,\ldots,y_n)$ is a holomorphic function in variables $y_1,\ldots,y_n$ in a small open 
neighborhood of the origin. Namely, the polynomial endomorphism $F$ induces a local isomorphism between 
the small open neighborhoods of the origin. Let $P=(a_1,\ldots,a_n)\in \A^n$ and let 
$Q=F(P)=(b_1,\ldots,b_n)$. Then after the change of variables
\[
x_i=x'_i+a_i,\quad y_i=y'_i+b_i,
\]
we have $y'_i=F'_i(x'_1,\ldots,x'_n)$ and $J(F';X')=J(F;X)\in \C^*$. Hence the polynomial endomorphism 
$F : \A^n \to \A^n$ induces a local isomorphism between small open neighborhoods of the point $P$ and 
its image $Q$. This follows from the observation that if $J(F;X)\in \C^*$, the polynomial endomorphism 
$F : \A^n \to \A^n$ induces an isomorphism $T_F : T_{\A^n,P} \to T_{A^n,Q}$ of the tangent spaces or 
an isomorphism $dF : \Omega^1_{\A^n,Q} \to \Omega^1_{\A^n,P}$ of the cotangent spaces. Hence the 
condition $J(F;X)\in \C^*$ is equivalent to the condition that $F: \A^n \to \A^n$ is unramified (or 
\'etale since $F$ is the morphism between two smooth varieties). 

Although $F$ is locally an isomorphism, one cannot conclude that the endomorphism $F$ gives, therefore, 
an isomorphism between $\A^n$ and its image $F(\A^n)$. In fact, there might be two or more points 
$P_1,\ldots,P_d$ of $\A^n$ mapped to the same image $Q$. Note that there is no subvariety mapped to 
a point by $F$, for otherwise a tangential direction along the subvariety is in the kernel of the 
tangential mapping $T_F$ above for a point $P$ of the subvariety. So, only finitely many points are 
mapped to $Q$. If one takes a point $Q$ to be sufficiently general then the number $d$ is equal to the 
degree of the field extension 
\[
[\C(x_1,\ldots,x_n):\C(f_1,\ldots,f_n)].
\]
A polynomial endomorphism $F : \A^n \to \A^n$ is called a {\em local isomorphism}\index{local isomorphism}. 
We also say that $F$ is {\em unramified}\index{endomorphism!unramified} or 
{\em \'etale}\index{endomorphism!\'etale} at every point of $\A^n$. With this terminology, the Jacobian 
conjecture asserts equivalently the following: 
\svskip

\noindent
{\bf Jacobian Conjecture (the third form).}\ \ {\em A polynomial endomorphism $F$ of $\A^n$ is an 
isomorphism if it is unramified everywhere on $\A^n$.}
\svskip

We can define an unramified morphism $\varphi : X \to Y$ between algebraic varieties. A morphism 
$\varphi$ is said to be {\em unramified}\index{morphism!unramified} if, for any point $P\in X$ and 
$Q:=\varphi(P) \in Y$, the morphism $\varphi$ induces a local isomorphism between small open 
neighborhoods of $P$ and $Q$ (the case $k=\C$), or if, for any point $P\in X$ and $Q:=\varphi(P) \in Y$, 
$\varphi$ induces an isomorphism $\wh{\varphi}^* : \wh{\SO}_{Y,Q} \to \wh{\SO}_{X,P}$ of the local 
rings (if $k$ is arbitrary), where $\wh{\SO}_{Y,Q}$ and $\wh{\SO}_{X,P}$ are respectively the completions 
of the local rings $\SO_{Y,Q}$ and $\SO_{X,P}$. An unramified morphism $\varphi: X \to Y$ is a {\em finite}
\index{morphism!finite} morphism if, for any point $Q\in \varphi(X)$, the inverse image $\varphi^{-1}(Q)$ 
consists of $d$ points for certain fixed number $d$. Then it follows that the image $\varphi(X)$ 
coincides with $Y$ itself because a finite morphism is a closed mapping. We then say that the morphism 
$\varphi: X \to Y$ is a {\em finite covering}\index{finite covering} of {\em degree}\index{degree} $d$. 
Such coverings are familiar in the study of the fundamental group. If one notes that $\C^n$ is simply 
connected, it suffices to show that $F$ is a finite morphism. Hence we can generalize the above form of 
the Jacobian conjecture and state the following more general conjecture. 
\svskip

\noindent
{\bf Generalized Jacobian Conjecture.}\ \ {\em Let $\varphi : X \to X$ be an unramified endomorphism. 
Then $\varphi$ is a finite morphism.}
\svskip

The Jacobian conjecture looks, at the first glance, much of analytic nature, but it is more 
algebro-geometric or topological. The author feels that the conjecture is deeply related to a 
(yet non-existent) theory of branched coverings (ramified along the hidden part, i.e., the boundary at 
infinity) of algebraic varieties. In other words, the author believes that the study of unramified 
(or \'etale) endomorphisms of algebraic varieties is significant in order to solve the Jacobian conjecture. 
We shall list up below some of known results on the Jacobian conjecture which have been proved by far. 

\begin{thm}\label{Theorem 2.0.1}
The following assertions hold for a polynomial mapping $F : \A^n \to \A^n$ with $J(F;X)\in \C^*$.
\begin{enumerate}
\item[{\rm (1)}]
The image of the polynomial mapping $F : \A^n \to \A^n$ is a Zariski open set and contains all 
codimension one points.
\item[{\rm (2)}]
If $\C(x_1,\ldots,x_n)=\C(f_1,\ldots,f_n)$ holds additionally, so does the conjecture, where 
$\C(x_1,\ldots,x_n)$ (resp. $\C(f_1,\ldots,f_n)$) is the quotient field of the polynomial ring 
$\C[x_1,\ldots,x_n]$ (resp. $\C[f_1,\ldots,f_n]$). 
\item[{\rm (3)}]
More generally, if $\C(x_1,\ldots,x_n)$ is a Galois extension of $\C(f_1,\ldots,f_n)$, then the 
conjecture holds. 
\end{enumerate}
\end{thm}
\Proof
A related result is taken up again in Lemma \ref{Lemma 2.1.1}. We shall prove only the assertions 
(1) and (3). The assertion (2) follows from (3). 

(1)\ If $\varphi : X \to X$ is an endomorphism of an algebraic variety, we write it as $\varphi : 
X_1 \to X_2$ to distinguish the source from the target, where $X_1=X_2=X$. Similarly, the induced 
endomorphism of the coordinate rings is written as $\varphi^* : R_2 \to R_1$, where $R_1=R_2= 
\Gamma(X,\SO_X)$. Let $R=\C[x_1,\ldots,x_n]$. Then $\varphi^*$ is injective and $\varphi^*(R_2) 
=\C[f_1,\ldots,f_n]$, which we denote by $S$. Let $\gp$ be a height $1$ prime ideal of $S$. Since 
$S$ is factorial, $\gp=(h)$ for $h \in S$. Since $R^*=S^*=\C^*$, $hR$ is a proper ideal of $R$. 
Let $\GP_1, \ldots, \GP_d$ be all the prime divisors of $hR$. Then $\GP_1, \ldots, \GP_d$ have height 
one. Let $\GP=\GP_1$ and let $\gq=\GP\cap S$. Then $\gq \supseteq \gp$. We have only to show that 
$\gq=\gp$. Suppose the contrary. Consider $S/\gq \to R/\GP$. The associated morphism $\Spec R/\GP 
\to \Spec S/\gq$ is the restriction $F_{\GP}$ of $F$ onto the subvariety $V(\GP)$. If 
$\height(\gq)\ge 2$, then the general fiber of $F_{\GP}$ has positive dimension. This contradicts the 
hypothesis that $F$ is unramified. The openness of $F(\A^n)$ follows from the fact that $F$ is 
\'etale (= unramified + flat) and that a flat morphism is an open mapping. 

(3)\ Set $X=\A^n$ and consider a polynomial endomorphism $F : X_1 \to X_2$. Let $G$ be the Galois group 
of the field extension $\C(x_1,\ldots,x_n)/\C(f_1,\ldots,f_n)$. Let $\wt{X}_2$ be the normalization of 
$X_2$ in the function field of $X_1$. Then $\wt{X}_2$ is a normal algebraic variety containing $X_1$ 
as a Zariski open set by Zariski's Main Theorem and the normalization morphism 
$\wt{F} : \wt{X}_2 \to X_2$ is a finite morphism. Then the Galois group $G$ acts on $\wt{X}_2$ 
in such a way that 
\[
\wt{F}\circ\sigma=\wt{F}, \quad \forall \sigma \in G.
\]
Furthermore, the quotient variety of $\wt{X_2}$ under the $G$-action is isomorphic to $X_2$. 
By the assertion (1), the image $F(X_1)$ is a Zariski open set of $X_2$, which contains all codimension 
one points. Consider a codimension one point $Q$ of $X_2$ and its closure 
$W$ in $X_2$ which is a hypersurface of $\A^n$. Since $Q\in F(X_1)$, there exists a codimension one 
point $P$ of $X_1$ such that $Q=F(P)$. If we set $\wt{F}^{-1}(Q)=\{P_1,\ldots,P_d\}$ with $P=P_1$, it is known that $G$ acts transitively on the set $\wt{F}^{-1}(Q)$. Namely, for any point 
$P_i\in \wt{F}^{-1}(Q)$, there exists $\sigma_i\in G$ such that $P_i=\sigma_i(P_1)$. Since $F$ induces 
a local isomorphism between the points $P_1$ and $Q$, $\wt{F}$ induces a local isomorphism between the automorphic image $P_i$ of $P_1$ and $Q$. Hence the finite morphism $\wt{F} : \wt{X_2}\to X_2$ is 
unramified at every codimension one point of $\wt{X_2}$. Since $X_2$ is smooth, the finite morphism 
$\wt{F} : \wt{X}_2\to X_2$ is then unramified by the {\em purity of branch loci}
\index{purity of branch loci}. Now since $\pi_1(X)=\{1\}$, it follows that $\wt{F}$ is an isomorphism. 
Hence $F : X_1 \to X_2$ is also an isomorphism. 
\QED

When $n=2$ we have more specified results. For more comprehensive treatments, see van den Essen 
\cite{E2}.

\begin{thm}\label{Theorem 2.0.2}
Suppose $n=2$. Let $\deg f_1=m, \deg f_2=n$, where $\deg$ stands for the total degree with respect to 
variables $x_1,x_2$. Then the following assertions hold. 
\begin{enumerate}
\item[{\rm (4)}]
In one of the following cases, $\C[f_1,f_2]=\C[x_1,x_2]$ holds {\em (Nakai-Baba \cite{NB})}:
\begin{enumerate}
\item[{\rm (i)}]
Either $m$ or $n$ is a prime number.
\item[{\rm (ii)}]
 Either $m=4$ or $n=4$. 
\item[{\rm (iii)}]
$m=2p>n$ and $p$ is an odd prime.
\end{enumerate}
\item[{\rm (5)}]
If either $m$ or $n$ is a product of at most two prime numbers, then $\C[x_1,x_2]=\C[f_1,f_2]$ holds 
{\em (Applegate-Onishi \cite{AO}).}
\item[{\rm (6)}]
If $\max(m,n)\leq 100$, then $\C[f_1,f_2]=\C[x_1,x_2]$ holds {\em (Moh \cite{Moh}).}
\item[{\rm (7)}]
If $\C[x_1,x_2]\propersupset \C[f_1,f_2]$, $\gcd(\deg f_1,\deg f_2) \geq 16$ {\em (Heitmann \cite{H}).}
\item[{\rm (8)}]
With the affine plane $\A^2$ embedded naturally into the projective plane $\BP^2$, we denote the line 
at infinity by $\ell_{\infty}$. Denote by $\ol{V(f_1)}$ the closure in $\BP^2$ of the affine plane curve 
$V(f_1)$ defined by $f_1=0$. Similarly, we define $\ol{V(f_2)}$. Suppose that 
$\ol{V(f_1)}\cap\ell_\infty$ as well as $\ol{V(f_2)}\cap\ell_\infty$ consists of a single point for 
every pair $(f_1,f_2)$ satisfying $J(F;X)\in\C^*$. Then the Jacobian conjecture in the case $n=2$ holds 
{\em (Abhyankar \cite{Ab}).}
\end{enumerate}
\end{thm}
\svskip

We shall explain below analytic, algebraic and geometric approaches. For a set of polynomials 
$F=(f_1,\ldots,f_n)$, we define a derivation $\delta_i$ of the polynomial 
ring $\C[x_1,\ldots,x_n]$ by $h \mapsto \delta_{f_i}(h)$, where we set 
\[
\delta_i(h)=J\left(\frac{f_1,\ldots,h,\ldots,f_n}{x_1,\ldots,x_n}\right)
\]
by substituting $h$ for $f_i$. If a derivation $\delta$ of  $\C[x_1,\ldots,x_n]$ satisfies the condition 
that $\delta^N(h)=0\ \ (N \gg 0)$ for any element $h \in \C[x_1,\ldots,x_n]$, we say that $\delta$ is 
{\em locally nilpotent}\index{$k$-derivation!locally nilpotent}. If $n=2$, the following result is known 
(cf. Bass-Connell-Wright \cite{BCW}).

\begin{lem}\label{Lemma 2.0.3}
With the above notations, the following assertion holds.
\begin{enumerate}
\item[{\rm (9)}] Suppose that $J(F;X)\in \C^*$ for $F=(f_1,f_2)$. Then $\C[x_1,x_2]=\C[f_1,f_2]$ 
if and only if $\delta_1$ or $\delta_2$ is a locally nilpotent derivation. 
\end{enumerate}
\end{lem}

We shall consider a formulation of the Jacobian conjecture in terms of partial differentiations. 
We consider the following operators $x_i$ and $\partial_i$ which act from the left on the polynomial 
ring $A_n:=\C[x_1,\ldots,x_n]$ in $n$ variables. They are defined by
\[
x_i: f \mapsto x_i f, \quad \partial_i : f \mapsto \frac{\partial f}{\partial x_i}
\]
The operator algebras over $\C$ generated by the operators $x_1,\ldots,x_n, \partial_1,\ldots,
\partial_n$
\[
\SD_n:=\C[x_1,\ldots,x_n,\partial_1,\ldots,\partial_n]
\]
is called the {\em Weyl algebra}\index{Weyl algebra}. It is easy to verify the following commuting 
relations among those operators: D
\[
[\partial_i,x_j]=\delta_{ij}, [\partial_i,\partial_j]=[x_i,x_j]=0\ (\forall i, j), 
\]
where the bracket product of the operators $\rho$ and $\sigma$ is defined by 
$[\rho,\sigma]=\rho\sigma-\sigma\rho$. An element of the Weyl algebra $\SD_n$ is written as
\[
P=\sum_\alpha a_\alpha\partial^\alpha, \quad a_\alpha \in \C[x_1,\ldots,x_n], 
\]
where we put $\partial^\alpha=\partial_1^{\alpha_1}\cdots\partial_n^{\alpha_n}$ for a 
multi-index $\alpha=(\alpha_1,\ldots,\alpha_n)$. Let $|\alpha|=\alpha_1+\cdots+\alpha_n$ and let 
\[
\SD_n(v)=\{\sum_\alpha a_\alpha\partial^\alpha\ ;\ |\alpha| \leq v\}.
\]
Then $\{\SD_n(v)\}_{v\geq 0}$ is an increasing sequence of left $A_n$-modules of the Weyl algebra 
$\SD_n$, which we call the {\em $\Gamma$-filter}\index{$\Gamma$-filter}. Let $\varphi : \SD_n \to \SD_n$ 
be a ring endomorphism of the Weyl algebra $\SD_n$, which is uniquely determined if the elements 
$\varphi(x_i)$ and $\varphi(\partial_j)$ are assigned so that the following relations are satisfied: 
\[
[\varphi(\partial_i),\varphi(x_j)]=\delta_{ij}, 
[\varphi(\partial_i),\varphi(\partial_j)]=[\varphi(x_i),\varphi(x_j)]=0.
\]
It is known that any ring endomorphism of $\SD_n$ is injective. We then have the following conjecture 
due to Diximier. 
\svskip

\noindent
{\bf Diximier Conjecture.}\ \ {\em Any ring endomorphism of the Weyl algebra $\SD_n$ is surjective.}
\svskip

A weakened version of the Diximier conjecture is stated as follows:
\svskip

\noindent
{\bf Weak Diximier Conjecture.}\ \ {\em Any ring endomorphism $\varphi$ of the Weyl algebra $\SD_n$ 
is surjective provided it preserves the $\Gamma$-filter, {\em i.e., }
\[
\varphi(\SD_n(v))\subseteq \SD_n(v)\ (\forall v \geq 0).
\]
}

There are mutual implications between the Diximier conjecture and the Jacobian conjecture.

\begin{lem}\label{Lemma 2.0.4}
The following assertions hold.
\begin{enumerate}
\item[{\rm (10)}]
If the Diximier conjecture holds, so does the Jacobian conjecture {\em (Vas\-erstein-Kac \cite{BCW})}.
\item[{\rm (11)}]
The weak Diximier conjecture holds if and only if the Jacobian conjecture holds {\em (van den Essen 
\cite{E})}.
\end{enumerate}
\end{lem}

We have observed the Jacobian conjecture over an algebraically closed field. But we can state the 
same conjecture over a non-closed field, e.g., over the real field $\R$. In the real case, there is a 
counterexample by Pinchuk \cite{Pin}.

\begin{thm}\label{Theorem 2.0.5}
We have the following result.
\begin{enumerate}
\item[{\rm (12)}]
There exists a pair of polynomials $(f_1,f_2)$ in two variables $x_1,x_2$ with real coefficients 
such that the polynomial mapping $F=(f_1,f_2) : \R^2\to \R^2$ is a local isomorphism 
but not a homeomorphism. The Jacobian determinant is not invertible but does not have solutions in 
$\R^2$.
\end{enumerate}
\end{thm}

Given a polynomial $f \in \C[x_1,\ldots,x_n]$, let $f=f_{(0)}+f_{(1)}+\cdots+f_{(d)}$ be the 
homogeneous decomposition with respect to the total degree, where $f_{(i)}$ is the $i$-th homogeneous 
part. For a system of polynomials $F=(f_1,\ldots,f_n)$, set $F_i=(f_{1,(i)},\ldots,f_{n,(i)})$, 
where $f_{j,(i)}$ is the $i$-th homogeneous part of $f_j$. Then we have a decomposition into the 
homogeneous parts $F=F_{(0)}+\ldots+F_{(d)}$. Set also $\deg F=\max(\deg f_1,\ldots,\deg f_n)$. 
We shall recall first a result of Wang \cite{Wang}.

\begin{thm}\label{Theorem 2.0.6}
Let $F=(f_1,\ldots,f_n) : \A^n\to \A^n$ be a polynomial endomorphism such that $J(F;X)\in\C^*$ 
and $\deg F \le 2$. Then $F$ is an automorphism. 
\end{thm}
\Proof
By a theorem of Ax \cite{Ax}, it suffices to show that $F: \A^n\to \A^n$ is injective. Suppose the 
contrary. Then there exist two points $P,Q$ of $\C^n$ such that $F(P)=F(Q)$. Define $G=(g_1,\ldots,g_n)$ 
by $G(X)=F(X+P)-F(P)$ and let $R:=Q-P$, where $X+P$ or $Q-P$ is the point obtained by the coordinate-wise 
addition or subtraction. Then $G(O)=G(R)=0$. Since $\deg F \le 2$, we have $\deg G \le 2$. Hence we can 
write $G=G_{(1)}+G_{(2)}$. Set $t_0=1/2$. Then we can compute as follows. 
\begin{eqnarray*}
0&=&G(R)=G_{(1)}(R)+G_{(2)}(R)\\
&=&G_{(1)}(R)+2t_0G_{(2)}(R)\\
&=&\frac{d}{dt}\left(tG_{(1)}(R)+t^2G_{(2)}(R)\right)_{t=t_0}\\
&=& \frac{d}{dt}\left(G_{(1)}(tR)+G_{(2)}(tR)\right)_{t=t_0}\\
&=&\left.\frac{d}{dt}G(tR)\right|_{t=t_0}=J(G;X)(t_0R)\cdot R.
\end{eqnarray*}
Since $J(G;X)=J(F;X)$, we have $J(G;X)(t_0R)\ne 0$ and hence $R=0$. But this is absurd. 
\QED

For a given polynomial endomorphism $F=(f_1,\ldots,f_n):\A^n \to \A^n$, we can define a new set of 
polynomials by adding new variables $x_{n+1}, \ldots, x_{n+m}$ and the associated polynomial mapping 
$F^{[m]}=(f_1,\ldots,f_n,x_{n+1},\ldots,x_{n+m}) : \A^{n+m} \to \A^{n+m}$. Then we have the following 
{\em Reduction Theorem}\index{Reduction Theorem} \cite{BCW}. 

\begin{thm}\label{Theorem 2.0.7}
Given a polynomial endomorphism $F :\A^n \to \A^n$, there exist an integer $m > 0$ and polynomial 
automorphisms $G, H : \A^{n+m} \to \A^{n+m}$ such that the composite $F':=G\circ F^{[m]}\circ H$ satisfies 
the condition $\deg F'\leq 3$. If $J(F;X)\in \C^*$ holds, so does the condition $J(F';X^{[m]})\in \C^*$. 
Hence the Jacobian conjecture holds in general if it holds for polynomial endomorphisms satisfying 
$\deg F \leq 3$. 
\end{thm}

This result can be elaborated in the following fashion. 

\begin{thm}\label{Theorem 2.0.8}
The Jacobian conjecture holds in general if it holds for a polynomial endomorphism $F:\A^n \to \A^n$ 
of the special form $F=X+K, K=(k_1,\ldots,k_n)$, where $k_i$ is either $0$ or a homogeneous 
polynomial of degree $3$.
\end{thm}

Dru$\dot{\rm z}$kowski \cite{D1} made the following improvement of Theorem \ref{Theorem 2.0.8}. 

\begin{thm}\label{Theorem 2.0.9}
If the Jacobian conjecture holds for the $k_i$ in the following form 
\[
k_i=\left(\sum_{j=1}^n a_{ji}x_j\right)^3, \quad a_{ji} \in \C
\]
then it holds in general.
\end{thm}

We say that a polynomial mapping $F=X+K$ is of {\em Dru$\dot{\rm z}$kowski type}
\index{Dru$\dot{\rm z}$kowski type} if the homogeneous part $K$ of degree $3$ is the third power of 
a linear form as in the above theorem. Dru$\dot{\rm z}$kowski \cite{D2} proved the following result. 

\begin{lem}\label{Lemma 2.0.10}
A polynomial endomorphism $F=X+K$ of Dru$\dot{\rm z}$kowski type is an automorphism if the coefficient 
matrix $A=(a_{ji})$ of the linear form in $K$ satisfies either $\rank A \leq 2$ or $n-\rank A \leq 2$.
\end{lem}

We shall restrict ourselves to the case of two variables in the following discussions. We employ the 
notations $f,g$ preferably to $f_1,f_2$ and $x,y$ to $x_1,x_2$. We set $\deg f=m, \deg g=n$ and let 
$f_m, g_n$ be the top homogeneous parts of $f,g$. If $m+n>2$, the hypothesis $J(F;X)\in\C^*$ implies 
\[
J\left(\frac{f_m,g_n}{x,y}\right)=0, 
\]
where $F=(f,g)$ and $X=(x,y)$. If $d=\gcd(m,n)$, this implies that 
{\em there exists a homogeneous polynomial $h$ of degree $d$ satisfying $f_m \sim h^{m/d}$ and 
$g_n \sim h^{n/d}$}, where $p \sim q$ signifies $p/q \in \C^*$. Hence if the hypothesis $J(F;X)\in \C^*$ 
implies either $m\mid n$ or $n\mid m$, it will follow that $g_n \sim f_m^{n/m}$ or $f_m\sim g_n^{m/n}$. 
Suppose, for example, $g_n\sim f_m^{n/m}$. Then we can choose $c \in \C$ in such a way that with 
$f_1=f, g_1=g-cf^{n/m}$, we have
\[
J\left(\frac{f_1,g_1}{x,y}\right)=J\left(\frac{f,g}{x,y}\right) \in \C^*
\]
and $\max(\deg f_1,\deg g_1)<\max(\deg f, \deg g)$. Namely the proof of the Jacobian conjecture will 
proceed by induction on $\max(\deg f, \deg g)$. 

Note that the homogeneous polynomials $f_m$ and $g_n$ are decomposed into the products of linear forms 
\[
f_m\sim \prod_{i=1}^m(a_ix+b_iy), \quad g_n\sim\prod_{j=1}^n(c_jx+d_jy).
\]
Now embed the affine plane $\A^2$ into the projective plane $\BP^2$ by $(x,y)\mapsto (x,y,1)$. Then the 
set of points $I(f):=\{(b_i,-a_i,0)\ ;\ 1\leq i\leq m\}$ is the set of points where the closure of 
the affine curve $V(f)$ meets the line at infinity $\ell_\infty:=\{z=0\}$ and the set of points $I(g):=
\{(d_j,-c_j,0)\ ;\ 1 \leq j\leq n\}$ is the set $\ol{V(g)}\cap \ell_\infty$, where $\ol{V(g)}$ is the 
closure of the affine curve $V(g)$. So, the above result on the existence of a homogeneous polynomial 
$h$ of degree $d$ implies that the hypothesis $J(F;X)\in\C^*$ yields the coincidence of two 
sets of intersection $I(f)$ and $I(g)$. 

In the above arguments, we used the weights $\deg x=\deg y=1$. But we may use different weights 
\[
\deg_\omega x=\omega_1,\quad \deg_\omega y=\omega_2,
\]
where $\omega=(\omega_1,\omega_2)$ is a pair of integers, and define the degree of a monomial 
$cx^\alpha y^\beta\ (c\ne 0)$ by 
\[
\deg_\omega cx^\alpha y^\beta=\alpha\omega_1+\beta\omega_2.
\]
With this degree, we can consider the {\em $\omega$-degree}\index{$\omega$-degree} of polynomials, 
{\em $\omega$-homogeneous polynomials}\index{$\omega$-homogeneous polynomials} and the decompositions 
of polynomials into the {\em $\omega$-homogeneous parts}\index{$\omega$-homogeneous parts}. We denote 
the top $\omega$-homogeneous parts of $f$ and $g$ by $f^+_\omega$ and $g^+_\omega$, respectively. 
If $f$ and $g$ are not $\omega$-homogeneous but satisfy the condition $J(F;X)\in\C^*$, then we know
\[
J\left(\frac{f^+_\omega,g^+_\omega}{x,y}\right)=0. 
\]
If we set $m=\deg_\omega f, n=\deg_\omega g, d=\gcd(m,n)$ then there exists, as in the previous case, 
a $\omega$-homogeneous polynomial $h$ such that 
\[
f^+_\omega\sim h^{m/d}, \quad g^+_\omega\sim h^{n/d}
\]
(see Abhyankar\cite{Ab}). If we change the value of $\omega$, the combinations of monomials appearing 
in $f^+_\omega$ and $g^+_\omega$ will change accordingly. Hence we can get more informations on the 
polynomials $f$ and $g$ satisfying the hypothesis $J(F;X)\in\C^*$. The following results of Magnus 
\cite{Magnus} is shown by using the idea of $\omega$-degree (see Nakai-Baba \cite{NB}).

\begin{lem}\label{Lemma 2.0.11}
The following assertion holds.
\begin{enumerate}
\item[{\rm (13)}]
If $J(F;X)\in\C^*$ and $\min(m,n)>1$, then $\gcd(m,n)>1$. Hence if either $m$ or $n$ is a prime number, 
then the Jacobian conjecture holds.
\end{enumerate}
\end{lem}

Given a polynomial $f=\sum_{i,j}c_{ij}x^iy^j$, the set $S(f)=\{(i,j)\ ;\ c_{ij}\ne 0\}$ is called the 
{\em support}\index{support} of $f$. In the first quadrant of the coordinate plane, consider 
the smallest convex polygon containing the origin $(0,0)$ and the points of $S(f)$. We call it the 
{\em Newton polygon}\index{Newton polygon} of $f$ and denote it by $N(f)$. The following result is 
due to Abhyankar \cite{Ab}. 

\begin{lem}\label{Lemma 2.0.12}
The Jacobian conjecture is equivalent to the following condition: For any pair of polynomials $F=(f,g)$ 
satisfying the condition $J(F;X)\in\C^*$, the newton polygon $N(f)$ (or $N(g)$) is a triangle with three points $(0,q), (0,0), (p,0)$ as summits, where $p, q$ are non-negative integers.
\end{lem}

\section{Generalized Jacobian Conjecture}

We shall consider the generalized Jacobian conjecture in a bit elaborated form:
\svskip

\noindent
{\bf Generalized Jacobian Conjecture.}\ \ {\em Let $X$ be a smooth algebraic variety 
defined over an algebraically closed field $k$ of characteristic zero and let 
$\varphi : X \to X$ be an \'{e}tale endomorphism. Then $\varphi$ is an \'etale finite morphism.}
\svskip

We shall begin with the following result, of which the first assertion is already given in Theorem 
\ref{Theorem 2.0.1}.

\begin{lem}\label{Lemma 2.1.1}
Let $X$ and $\varphi$ be the same as in the above conjecture. Then the following assertions 
hold:
\begin{enumerate}
\item[{\rm (1)}]
If $\varphi$ is injective or if $\varphi$ is birational, then $\varphi$ is an automorphism. 
\item[{\rm (2)}]
If the logarithmic Kodaira dimension $\lkd(X)=\dim X$, i.e., $X$ is of general type, 
then $\varphi$ is an automorphism. {\em We refer to \cite{M2} for the definition of logarithmic Kodaira 
dimension and relevant results.}
\item[{\rm (3)}]
If $X$ is complete and has nonzero Euler number $e(X)$ then $\varphi$ is an automorphism.
\end{enumerate}
\end{lem}
\Proof
(1) If $\varphi$ is injective, the assertion follows from Ax's theorem \cite{Ax,Borel}. If 
$\varphi$ is birational, Zariski's Main Theorem implies that $\varphi$ is injective because $X$ 
is smooth and $\varphi$ is quasi-finite. Hence follows the assertion.

\noindent
(2) The assertion follows from \cite{Iitaka}.

\noindent
(3) Let $d:=\deg\varphi$. Since $\varphi$ is a finite morphism when $X$ is complete, we have 
$e(X)=de(X)$. If $e(X) \ne 0$ then $d=1$. Namely $\varphi$ is birational. Hence $\varphi$ is 
an automorphism by (1).
\QED

\noindent
{\bf Remark.}\ (1)\ If $X$ has logarithmic Kodaira dimension $\lkd(X) \ge 0$ and if 
$\varphi : X \to X$ is a dominant morphism, then $\varphi$ is an \'etale endomorphism 
(cf. \cite[Th. 2]{Iitaka}).

(2)\ If $X$ is a commutative group variety, then the {\em multiplication by $m$ }
\index{endomorphism!multiplication by $m$ } endomorphism for a positive integer $m > 1$ is a finite 
\'etale morphism, and it is an automorphism if and only if $X$ is a vector group. Furthermore, 
$\lkd(X)=0$ if and only if $X$ has no unipotent subgroups. It is well-known that $X$ has no unipotent 
groups if and only if $X$ is an extension of an abelian variety by an algebraic torus. 
\svskip

If $(X, \varphi)$ is a pair of a smooth algebraic variety $X$ and a non-finite \'etale endomorphism 
$\varphi$ of $X$, then we say that $(X,\varphi)$ is a counterexample to the generalized Jacobian 
Conjecture (GJC in short). It is clear that if $(X,\varphi)$ is a counterexample to GJC then so is 
$(X\times Y,\varphi\times \id_Y)$ for any smooth variety $Y$. 

The appearance of a finite \'etale endomorphism for an abelian variety or an algebraic torus is 
related to the group structure on $X$ which comes from a lattice structure on the universal covering 
space. Thus we are tempted to look for a counterexample to GJC which is not related to a group structure. 
The first of such examples was constructed in \cite{ME}.

\begin{example}\label{Example 2.1.2}
Let $C$ be a nonsingular cubic curve in $\BP^2$ and let $X=\BP^2-C$. Then the following assertions hold:
\begin{enumerate}
\item[{\rm (1)}]
$\lkd(X)=0$ and $\Pic(X)\cong \Z/3\Z$.
\item[{\rm (2)}]
There exists a surjective, non-finite, \'etale endomorphism $\varphi : X \to X$ of degree $3$.
\item[{\rm (3)}]
Let $\wt{X}$ be the normalization of the lower $X$ in the function field of the upper $X$. Then 
$\wt{X}$ is a smooth affine surface containing $X$ as a Zariski open set, and $\wt{X}-X$ is a disjoint 
union of six affine lines.
\end{enumerate} 
\end{example}
\Proof
We shall prove the assertion (2). Let $\pi : W \to \BP^2$ be a triple covering which ramifies totally 
over the curve $C$. Then $W$ is a smooth projective surface with 
\[
K_W \sim \pi^*(K_{\BP^2}+2H),
\]
where $H$ is a hyperplane in $\BP^2$. Then $K_W \sim -\pi^*(H)$, hence $-K_W$ is ample and $\sis{K_W}=3$. 
So, $W$ is a del Pezzo surface of degree $3$. It is well-known that a del Pezzo surface of degree $3$ is 
obtained from $\BP^2$ by blowing up six points in general position and that there are $27$ lines contained 
in the surface. Conversely, if one contracts mutually disjoint six lines on $W$, then one obtains a 
birational morphism from $W$ to $\BP^2$. The $27$ lines on $W$ are obtained as follows. The curve $C$ has 
$9$ flexes $P_1,\ldots,P_9$. Let $\ell_i\ (1 \le i \le 9)$ be the tangent line of $C$ at the point $P_i$. 
Let $\wt{C}$ be the inverse image of $C$ on $W$ and let $\wt{P}_i$ be a point on $\wt{C}$ lying over $P_i$. 
Then $\pi^{-1}(\ell_i)$ consists of three lines $L^{(j)}_i\ (j=1,2,3)$ such that 
$L^{(1)}_i, L^{(2)}_i, L^{(3)}_i$ and $\wt{C}$ meet each other transversally only in the point $\wt{P}_i$. 
Then $\{L^{(j)}_i\mid 1 \le i \le 9, 1 \le j \le 3\}$ are the $27$ lines on $W$. Choose six disjoint 
lines $L_1,\ldots,L_6$ among them and contract them to obtain a birational morphism 
$\rho : W \to \BP^2$. Let $\wh{C}=\rho(\wt{C})$. Then $\BP^2-\wh{C}$ is isomorphic to $X$ and $\rho$ 
induces an isomorphism \footnote{A closed immersion of an elliptic curve $C$ into $\BP^2$ is given 
by $\Phi_{|3P|} : C \hookrightarrow \BP^2$  with a closed point $P$. If $f : C \to C'$ is an isomorphism 
of elliptic curves with $P':=f(P)$, then there exists an automorphism $g : \BP^2 \to \BP^2$ such that 
$g\cdot \Phi_{|3P|}=\Phi_{|3P'|}\cdot f$. Suppose that $C$ (resp. $C'$) is written as 
$Y^2Z=4X^3-g_2XZ^2-g_3Z^3$ (resp. ${Y'}^2Z'=4{X'}^3-g'_2X'{Z'}^2-g'_3{Z'}^3$ with respect to a system of 
homogeneous coordinates $(X,Y,Z)$ (resp. $(X',Y',Z')$) such that the point $P$ (resp. $P')$ is given 
as $(0,1,0)$. Then the $j$-invariant of $C$ (resp. $C'$) is given as $j= g_2^3(g_2^3-27g_3^2)^{-1}$ 
(resp. $j'= {g'}^3_2({g'}^3_2-27{g'}^2_3)^{-1}$). Since $j=j'$, we have $(g_/g'_2)^3=(g_3/g'_3)^2$. 
Namely, $g_2=cg'_2$ and $g_3=dg'_3$ with $c, d \in k^*$. By a change of coordinates $(X,Y,Z) 
\mapsto (cX, \sqrt{c^3/d}Y,dZ)$, we may assume $g_2=g'_2$ and $g_3=g'_3$. Then the complements 
$\BP^2-C$ and $\BP^2-C'$ are isomorphic to each other. } 
\[
\rho' : W-(\wt{C}+ L_1+\cdots+L_6) \st{\cong} \BP^2-\wh{C}.
\]
Let $\varphi$ be the composite of ${\rho'}^{-1}$ and $\pi |_{W-(\wt{C}+L_1+\cdots+L_6)}$. 
Then $\varphi : X \to X$ is a required \'etale endomorphism. Let $\wt{X}= W-\wt{C}$. Then 
$\rho |_{\wt{X}} : \wt{X} \to X$ is the normalization morphism of $X$ in the function field 
of $W$. Hence $\wt{X}-X$ consists of six affine lines.
\QED

This example with due specialization gives rise to a new counterexample with a $\Q$-homology plane 
\cite{GM2}.

\begin{example}\label{Example 2.1.3}
Let $C:=\{X^3+Y^3+Z^3=0\}$ be a cubic curve in ${\BP}^2$ and let $S:={\BP}^2-C$. Let 
$T:=\{X^3+Y^3+Z^3+W^3=0\}$ be a cubic hypersurface in $\BP^3$ and let $\pi:{\BP}^3 \to {\BP}^2$ be 
the projection $\pi([x,y,z,w])=[x,y,z]$. Denote by $T_0$ the curve $\{W=0\}$ on $T$. Then the following 
assertions hold:
\begin{enumerate}
\item[{\rm (1)}]
$\pi:T-T_0 \to S$ is the universal covering morphism. 
\item[{\rm (2)}]
$T$ is a del Pezzo surface with $K_T={\cal O}(-T_0)$. There exist six points in general position on 
${\BP}^2$, say $P_1,\ldots,P_6$, such that $T$ is obtained from ${\BP}^2$ by blowing up these points. 
Let $E_i$ be the exceptional curve lying over $P_i$. Then each $E_i$ meets $T_0$ transversally in 
exactly one point and hence each $E_i$ is a straight line in $\BP^3$. Let $\tau:T \to {\BP}^2$ be the 
blowing-up morphism and let $C'=\tau(T_0)$. We have an isomorphism 
\[
T-(T_0\cup E_1\cup \cdots\cup E_6) \to {\BP}^2-C'.
\] 
Since $T_0\cong C'$, there is an isomorphism $C\cong C'$. Hence there exists an automorphism of 
${\BP}^2$ which maps $C$ onto $C'$. Clearly, the morphism $\pi:T-(T_0\cup E_1\cup \cdots\cup E_6) \to S$ 
is an {\'e}tale but non-finite morphism. Thus we obtain a non-finite {\'e}tale endomorphism 
$\wt{f}:=\pi\circ\tau^{-1}:S \to S$.
\item[{\rm (3)}]
Consider the following action of ${\Z}/3\Z=\langle\sigma\rangle$ on $T$ induced by an action on $\BP^3$, 
\[
\sigma([x,y,z,w])=[\theta x,\theta^2 y,\theta z,w].
\] 
Here $\theta$ is a primitive cube root of unity. This action has no fixed points on $T-T_0$ and it 
commutes with the projection $\pi : \BP^3 \to \BP^2$ and the covering transformation $h:T-T_0 \to T-T_0$, 
where
\[
\sigma([x,y,z])=[\theta x,\theta^2 y,\theta z]\quad\mbox{and}\quad 
h([x,y,z,w])=[x,y,z,\theta w].
\]
\item[{\rm (4)}]
There exist six disjoint lines $F_1,\ldots,F_6$ on $T$ such that $\sigma$ keeps this set of lines stable. 
We can use these six skew lines on $T$ to get a morphism $g:T \to {\BP}^2$ such that these six lines are 
the exceptional curves for the blowing up morphism $g$. Then $\sigma$ acts fixed-point freely on 
$T-(T_0\cup E_1\cup\ldots\cup E_6)$ and commutes with the projection $\pi : T-T_0 \to S$. Thus we get 
an induced non-finite {\'e}tale endomorphism from $S/(\sigma)$ to itself.
\item[{\rm (5)}]
$V:=S/(\sigma)$ is an NC-minimal ${\Q}$-homology plane with $\lkd(V)=0$.
\end{enumerate}
\end{example}
\Proof
We shall prove the assertions (4) and (5). The rest is rather straightforward to verify.
\svskip

\noindent
(4)\ \ Let $F_1:=\{X+Y=0=Z+W\}$. Then $F_2:=\sigma(F_1)=\{\theta X+Y=0=Z+\theta W\}$ and 
$F_3:=\sigma^2(F_1)=\{X+\theta Y=0=Z+\theta^2 W\}$. It is easy to see that these three lines are 
mutually disjoint. Let $F_4:=\{X+\theta W=0=Y+Z\}$. Then $F_5:=\sigma(F_4)=\{\theta X+W=0=Y+\theta Z\}$ 
and $F_6:=\sigma^2(F_4)=\{X+W=0=\theta Y+Z\}$. Then $F_4,F_5,F_6$ are mutually disjoint and 
$F_1$ is disjoint from each of them. It follows that the six lines $F_i$ are mutually disjoint and 
$\sigma$ keeps the set of these lines stable.
\svskip

\noindent
(5)\ \ Recall from \cite[\S 6]{Fujita} that a pair $(Y,D)$ of a smooth projective surface $Y$ and 
a connected simple normal crossing divisor $D$ on $Y$ with $\lkd(Y-D)\geq 0$ is said to be 
{\em NC-minimal}\index{NC-minimal} if in the Zariski decomposition $K_Y+D=P+N$, the negative part 
$N=\Bk^*(D)$. Fujita has proved that if $(Y,D)$ is not NC-minimal then $Y$ contains a $(-1)$-curve $L$ 
which meets $D$ transversally in only one point which is a point on a maximal rational twig of $D$. 
Further, $\lkd(Y-D-L)=\lkd(Y-D)$ (cf. \cite[Lemma 6.20]{Fujita}). Clearly, since $\lkd(S)=0$ and 
$\sigma$ acts 
fixed-point freely on $S$ we see that $\lkd(V)=0$ by Iitaka's result. Since $e(S)=3$, we get $e(V)=1$. 
Since $\pi_1(S)\cong{\Z}/(3)$, we see that $b_1(V)=0$. Hence $b_2(V)=0$, which shows that $V$ is a 
${\Q}$-homology plane. If $V$ is not NC-minimal, then by Fujita's lemma $V$ contains a curve 
$L\cong{\A^1}$. Then for some integer $\ell>0$ the divisor $\ell L$ is the divisor of a regular function 
$\varphi$ on $V$. The morphism $\varphi:V \to \A^1$ is a ${\A^1}_*$-fibration because $\lkd(V-L)=0$ by 
Fujita's lemma and using Kawamata's inequality \cite{M2}, $\lkd(V-L)\geq \lkd(F)+\lkd({\A}^1_*)$ to the 
morphism $V-L \to {\A^1}_*$, where $F$ is a general fiber. The pull-back of this morphism to $S$ gives 
an ${\A^1}_*$-fibration on $S$. We will show that $S$ does not admit any ${\A^1}_*$-fibration. Suppose 
that $\tau:S \to B$ is an ${\A}^1_*$-fibration, where $B$ is a smooth curve. Since there is no 
non-constant morphism from ${\BP}^2$ to a curve, there are one or two points of indeterminacy in 
${\BP}^2$. After resolving the indeterminacies we get a ${\BP}^1$-fibration on the blown-up surface 
such that the proper transform of the curve $C$ is contained in a singular fiber of the
${\BP}^1$-fibration. But it is well-known that every irreducible component of a ${\BP}^1$-fibration 
on a smooth projective surface is a (smooth) rational curve. This contradiction shows that $V$ is 
NC-minimal.
\QED

In the following two examples the endomorphisms $\varphi$ are finite \'etale endomorphisms, 
but its appearance is less related to the group structure of a surface which has logarithmic 
Kodaira dimension zero. We shall recall some definitions on $\A^1_*$-fibrations which are necessary in 
the subsequent arguments. See \cite{M2} for more explanations and relevant results.

Let $X$ be a smooth affine surface and let $\pi : X \to B$ be an $\A^1_*$-fibration. Here we note 
that $X$ has always an $\A^1_*$-fibration provided $\lkd(X)=1$ and that there are no other 
$\A^1_*$-fibrations. There is a smooth normal compactification $V$ of $X$ such that $\pi$ extends to 
a $\BP^1$-fibration $p : V \to \ol{B}$. Let $D=V-X$, which is a reduced effective divisor with simple 
normal crossings. Since the general fibers of $\pi$ are isomorphic to $\A^1_*$, $D$ contains components 
which are transverse to the fibration $p$. There are two cases to consider. The first case is that 
there are two cross-sections $H_1, H_2$ of $p$, and the second case is that $D$ has a unique 
component $H$ which is transverse to $p$ and is a $2$-section. The $\A^1_*$-fibration $\pi$ is called 
{\em untwisted}\index{$\A^1_*$-fibration!untwisted} (resp. {\em twisted}
\index{$\A^1_*$-fibration!twisted}) in the first (resp. second) case.

\begin{example}\label{Example 2.1.4}
Let $\pi : X \to B$ be an untwisted $\A^1_*$-fibration such that every fiber is isomorphic to $\A^1_*$ if taken with the reduced structure. 
Assume that $X$ has a compactification $V$ such that $V$ is normal and $V-X$ consists of one point 
and a connected divisor. Assume furthermore that $B$ is an elliptic curve and that $\pi$ has no multiple 
fibers. Then there exists a nontrivial \'{e}tale endomorphism $\varphi : X \to X$ if and only if 
there exist an endomorphism $\beta$ of $B$, an invertible sheaf $\SL$ on $B$ of positive degree and 
a positive integer $d > 1$ satisfying the condition: $\SL^{\otimes d}\cong \beta^*(\SL)$. We have 
necessarily $d=\deg\beta$ and $\deg\varphi=d^2$. Furthermore, $\varphi$ is a finite morphism.
\end{example}
\Proof
There exists an invertible sheaf $\SL$ such that $X= \BP(\SO_B\oplus\SL)-(S_0\cup S_1)$, where $S_0, S_1$ 
are the cross-sections of the $\BP^1$-bundle. The invertible sheaf $\SL$ is determined uniquely up to 
isomorphisms by the $\A^1_*$-bundle $X$. We may assume that $(S_0^2) < 0$. Then $S_0$ (resp. $S_1$) 
corresponds to the projection from $\SO_B\oplus\SL$ to $\SO_B$ (resp. $\SL$). Suppose there exists an 
\'{e}tale endomorphism $\varphi : X \to X$. Then $\pi\circ\varphi=\beta\circ\pi$ for an endomorphism 
$\beta$ of $B$. Let now $Y$ be a fiber product $X\times_B{(B,\beta)}$ and $\pi_Y : Y \to B$ be the 
induced $\A^1_*$-fibration. Then $Y= \BP(\SO_B\oplus\beta^*(\SL))-(S_{0,Y}\cup S_{1,Y})$, where 
$S_{0,Y}=S_0\times_B{(B,\beta)}$ and $S_{1,Y}=S_1\times_B{(B,\beta)}$. Let $\beta_Y : Y \to X$ 
be the base change of $\beta$. Then $\varphi$ splits as a composite $\varphi=\beta_Y\circ\psi$, where 
$\psi : X \to Y$ is an \'{e}tale $B$-morphism. Then $\psi$ extends to a $B$-morphism of 
$\BP^1$-bundles $\ol{\psi} : \BP(\SO_B\oplus \SL) \to \BP(\SO_B\oplus \beta^*(\SL))$ such that 
$\ol{\psi}^{-1}(S_{0,Y})=S_0$ and $\ol{\psi}^{-1}(S_{1,Y})=S_1$. Take an affine open covering 
$\SU=\{U_i\}_{i\in I}$ of $B$ such that $\SL\mid_{U_i}=\SO_{U_i}e_i$ for every $i \in I$. 
Write $\SO_B=\SO_Be$. Then $\pi^{-1}(U_i)=\Spec A_i[z_i,z_i^{-1}]$, where $z_i=e_i/e$ and 
$A_i=\Gamma(U_i,\SO_B)$. Let $\{a_{ij}\}$ be the transition functions of $\SL$ with respect to 
$\SU$. Then $\psi\mid_{U_i} : \pi^{-1}(U_i) \to \pi_Y^{-1}(U_i)$ is given by $\beta^*(z_i)=
\lambda_iz_i^d$ with $\lambda_i \in A_i^*$ and $d:=\deg\psi$. Then it is easy to deduce a 
relation $\beta^*(a_{ji})=a_{ji}^d\lambda_j\lambda_i^{-1}$ for $i, j \in I$. This implies that 
$\SL^{\otimes d} \cong \beta^*(\SL)$. Since $\deg \SL^{\otimes d}=d\deg\SL$ and 
$\deg\beta^*(\SL)= \deg\beta\deg\SL$, we have necessarily $d=\deg \beta$.

Conversely, if we are given $\beta, \SL, d$ as above, we can construct an endomorphism 
$\varphi : X \to X$ as a composite of $\beta_Y$ and $\psi$, where $\psi$ is a $d$-th power 
morphism on the generic fibers. This $d$-th power morphism is extended to a $B$-morphism 
$\psi : X \to Y$ because of the condition $\SL^{\otimes d}\cong \beta^*(\SL)$. Since 
$\varphi=\beta_Y\circ\psi$, $\deg \beta_Y=\deg\beta$ and $\deg\psi=d$, we have 
$\deg\varphi=d^2$.
\QED

Taking the quotients of $X$ and $B$ in Example \ref{Example 2.1.4} with respect to suitable involutions, 
we can construct a new example. More precisely, we have the following example (cf. \cite{MM}).

\begin{example}\label{Example 2.1.5}
Let $n$ be a positive integer $> 1$ such that $d:= (n^2-n)/2$ is odd. Let $C$ be an elliptic curve 
which we view as an abelian variety by fixing one point as a point of origin. Let $\SL$ be an invertible 
sheaf on $C$ such that $\deg \SL > 0$ and $n_C^*(\SL) \cong \SL^{\otimes n^2}$. Let 
$\SM := \SL^{\otimes (n^2-n)/2}$. Let $(-1)_C : C \to C$ be the involution $p \mapsto -p$. Then the 
following assertions hold:
\begin{enumerate}
\item[{\rm (1)}]
$n_C^*(\SM) \cong \SM^{\otimes n^2}$ and $(-1)_C^*(\SM) \cong \SM$.
\item[{\rm (2)}]
Let $Y=\BP(\SO_C\oplus \SM)-(S_0\cup S_1)$, where $S_0$ (resp. $S_1$) is the cross-section corresponding 
to the projection $\SO_C\oplus \SM \to \SO_C$ (resp. $\SO_C\oplus\SM \to \SM$). Let $q : Y \to C$ be the 
canonical projection of $\A^1_*$-bundle. Then there exists  an \'{e}tale endomorphism 
$\psi : Y \to Y$ such that $n_C\circ q=q\circ\psi$. 
\item[{\rm (3)}]
Let $\{U_i\}_{i\in I}$ be an open covering of $C$ such that $\SM\mid_{U_i} = \SO_{U_i}e_i$ for every 
$i \in I$. Then $q^{-1}(U_i)=\Spec \SO_{U_i}[e_i,e_i^{-1}]$. Define an involution $\iota : Y \to Y$ 
locally by $\iota^*(e_i)=-e_i$ for $i \in I$. Then the involution $\iota$ is well-defined and satisfies 
the conditions $(-1)_C\circ q = q\circ\iota$ and $\psi\circ\iota=\iota\circ\psi$. 
\item[{\rm (4)}]
Let $B:= C/\langle (-1)_C\rangle$ and let $X:= Y/\langle\iota\rangle$. Then the projection 
$q : Y \to C$ induces an $\A^1_*$-fibration $\pi : X \to B$, where $B \cong \BP^1$ and $\pi$ has four 
multiple fibers of multiplicity two. Furthermore, the multiplication by $n$ endomorphism $n_C$ of $C$ 
induces an endomorphism $\beta : B \to B$.
\item[{\rm (5)}]
The \'{e}tale endomorphism $\psi : Y \to Y$ induces an \'{e}tale endomorphism $\varphi : X \to X$ 
such that $\pi\circ\varphi=\beta\circ\pi$. The endomorphism $\varphi$ is not an automorphism.
\end{enumerate}
\end{example}
\Proof
(1)\ By \cite[Corollary 3, page 59]{Mu}, we have for any invertible sheaf $\SL$
\[
n_C^*(\SL) \cong \SL^{\left(\frac{n^2+n}{2}\right)}\otimes (-1)_C^*\SL^{\left(\frac{n^2-n}{2}\right)}\ .
\]
Since $n_C^*(\SL) \cong \SL^{\otimes n^2}$ we have 
\[
(-1)_C^*\SL^{\left(\frac{n^2-n}{2}\right)}\cong \SL^{\left(\frac{n^2-n}{2}\right)}\ .
\]
Hence we have $n_C^*(\SM) \cong \SM^{\otimes n^2}$ and $(-1)_C^*\SM \cong \SM$.

\noindent
(2)\ See Lemmas 5.3 and 5.4 in \cite{MM}.

\noindent
(3)\ Clear by the definition.

\noindent
(4)\ The fixed points in $C$ by $(-1)_C$ are the $2$-torsion points, and there are four of them. 
The quotient morphism $C \to B$ ramifies at these four points, and $B$ is therefore isomorphic 
to $\BP^1$. Since $(-1)_C\circ n_C= n_C\circ (-1)_C$, $n_C$ induces an endomorphism 
$\beta : B \to B$. 

\noindent
(5)\ The involution $\iota : Y \to Y$ has no fixed points by the definition. Hence the quotient morphism 
$Y \to X$ is a finite \'{e}tale morphism. This implies that the fibers of $q : Y \to C$ lying over the 
four $2$-torsion points of $C$ produce the four multiple fibers of $\pi : X \to B$ with multiplicity $2$. 
The rest of the assertion is now readily ascertained.
\QED

The following two examples show that there do exist counterexamples to GJC in the cases where 
the logarithmic Kodaira dimension is not necessarily zero.

\begin{example}\label{Example 2.1.6}
Let $C$ be a smooth complete curve of genus $g$ and let $T = \Spec \C[\xi,\xi^{-1}]$ be a one-dimensional 
algebraic torus. Let $Q_1$ and $Q_2$ be the points of $T$ defined by $\xi =1$ and $\xi=-1$, respectively. 
Let $P_1$ and $P_2$ be two distinct points of $C$. Let $Y=C\times T$, let $C_i=C\times \{Q_i\}$ and let 
$T_i=\{P_i\}\times T$ for $i=1, 2$. Let $\sigma : Z \to Y$ be the blowing-up of the points $(P_1,Q_1)$ 
and $(P_2,Q_2)$ and let $E_i=\sigma^{-1}((P_i,Q_i))\ (i=1,2)$. Let $X=Z-\sigma'T_1-\sigma'T_2$, where 
$\sigma'T_i$ is the proper transform of $T_i$ by $\sigma$. Let $q : X \to T$ be a morphism induced by 
the projection $C\times T \to T$. 

Let $g : T \to T$ be the endomorphism defined by $g^*(\xi)=\xi^n$ for odd $n > 2$, let 
$\wt{X}=X\times_T(T,g)$ and let $\wt{q} : \wt{X} \to T$ be the canonical projection. Then $\wt{q}$ has 
$2n$ singular fibers $L_{1j}=\wt{q}^*(Q_{1j})$ and $L_{2j}=\wt{q}^*(Q_{2j})\ (1 \le j \le n)$, where 
$Q_{1j}$ and $Q_{2j}$ are defined respectively by $\xi=\omega^{j-1}$ and $\xi=-\omega^{j-1}$ with 
$\omega$ being a primitive $n$-th root of the unity. The fibers $L_{1j}$ and $L_{2j}$ have the same 
forms as the fibers $L_1:=q*(Q_1)$ and $L_2:=q^*(Q_2)$, respectively. Write $L_{1j}=M_{1j}+\Delta_{1j}$ 
and $L_{2j}=M_{2j}+\Delta_{2j}$, where $\Delta_{1j}\cong \Delta_{2j}\cong \A^1$ and $M_{1j}$ and $M_{2j}$ 
are considered as open sets of $C$. Let $X_1:=\wt{X}-\sum_{i=1}^2\sum_{j=2}^n\Delta_{ij}$. Then the 
following assertions hold:
\begin{enumerate}
\item[{\rm (1)}]
$X_1$ is isomorphic to $X$ and the composite $\varphi$ of the open immersion $X_1 \hookrightarrow \wt{X}$ 
and the canonical projection $\wt{X} \to X$ is a surjective non-finite \'etale endomorphism of 
degree $n$. 
\item[{\rm (2)}]
We have $\lkd(X)=1$ if $g > 0$ and $\lkd(X)=0$ if $g=0$. If $g=0$ then $X \cong F_0-(D_1+D_2)$, 
where $F_0 \cong \BP^1\times \BP^1$ and $D_1, D_2$ are the curves of type $(1,1)$.
\item[{\rm (3)}]
In the case $g=0$, $X$ is isomorphic to $\Spec k[x,y,z,z^{-1}]/(xy=z^2-1)$. A surjective \'etale 
endomorphism $\varphi : X:=\Spec k[x',y',z',{z'}^{-1}]/(x'y'={z'}^2-1) \to 
X:= \Spec k[x,y,z,z^{-1}]/(xy=z^2-1)$ is given by 
\[
x=x', \ y=y'({z'}^{2(n-1)}+{z'}^{2(n-2)}+\cdots +{z'}^2+1), \ z={z'}^n,
\]
where $n$ is a positive integer.
\item[{\rm (4)}]
$X$ has an untwisted $\A^1_*$-fibration $\phi : X \to C$ induced by the projection $Y \to C$.
\end{enumerate}
\end{example} 

\noindent
{\bf Remark.}\ \ The non-finite \'etale endomorphism given in the assertion (3) above can be 
generalized to the universal coverings of some classes of affine pseudo-planes. Let $V(d,r)$ be 
tom Dieck's affine pseudo-plane of type $(d,r)$ (cf. Lemma 1.4.18) and let $\wt{V}(d,r)$ be its universal 
covering which is defined by $x^rz+y^d=1$ in $\A^3$. Let $X=\Spec k[x,y,z,y^{-1}]/(x^rz+y^d=1)$. Define 
a surjective \'etale endomorphism $\varphi : X=\Spec k[x',y',z',{y'}^{-1}]/({x'}^rz'+{y'}^d=1) \to 
X=\Spec k[x,y,z,y^{-1}]/(x^rz+y^d=1)$ by $x=x', y={y'}^n, z=z'\left({y'}^{d(n-1)}+{y'}^{d(n-2)}+\cdots
+{y'}^d+1\right)$. Then $\varphi$ has degree $n$, and there are $n^2-n$ affine lines missing in the above 
$X$ for $\varphi$ to be finite.
\svskip

In Example \ref{Example 2.1.6}, the surface $X$ has an untwisted $\A^1_*$-fibration. By taking a quotient 
of the surface $X$ by an involution, we obtain a surface $\wh{X}$ with a twisted $\A^1_*$-fibration.

\begin{example}\label{Example 2.1.7}
Take $C=\BP^1$ in Example \ref{Example 2.1.6}. Let $\eta$ be an inhomogeneous coordinate on $C$ such that 
$\eta=0$ (resp. $\infty$) at $P_1$ (resp. $P_\infty$). Let $\iota : Y \to Y$ be an involution defined by 
$\iota^*(\xi)=\xi^{-1}$ and $\iota^*(\eta)=-\eta$. Since $\iota((P_i,Q_i))=(P_i,Q_i)$ for $i=1, 2$, 
the involution $\iota$ lifts to an involution $\iota$ of $X$ such that 
$\iota\cdot\varphi=\varphi\cdot\iota$. Let $\wh{X}=X/\langle\iota\rangle$. 
Then the following assertions hold:
\begin{enumerate}
\item[{\rm (1)}]
The surface $\wh{X}$ is a smooth affine surface with a twisted $\A^1_*$-fibration 
$\wh{\phi} : \wh{X} \to \BP^1$ which is induced by the untwisted $\A^1_*$-fibration $\phi : X \to C$. 
\item[{\rm (2)}]
The \'etale endomorphism $\varphi : X \to X$ descends to an \'etale endomorphism 
$\wh{\varphi} : \wh{X} \to \wh{X}$ of degree $n$. Furthermore, $\wh{\varphi}$ is not finite. 
\item[{\rm (3)}]
The surface $\wh{X}$ is isomorphic to $F_0-D$ with an irreducible curve $D$, where 
$F_0 \cong \BP^1\times \BP^1$ and $D \sim 2M+\ell$ with $M$ and $\ell$ being the fibers of two distinct 
projections from $F_0$ to $\BP^1$. The $\A^1_*$-fibration $\wh{\phi}$ is given by restricting the 
$\BP^1$-fibration $p:=p_{|\ell|}$ onto $\wh{X}$.
\item[{\rm (4)}]
If $\theta : F_0 \to F_0$ is a double covering ramifying over $\ell_1, \ell_2$ which are the fibers of 
$p$ passing through the points of ramifications of $p|_D : D \to \BP^1$. Then $\theta^*(D)=D_1+D_2$, 
where $D_1\sim D_2 \sim M+\ell$ and $X=F_0-(D_1+D_2)$.
\item[{\rm (5)}]
$\lkd(\wh{X})=-\infty, \Pic(\wh{X})=\Z$ and  $\Gamma(\wh{X},\SO_{\wh{X}})^*=k^*$.
\item[{\rm (6)}]
An explicit construction of $X, \iota, \wh{X}$ and $\wh{\varphi}$ are given as follows:
\[
X=\Spec k[x,y,z,z^{-1}]/(xy=z^2-1), \ \iota\left(\begin{array}{c}x\\ y\\ z\end{array}\right)
=\left(\begin{array}{c}x\\ -yz^{-2}\\ z^{-1}\end{array}\right)
\]
\[
\wh{X}=\Spec k[x,t,u]/(x^2u=t^2-4), \ \wh{\varphi}\left(\begin{array}{c}x\\ t\\ u \end{array}
\right)=\left(\begin{array}{c}x\\ g(t) \\ uh(t)\end{array}\right),
\]
where $g(t), h(t)$ are polynomials in $t$ such that 
\[
g\left(z+\frac{1}{z}\right)=\frac{z^{2n}+1}{z^n},\ h\left(z+\frac{1}{z}\right)=
\left(\frac{z^{2n-2}+\cdots+z^2+1}{z^{n-1}}\right)^2.
\]
\end{enumerate}
\end{example}

In the assertion (4), it is straightforward to ascertain the commutativity $\varphi\cdot\iota=
\iota\cdot\varphi$. Set $A=k[x,y,z,z^{-1}]/(xy=z^2-1)$. Let
\[
t= z+\frac{1}{z} \quad\mbox{and}\quad u=\frac{y^2}{z^2}.
\]
Then the $\iota$-invariant subring $B$ is generated by three elements $x,t,u$ which satisfy a relation 
$x^2u=t^2-4$. In fact, taking the squares of both sides of the relation $xy=z^2-1$, we obtain the said 
relation in $x,t, u$. 

As indicated in the above examples, given a pair $(X,\varphi)$ of a smooth surface $X$ and an 
\'etale endomorphism $\varphi$ together with a finite group action $G$ on $X$ which commutes 
with the endomorphism $\varphi$ and whose fixed point locus is either the empty set or a disjoint 
union of smooth one-dimensional components, we can produce another pair $(\wh{X},\wh{\varphi})$ 
of a smooth surface $\wh{X}$ and an \'etale endomorphism $\wh{\varphi}$, where $\wh{X}=X/G$. 
Reversing this construction, we have the following result:

\begin{thm}\label{Theorem 2.1.8}
Let $X$ be a smooth algebraic variety and let $\varphi : X \to X$ be a non-finite \'etale endomorphism. 
Suppose that $\pi_1(X)$ is a finite group. Let $q : \wt{X} \to X$ be the universal covering morphism. 
Then the following assertions hold:
\begin{enumerate}
\item[{\rm (1)}]
There exists a non-finite \'etale endomorphism $\wt{\varphi} : \wt{X} \to \wt{X}$ such that 
$q\circ\wt{\varphi}=\varphi\circ q$.
\item[{\rm (2)}]
There exists a group endomorphism $\chi : \pi_1(X) \to \pi_1(X)$ such that 
$\wt{\varphi}(gu)=\chi(g)\wt{\varphi}(u)$ for any $g \in \pi_1(X)$ and $u \in \wt{X}$.
\item[{\rm (3)}]
Conversely, if there exists a non-finite \'etale endomorphism $\wt{\varphi} : \wt{X} \to \wt{X}$ 
satisfying the condition that $\wt{\varphi}(gu)=\chi(g)\wt{\varphi}(u)$ with a group endomorphism 
$\chi : \pi_1(X) \to \pi_1(X)$ for $g \in \pi_1(X)$ and $u \in \wt{X}$, then there exists a non-finite 
\'etale endomorphism $\varphi : X \to X$ such that $q\circ\wt{\varphi}=\varphi\circ q$. 
\end{enumerate}
\end{thm}
\Proof
(1)\ Let $Z=(X,\varphi)\times_X(\wt{X},q)$ and let $Z=Z_1\coprod \cdots \coprod Z_r$ be the decomposition 
into connected components. Let $\rho_i : Z_i \to X$ and $\sigma_i : Z_i \to \wt{X}$ be respectively the 
restriction onto $Z_i$ of the projections from $Z$ to the first factor $X$ and the second factor 
$\wt{X}$ of the fiber product. Then $\rho_i$ is a finite \'etale morphism but $\sigma_i$ is a non-finite 
\'etale morphism. In fact, if $\sigma_i$ is finite then $q\circ\sigma_i=\varphi\circ\rho_i$ is a finite 
morphism. Hence $\varphi$ is a finite morphism as well, a contradiction. Since $\rho_i : Z_i \to X$ is a 
finite \'etale covering and since $\wt{X}$ is the universal covering space of $X$ there exists a finite 
\'etale morphism $\tau_i : \wt{X} \to Z_i$ such that $q=\rho_i\circ\tau_i$. Let $\wt{\varphi} 
=\sigma_i\circ\tau_i$. Then $q\circ\wt{\varphi}=\varphi\circ q$ and $\wt{\varphi}$ is a non-finite 
\'etale endomorphism.

\noindent
(2)\ Since $q\circ\wt{\varphi}(gu)=\varphi\circ q(gu)=\varphi\circ q(u)=q\circ\wt{\varphi}(u)$, 
we have $\wt{\varphi}(gu)=\chi(g,u)\wt{\varphi}(u)$ for any $g \in \pi_1(X)$ and $u \in \wt{X}$. 
Fix an element $g \in \pi_1(X)$ and move $u \in \wt{X}$. Then $u \mapsto \chi(g,u)$ is a morphism from 
$\wt{X}$ to a finite group $\pi_1(X)$. Hence $\chi(g,u)$ is independent of the choice of $u \in \wt{X}$. 
So, we have $\wt{\varphi}(gu)=\chi(g)\wt{\varphi}(u)$. Then it is easy to show that $\chi$ is a group 
endomorphism. 
\QED

This result inspires us to think of the following:
\svskip

\noindent
{\bf Conjecture.}\ \ {\em Let $X$ be a smooth algebraic variety and let $q : U \to X$ be the universal 
covering, where $U$ is not necessarily an algebraic variety. Then $X$ has a non-finite \'etale 
endomorphism if and only if $U$ has a non-finite unramified endomorphism.}
\svskip

Theorem 2.1.8 implies that some of the counterexamples already treated above have \'etale, non-finite 
endomorphisms.
 
\begin{thm}\label{Theorem 2.1.9}
Let $\wt{X}$ be either the universal covering of $\BP^2-C$ in Example \ref{Example 2.1.2} or $T-T_0$ 
in Example \ref{Example 2.1.3}. Then $\wt{X}$ has an \'etale, non-finite endomorphism.
\end{thm}

In Theorem \ref{Theorem 2.1.9}, the universal covering $\wt{X}$ is not rational as $\wt{X}$ is the 
complement of a smooth K3-surface. Nevertheless, the following example shows that there exists a 
simply connected, smooth, affine surface $\wt{X}$ which has an \'etale non-finite endomorphism. 

\begin{example}\label{Example 2.1.10}
Let $V_0=\BP^1\times\BP^1$. Let $M_0$ be a cross-section and let $\ell_0, \ell_1, \ell_\infty$ 
be distinct three fibers with respect to the second projection $\pi_2 : \BP^1\times\BP^1 \to 
\BP^1$. Let $\varphi : V \to V_0$ be a sequence of blowing-ups with centers at $\ell_0\cap M_0, 
\ell_1\cap M_0$ and their infinitely near points such that $\varphi^*(\ell_0)=\ell'_0+E_1+2E_2+
2E_3$ and $\varphi^*(\ell_1)=\ell'_1+F_1+2F_2+2F_3$, where $\sis{\ell'_0}=\sis{\ell'_1}=
\sis{E_i}=\sis{F_i}=-2$ for $i=1,2$ and $\sis{E_3}=\sis{F_3}=-1$. Let 
\[
X:= V-\left(\ell_\infty+M'_0+\ell'_0+\ell'_1+E'_1+F'_1+E'_2+F'_2\right).
\]
Hence $X$ has an $\A^1$-fibration $\rho : X \to B$ with two multiple fibers $2E_3\cap X, 
2F_3\cap X$ of multiplicity $2$. Then $X$ has a degree two, non-finite \'etale endomorphism.
\end{example}
\Proof
Let $\sigma : B' \to B$ be a degree two covering ramifying over the point at infinity 
$p_\infty$ and $p_0$, where $p_0=\rho(E_3\cap X)$. Let $\wt{X}$ be the normalization of 
$X\times_BB'$, let $\tau : \wt{X} \to X$ be the composite of the normalization morphism and 
the first projection $X\times_BB' \to X$ and let $\wt{\rho} : \wt{X} \to B'$ be the 
$\A^1$-fibration induced naturally by $\rho$. Then $\wt{\rho}^*(q_0)$ is a disjoint sum 
$G_1+G_2$ of two affine lines and $\tau : \wt{X} \to X$ is a finite \'etale morphism, where 
$q_0$ is a point of $B'$ lying over $p_0$. Then $\wt{X}-G_1 \cong \wt{X}-G_2 \cong X$, and 
$\tau \mid_{\wt{X}-G_1}$ and $\tau\mid_{\wt{X}-G_2}$ induce a non-finite \'etale endomorphism 
of $X$.
\QED

\section{Affirmative results}

Most results ever known in the affirmative are in the surface case. So we restrict ourselves to the 
surface case.

\begin{lem}\label{Lemma 2.2.1}
Let $X$ be a smooth affine surface with $\lkd(X)=-\infty$. Suppose that one of the following conditions 
is satisfied:
\begin{enumerate}
\item[{\rm (1)}]
$X$ is irrational but not elliptic ruled.
\item[{\rm (2)}]
$\Gamma(X,\SO_X)^*\ne k^*$, and $\rank (\Gamma(X,\SO_X)^*/k^*)\geq 2$ if $X$ is rational.
\end{enumerate}
Then any \'etale endomorphism $\varphi : X \to X$ is an automorphism.
\end{lem}
\Proof
{\sc Case (1)}\ \ Since $\lkd(X)=-\infty$, there exists an $\A^1$-fibration $\rho : X \to C$, where 
$C$ is a smooth curve of positive genus and $C=\rho(X)$. Let $\varphi : X \to X$ be an \'etale 
endomorphism. Then there exists an endomorphism $\beta : C \to C$ such that 
$\rho\cdot\varphi=\beta\cdot\rho$. Since the genus of $C$ is greater than $1$ by the hypothesis, 
it follows from the Riemann-Hurwitz formula that $\beta$ is an automorphism. Since $\beta$ is of finite 
order, we may assume after replacing $\varphi$ by its suitable iteration that $\beta$ is the identity. 
Hence $\rho=\rho\circ\varphi$.  Now consider the generic fiber $X_K$ of $\rho$ with the function field 
$K$ of $C$. The endomorphism $\varphi$ induces an \'etale endomorphism $\varphi_K : X_K \to X_K$. Since 
$X_K$ is isomorphic to $\A^1_K$, $\varphi_K$ is an automorphism. This implies that $\varphi$ is a 
birational morphism. Hence it is injective by Zariski's main theorem. So, $\varphi$ is an automorphism 
by Ax's theorem \cite{Ax, Borel}.

\noindent
{\sc Case (2)}\ \ Let $A=\Gamma(X,\SO_X)$. Since $\lkd(X)=-\infty$, $X$ contains a cylinder-like open 
set $U_0\times \A^1=\Spec B[x]$, where $U_0=\Spec B$ is an affine normal curve. We have $A\subset B[x]$ 
and $A^* \subset B^*$. Let $R_0$ be the $k$-subalgebra of $A$ generated by all elements of $A^*$. Since 
$\rank A^*/k^* < \infty$ \footnote{Consider a compactification $V$ of $X$ such that $V-X$ has pure 
codimension one. Write $V-X=\sum_{i=1}^r D_i$. If $f \in A^*$, then the divisor $(f)$ is supported by the 
components $D_i$ and $(f)=(g)$ if and only if $g=cf$ for some $c \in k^*$. Hence $A^*/k^*$ is 
embedded into a free abelian group generated by the $D_i$. Hence $A^*/k^*$ is also a finitely generated 
free abelian group.} , the algebra $R_0$ is finitely generated over $k$. Let $R$ be the normalization 
of $R_0$ in $A$. We have $R \subseteq B$. Let $\ol{C}=\Spec R$ and let $\rho : X \to C \subseteq \ol{C}$ 
be the morphism induced by the injection $R \subset A$, where $C=\rho(X)$. Since $R^* \supseteq A^* 
\propersupset k^*$,  it follows that $\lkd(C) \ge 0$. Let $F$ be a general fiber of $\rho$. By 
Kawamata's addition theorem, we know that $\lkd(F)=-\infty$. Namely, $F \cong \A^1$. Hence $\rho$ is an 
$\A^1$-fibration. Let $\varphi : X \to X$ be an \'etale endomorphism. Since $C$ is irrational or $C$ is a 
rational curve with at least three places at infinity, $\varphi$ induces an endomorphism $\beta$ such 
that $\rho\cdot\varphi=\beta\cdot\rho$. In fact, $\beta$ is an automorphism by the Riemann-Hurwitz 
theorem. The rest of the proof is the same as in the case (1).
\QED

Let $X$ be a smooth affine surface with an $\A^1_*$-fibration $\rho : X \to C$. A singular fiber $S$ is, 
by definition, a fiber which is not isomorphic to $\A^1_*$ as a scheme. Let $S$ be a singular fiber. 
Then, by \cite[Lemma 4]{ME}, $S$ is written as $S=\Gamma+\Delta$, where 
\begin{enumerate}
\item[(1)]
$\Gamma=0, \Gamma=\alpha\Gamma_1$ with $\alpha \ge 1$ and $\Gamma_1\cong \A^1_*$, or 
$\Gamma=\alpha_1\Gamma_1+\alpha_2\Gamma_2$, where $\alpha_1\ge 1, \alpha_2\ge 1, 
\Gamma_1\cong\Gamma_2\cong \A^1$ and $\Gamma_1$ and $\Gamma_2$ meet each other in one point transversally.
\item[(2)]
$\Delta \ge 0$ and $\Supp \Delta$ is a disjoint union of curves isomorphic 
to $\A^1$.
\end{enumerate}

When an \'etale endomorphism $\varphi : X \to X$ is given, we write it $\alpha : X_1 \to X_2$ to 
distinguish the source $X$ from the target $X$.

\begin{lem}\label{Lemma 2.2.2}
Let $\rho : X \to C$ be an untwisted $\A^1_*$-fibration and let $\varphi : X\to X$ be an \'etale 
endomorphism such that $\rho\cdot\varphi=\rho$ and $\codim_X(X-\varphi(X))\ge 2$. Let 
$\nu : \wt{X} \to X$ be the normalization of $X_2$ in the function field of $X_1$. Then we have:
\begin{enumerate}
\item[{\rm (1)}]
There exists an open immersion $\iota : X \hookrightarrow \wt{X}$ such that $\varphi=\nu\cdot\iota$.
\item[{\rm (2)}]
$\nu : \wt{X} \to X$ is an \'etale Galois covering of $X$ with a cyclic group $G$ of order $n$ as the 
Galois group, where $n=\deg \varphi$.
\item[{\rm (3)}]
Let $S=\Gamma+\Delta$ be a singular fiber of $\rho$. Then $\varphi$ is finite over $\Gamma$, {\em i.e.}, 
$\varphi^*(\Gamma)$ is $G$-invariant, and $\varphi$ is totally decomposable over $\Delta$, {\em i.e.}, 
the stabilizer group of each connected component of $\Delta$ is trivial.
\end{enumerate}
\end{lem}

\Proof
Our proof consists of several steps.

\noindent
(1)\ \ Let $K$ be the function field of $C$ and let $X_K$ be the generic fiber of $\rho$. Then 
$X_K=\Spec K[x,x^{-1}]$, and $\varphi$ induces an \'etale $K$-endomorphism 
$\varphi_K : X_{1,K} \to X_{2,K}$. Clearly, $\varphi_K$ is given by a $K$-algebra endomorphism 
$\theta_K : x \mapsto ax^{\pm n}$ of $K[x,x^{-1}]$, where $a \in K$ and $n=\deg\varphi$. Let $G$ be 
the group of $n$-th roots of the unity in $\C$, which is a cyclic group of order $n$. The group $G$ 
acts on $X_{1,K}$ by $(x,\zeta) \mapsto x\zeta$, where $\zeta\in G$, and $X_{2,K}$ is clearly the quotient 
variety $X_{1,K}/G$. The $G$-action on $X_{1,K}$ is extended to an action on the normalization $\wt{X}$ 
and $X_1$ is embedded into $\wt{X}$ as an open set. Note that $\wt{X}$ is $G$-equivariant.

\noindent
(2)\ \ Let $P$ be a closed point of $C$ such that the fiber $F:= \rho^*(P)$ is a smooth fiber. Let 
$\SO:= \SO_{C,P}$ and let $X_{\SO}:= X\times_C\Spec\SO$. Then we can choose the element $x$ above so 
that $X_{\SO}=\Spec\SO[x,x^{-1}]$ and the induced endomorphism $\varphi_{\SO} : X_{\SO} \to X_{\SO}$ is 
given by an $\SO$-endomorphism $x \mapsto ax^{\pm n}$, where $a \in \SO^*$. So the $G$-action extends 
over $X_{\SO}$ and $X_{1,\SO}/G = X_{2,\SO}$. In this case, we have $\wt{X}_{\SO}= X_{\SO}$, 
where $\wt{X}_{\SO}=\wt{X}\times_C\Spec\SO$.

\noindent
(3)\ \ Now let $S:= \rho^*(P)$ be a singular fiber and write $S=\Gamma+\Delta$ as above. Since there 
are no non-constant morphisms from $\A^1$ to $\A^1_*$ and since $\codim_X(X-\rho(X))\ge 2$ by the 
hypothesis, it follows that $\varphi_*(\Gamma)=\Gamma$ and $\varphi(\Delta)=\Delta$ as cycles. In 
particular, $\varphi$ is surjective. We claim that if $\Gamma\ne 0$ then $\Gamma$ is $G$-invariant. 
In fact, take a nonsingular completion $p : V \to B$ as follows. Let $ \mu : \wh{X} \to \wt{X}$ be a 
$G$-equivariant resolution of singularities of $\wt{X}$ such that $X_1$ is an open set of $\wh{X}$ and 
that $\wh{X}-X_1$ is a divisor with simple normal crossings. Then $\wh{X}$ has an $\A^1_*$-fibration 
$\rho\cdot\mu : \wh{X} \to C$. We can find a $G$-equivariant smooth completion $p : V \to B$ of this 
$\A^1_*$-fibration such that $V$ is $G$-equivalent, $V-\wh{X}$ is a divisor with simple normal 
crossings, $B$ is a smooth complete curve containing $C$ as an open set and $p |_{\wh{X}}= \rho\cdot\mu$. 
The $G$-equivariant completion is possible by Sumihiro's theorem \cite{Sumihiro}. Let $\Sigma= p^*(P)$, 
where $P=\rho(S)$. Then $\Sigma\cap X_1=S$ and $\Sigma$ is $G$-invariant. If $\Gamma$ were not 
$G$-invariant then the translation $g^*\Gamma$ of $\Gamma$ by some element $g$ of $G$ would be a divisor 
disjoint from $\Gamma$ and $\Sigma$ would therefore contain a loop. This is impossible.  So, 
$\Gamma$ is $G$-invariant. Now suppose $\Delta\ne 0$, and let $\Delta_1$ be an irreducible component of 
$\Delta$. Since $\Delta_1$ and $\varphi(\Delta_1)$ are isomorphic to $\A^1$ and since 
$\varphi_{\Delta_1} : \Delta_1 \to \varphi(\Delta_1)$ is an \'etale morphism, it is an isomorphism. 
Note that $G$ acts transitively on the components of $\nu^{-1}(\varphi(\Delta_1))$ and that the isotropy 
group of $\Delta_1$ is trivial by the previous remark. Hence $g(\Delta_1)\ne \Delta_1$ for any non-unit 
element $g$ of $G$ and $\nu : \wt{X} \to X$ is \'etale above the singular fiber $S$. So, 
$\nu : \wt{X} \to X$ is a finite \'etale morphism.
\QED

We can obtain some noteworthy consequences from Lemma \ref{Lemma 2.2.2}. We need one preparatory 
result.

\begin{lem}\label{Lemma 2.2.3}
Let $\pi : Y \to X$ be an \'etale finite Galois covering of a smooth algebraic variety $X$ with a 
cyclic group $G$ of order $n$ as the Galois group. Then there exists an invertible $\SO_X$-module 
$\SL$ such that $\SL^{\otimes n} \cong \SO_X$ and $X \cong \Spec(\bigoplus_{i=0}^{n-1}\SL^{\otimes i})$. 
Furthermore, we have $\pi^*(\SL)\cong \SO_Y$.
\end{lem}
\Proof
Since the assertion is of local nature we may assume that $X$ is affine. So, let $X=\Spec A$ and let 
$Y=\Spec B$, where we view $A$ as a subalgebra of $B$. It is well-known that $G$ is written as a group 
scheme in the following form:
\[
\begin{array}{c}
G=\Spec \C[t] \quad\mbox{with}\quad t^n=1, \mu(t)=t\otimes t, \\
\varepsilon(t)=1 \quad\mbox{and}\quad \eta(t)=t^{-1},
\end{array}
\]
where $\mu, \varepsilon$ and $\eta$ are respectively the comultiplication, the augmentation and the 
coinverse. The action of $G$ on $X$ is given in terms of the following coaction:
\[
\Delta : B \to B[t], \ b \mapsto \Delta(b)=\sum_{i=0}^{n-1}\Delta_i(b)t^i.
\]
The property that $\Delta$ is a coaction is equivalent to the following properties:
\begin{enumerate}
\item[(1)]
The mapping $\Delta_i$ defined by $b \mapsto \Delta_i(b)$ is a $\C$-endomorphism of $B$.
\item[(2)]
$\Delta_i\Delta_j=\delta_{ij}\Delta_j$, where $\delta_{ij}=1$ if $i=j$ and $0$ if $i\ne j$, and 
$\sum_{i=0}^{n-1}\Delta_i=1$.
\item[(3)]
$\Delta_i(b_1)\Delta_j(b_2)\in \Delta_{i+j}(B)$ for $b_1, b_2 \in B$, where we replace $i+j$ by an 
integer $\ell$ with $0\le \ell < n$ and $\ell \equiv i+j \pmod{n}$ if $i+j \ge n$.
\end{enumerate}
Set $B_i=\Delta_i(B)$ for $0 \le i < n$. Then $B_i$ is an $A$-module and $B_0=A$, which is the 
$G$-invariant subalgebra of $B$. In view of the above properties, we have  $B=\sum_{i=0}^{n-1}B_i$ 
and $B_i\cdot B_j \subseteq B_{i+j}$. Now the property that $\pi$ is \'etale implies that $B_1$ is a 
projective $A$-module of rank $1$, $B_i\cong B_1^{\otimes i}\ (1 \le i < n)$ and $B_1^{\otimes n}\cong A$. 
Conversely, if $B_1$ is a projective $A$-module of rank $1$ such that $B_1^{\otimes n} \cong A$ then 
$B:=\sum_{i=0}^{n-1}B_1^{\otimes i}$ is given an \'etale $A$-algebra structure if an isomorphism 
$\theta : B_1^{\otimes n}\stackrel{\sim}{\to} A$ is assigned. The group $G$ acts on $B$ as follows: 
$(\sum_{i=0}^{n-1}b_i)^{\zeta}=\sum_{i=0}^{n-1}b_i\zeta^i$ if $b_i \in B_1^{\otimes i}$ and $\zeta$ is a 
primitive $n$-th root of unity. Clearly, $\pi^*\SL \cong \SO_Y$ because $B_1B\cong B$.

\begin{thm}\label{Theorem 2.2.4}
Let $X$ be a smooth affine surface with an $\A^1_*$-fibration $\rho : X \to C$. Assume that $\Pic X=(0)$ 
and $\Gamma(X,\SO_X)^*=\C^*$. Then any \'etale endomorphism $\varphi : X \to X$ is an automorphism 
provided $\rho\circ\varphi=\rho$. 
\end{thm}
\Proof
Since $\Pic(X)=(0)$ and $\Gamma(X,\SO_X)^*=\C^*$, it follows that the $\A^1_*$-fibration $\rho$ is 
untwisted \footnote{Let $X$ be a smooth affine surface with a twisted $\A^1_*$-fibration $\rho : X \to C$.
Since there is only a $2$-section $H$ in the boundary at infinity of $X$ which is transverse to 
the $\BP^1$-fibration extending $\rho$, $\Pic(X)$ modulo the subgroup generated by all fiber components of $\rho$ has a nonzero $2$-torsion element, and $\Pic(X)$ is not zero.} and $\codim_X(X-\varphi(X))\ge 2$ 
\footnote{Suppose that there is a curve $A$ on $X$ such that $A\cap \varphi(X)=\emptyset$. Since 
$\Pic(X)=(0)$, there is an element $a$ of $\Gamma(X,\SO_X)$ such that $A=V(a)$. Then the element 
$\varphi^*(a)$ is invertible on the upper $X$, which is impossible because $\Gamma(X,\SO_X)^*=\C^*$.}. 
Let $\wt{X}$ be the normalization of $X_2$ in the function 
field of $X_1$. By Lemma \ref{Lemma 2.2.2}, the normalization morphism $\nu : \wt{X} \to X$ 
is an \'etale finite Galois covering with a cyclic group of order $n$ as the Galois group. By Lemma 
\ref{Lemma 2.2.3} there exists an invertible sheaf $\SL$ such that $\SL^{\otimes n} \cong \SO_X$ 
and $\wt{X}\cong \Spec(\bigoplus_{i=0}^{n-1}\SL^{\otimes i})$. Since $\Pic X=(0)$ the invertible sheaf 
$\SL$ is trivial. Then $\wt{X}$ is a disjoint union of $n$-copies of $X$. Since $X$ is an dense open 
set of $\wt{X}$, we know that $n=1$. Hence $\varphi$ is an automorphism.
\QED 

In order to look more closely into the case of a smooth affine surface with an $\A^1_*$-fibration, 
we say that a singular fiber $S=\Gamma+\Delta$ is of the {\em simple type}
\index{singular fiber!simple type} (resp. {\em non-simple type}\index{singular fiber!non-simple type}) 
if $\Gamma=\alpha\Gamma_1$ with $\Gamma_1\cong \A^1_*$ and $\Delta=0$ (resp. otherwise). 
Let $\rho : X \to C$ be an untwisted $\A^1_*$-fibration and let $\varphi : X \to X$ be an \'etale 
endomorphism. We need the following result.

\begin{lem}\label{Lemma 2.2.5}
Let $X$ and $\varphi$ be the same as above. Suppose that $\lkd(X)=1$. Then the following assertions hold: 
\begin{enumerate}
\item[{\rm (1)}]
The endomorphism $\varphi$ preserves the $\A^1_*$-fibration $\rho$, i.e., there exists an endomorphism 
$\beta : C \to C$ such that $\beta\cdot\rho=\rho\cdot\varphi$ holds. 
\item[{\rm (2)}]
$\varphi$ sends a singular fiber of non-simple type to a singular fiber of non-simple type. 
\item[{\rm (3)}]
After replacing $\varphi$ by its iteration, we may view the endomorphism $\beta$ as the identity morphism 
in one of the following cases:
\begin{enumerate}
\item[{\rm (i)}]
$\lkd(C)=1$.
\item[{\rm (ii)}]
$\lkd(C)=0$ and there exists a singular fiber in the fibration $\rho$.
\end{enumerate}
\end{enumerate}
\end{lem}
\Proof
The assertion follows from \cite{GM}. In order to prove the assertion (2), note that there is no 
non-constant morphism from $\A^1$ to $\A^1_*$. Hence it follows that $\varphi$ sends any singular fiber 
of non-simple type to a fiber of non-simple type. Suppose that there exists a singular fiber of non-simple 
type. Let $T$ be the set of points $P$ of $C$ such that $\rho^*(P)$ is a fiber of non-simple type. 
Then it is easy to see that there exists a nonempty subset $T_1$ of $T$ such that $\beta$ induces an 
automorphism on $T_1$. Hence replacing $\varphi$ by its iteration, we may assume that $\beta$ induces the 
identity automorphism on the subset $T_1$.  Suppose that there are only singular fibers of simple type. 
Then the same argument works in this case as well because $\varphi$ sends a multiple fiber to a 
multiple fiber. 

Consider the assertion (3). If $\lkd(C)=1$ then $C$ is a log curve of general type. Hence $\beta$ is an 
automorphism. Furthermore, $\beta$ is of finite order. Suppose  $\lkd(C)=0$ and there exists a singular 
fiber of non-simple type. By the above observation, there exists a point $P \in C$ such that 
$\rho^*(P)$ is a singular fiber of non-simple type and $\beta(P)=P$. Then $\varphi$ induces 
an endomorphism of $C-\{P\}$. Since $\lkd(C-\{P\})=1$ we are done by the case $\lkd(C)=1$. If there 
are only singular fibers of simple type, we are done by the above remark. 
\QED

With the setting before Lemma \ref{Lemma 2.2.5}, we consider the case where $C$ is a rational curve. 
Suppose further that $\rank\Pic(X)=0$ and $\Gamma(X,\SO_X)^*=\C^*$. Then we can specify the structure 
of the $\A^1_*$-fibration $\rho : X \to C$.

\begin{lem}\label{Lemma 2.2.6}
Under the above assumptions the following assertions hold:
\begin{enumerate}
\item[{\rm (1)}]
$C$ is isomorphic to $\BP^1$ or $\A^1$.
\item[{\rm (2)}]
If $\C \cong \BP^1$ then every fiber of $\rho$ is irreducible. Hence every fiber is isomorphic to 
$\A^1_*$ or $\A^1$ if taken with the reduced structure.
\item[{\rm (3)}]
If $\C \cong \A^1$ then every fiber of $\rho$ except for only one fiber $S$ is irreducible and isomorphic 
to either $\A^1_*$ or $\A^1$ if taken with the reduced structure. The unique reducible fiber $S$ consists 
of two irreducible components $F_0, F_1$ such that either $\Gamma= F_0+F_1$ or $\Gamma = F_0$ and 
$\Delta=F_1$. 
\end{enumerate}
\end{lem}
\Proof
(1)\ \ Since $C$ is rational by the hypothesis and since $\Gamma(X,\SO_X)^*=\C^*$, we have $C \cong \BP^1$ 
or $C \cong \A^1$.

\noindent
(2)\ \ The dimension counting of $\Pic(X)_{\Q}$ by making use of the $\A^1_*$-fibration proves the 
assertions (2) and (3). See \cite{MS}.
\QED

\begin{lem}\label{Lemma 2.2.7}
With the same assumptions as in Lemma \ref{Lemma 2.2.6}, suppose further that $C$ is isomorphic to $\BP^1$. 
Let $m_1\Gamma_1, \ldots, m_s\Gamma_s, m_{s+1}\Delta_{s+1}, \ldots, \\
m_{s+t}\Delta_{s+t}$ exhaust all 
singular fibers of $\rho$, where $\Gamma_i \cong \A^1_*\ (1 \le i \le s)$ and $\Delta_j \cong 
\A^1\ (s+1 \le j \le s+t)$. We set $r=s+t$. Let $\varphi$ be an \'etale endomorphism of $X$ such that 
$\rho\cdot\varphi=\beta\cdot\rho$ for an endomorphism $\beta$ of $C$. Suppose $r \ge 3$. Then $\beta$ is 
an automorphism except when the multiplicity sequence is one of the following:
\[
\{m_1,\ldots,m_r\}= \{2,2,2,2\}, \{2,3,6\}, \{2,4,4\},\{3,3,3\}
\]
\end{lem}
\Proof
The proof is the same as in \cite{MM}.
\QED

Let $P_i=\rho(\Gamma_i)$ and $Q_j=\rho(\Delta_j)$ with the above notations. If $\beta$ is an automorphism, 
then $\beta$ induces a permutation on the finite set $\{P_1, \ldots, P_s, Q_{s+1},\\
\ldots, Q_r\}$. Hence, 
by replacing $\varphi$ by its suitable iteration, we may assume that $\beta$ is the identity morphism. 
Namely, we have $\rho\cdot\varphi=\rho$. Then Lemma \ref{Lemma 2.2.2} yields the following result.

\begin{thm}\label{Theorem 2.2.8}
Let $X$ be a smooth affine surface with an $\A^1_*$-fibration $\rho : X \to \BP^1$. Suppose that 
$\lkd(X)=1, \Pic(X)=0$ and $\Gamma(X,\SO_X)^*=\C^*$. Suppose further that there are at least three 
singular fibers and that the multiplicity sequence of the singular fibers of $\rho$ is none of the 
following:
\[
\{m_1,\ldots,m_r\}= \{2,2,2,2\}, \{2,3,6\}, \{2,4,4\},\{3,3,3\}
\]
Then any \'etale endomorphism of $X$ is an automorphism.
\end{thm}
\Proof
Since $\lkd(X)=1$, an \'etale endomorphism $\varphi$ of $X$ preserves the $\A^1_*$-fibration $\rho$. 
Namely there exists an endomorphism $\beta$ of the base curve $\BP^1$ such that 
$\rho\cdot\varphi=\beta\cdot\rho$. Then $\beta$ is an automorphism by Lemma \ref{Lemma 2.2.7}. After 
replacing $\varphi$ by its iteration, we may assume that $\beta$ is the identity morphism. The rest of 
the proof is the same as in Theorem \ref{Theorem 2.2.4}.
\QED

We consider the case of affine surfaces with $\A^1$-fibrations. In \cite{GM2}, the generalized Jacobian 
conjecture for $\Q$-homology planes is considered. It is shown that any \'etale endomorphism of a 
$\Q$-homology plane $X$ is an automorphism if one of the following conditions is satisfied:
\begin{enumerate}
\item[(1)]
$\lkd(X)=2$ or $1$.
\item[(2)]
$\lkd(X)=-\infty$ and $X$ has an $\A^1$-fibration $\rho : X \to B$ with at least two multiple fibers. 
\end{enumerate}
The second case above has one case slipping off. Namely, the surface $X$ in Example \ref{Example 2.1.10} 
is a $\Q$-homology plane with an $\A^1$-fibration which has two multiple fibers of multiplicity $2$. 
It is a counterexample to the generalized Jacobian conjecture. We shall here consider the case (2) 
above. We recall the following two lemmas (cf. \cite[Lemma 6.1]{GM2} and \cite[Lemma 3.1]{GM2,MM}).

\begin{lem}\label{Lemma 2.2.9}
Let $\rho : X \to B$ be an $\A^1$-fibration on a $\Q$-homology plane. Suppose that $\rho$ has 
at least two singular fibers. Let $g : \A^1 \to X$ be a non-constant morphism. Then the image 
of $g$ is a fiber of $\rho$.
\end{lem}
\Proof
In fact, note that the base curve $B$ is isomorphic to $\A^1$ (cf. \cite[Lemma 2.5]{MS}). Then 
the assertion follows from Lemma \ref{Lemma 1.4.16}.
\QED

\begin{lem}\label{Lemma 2.2.10}
For $i=1,2$, let $\rho_i : X_i \to B_i$ be $\A^1$-fibrations on $\Q$-homology planes. Let 
$\varphi : X_1 \to X_2$ and $\beta : B_1 \to B_2$ be dominant morphisms such that $\rho_2\cdot\varphi
=\beta\cdot\rho_1$. Let $m\Gamma$ be an irreducible fiber of $\rho_2$ lying over a point 
$P \in B_2$ with $m \ge 1$ and $\Gamma$ reduced, and let $Q \in B_1$ be a point such that 
$\beta(Q)=P$. Suppose $\rho_1^*(Q)=\ell\Delta$, where $\Delta$ is reduced and irreducible and 
$\ell$ is its multiplicity. Suppose furthermore that $\varphi$ is an \'etale morphism. If the 
ramification index of $\beta$ at $Q$ is $e$ then $\ell e=m$. In particular, if $m=1$ then 
$\ell=e=1$.
\end{lem}

Applying these lemmas, we shall show the following result.

\begin{lem}\label{Lemma 2.2.11}
Let $X$ be a $\Q$-homology plane with an $\A^1$-fibration $\rho : X \to B$. Let 
$m_1A_1, \ldots, m_nA_n$ exhaust all multiple fibers of $\rho$. Let $\varphi : X \to X$ be 
an \'etale endomorphism. Then the following assertions hold:
\begin{enumerate}
\item[{\rm (1)}]
If $n \ge 2$, then there exists an endomorphism $\beta$ of $B$ such that $\rho\cdot\varphi=
\beta\cdot\rho$.
\item[{\rm (2)}]
The above endomorphism $\beta$ is an automorphism provided $n \ge 3$ or $n=2$ and $\{m_1,m_2\}
\ne \{2,2\}$.
\end{enumerate}
\end{lem}

\Proof
The first assertion is an immediate consequence of Lemma \ref{Lemma 2.2.9}. 
So, we consider the second assertion. We employ the arguments in \cite[Lemmas 3.1 and 3.3]{MM}. 
Note that $\beta : B \to B$ is a finite morphism because $B$ is the affine line. By Lemma 
\ref{Lemma 2.2.10}, the set $\{P_1,\ldots,P_n\}$ is mapped to itself by $\beta$, where $P_i=\rho(A_i)$. 
Suppose, furthermore, that the points $Q_1, \ldots, Q_s$, none of which belongs to $\{P_1,\ldots,P_n\}$, 
are mapped to $\{Pp_1,\ldots,P_n\}$. Then, by Lemma \ref{Lemma 2.2.10}, the ramification index of $\beta$ 
at $Q_j$, say $e_j$, is larger than $1$. In fact, if $\beta(Q_j)=P_i$ then $e_j=m_i$.

Since $\beta$ induces an \'etale finite morphism 
\[
\beta : B-\{P_1,\ldots,P_n,Q_1,\ldots,Q_s\} \lto B-\{P_1,\ldots,P_n\},
\]
the comparison of the Euler numbers gives rise to an equality
$$
1-(n+s)=d(1-n), \eqno{(1)}
$$
where $d=\deg\beta$. On the other hand, by summing up the ramification indices, we have an inequality 
$$
2s+n \le dn. \eqno{(2)}
$$
So, by combining (1) and (2) together, we have an inequality 
$$
2(d-1)(n-1) = 2s \le (d-1)n. \eqno{(3)}
$$
Suppose $d > 1$. Then $n \le 2$. Hence, if $n \ge 3$ then $d=1$ and $\beta$ is an automorphism. 
Suppose that $d > 1$ and $n=2$. Then the equality occurs in (3), and hence the equality occurs in (2). 
Namely, the ramification index $e_j$ at $Q_j$ is two for all $j$, and $s=d-1$. Since $d > 1$ implies 
$s > 0$, we may assume that $Q_1$ is mapped to $P_1$. Then $m_1=2$. Suppose $d \ge 3$. Then 
$2s =2(d-1) > d$. Hence one of the $Q_j$ is mapped to $P_2, \ldots, P_n$, say $P_2$. Hence $m_2=2$. 
In this case, after a suitable change of indices, one of the following two cases is possible:
\begin{enumerate}
\item[(1)]
$s=s_1+s_2=d-1$, and $Q_1,\ldots,Q_{s_1}, P_1$ (or $P_2$) (resp. $Q_{s_1+1},\ldots,Q_s, P_2$ 
(or $P_1$) are mapped to $P_1$ (resp. $P_2$).
\item[(2)]
$s=s_1+s_2, d=2s_1=2s_2+2$, and $Q_1,\ldots,Q_{s_1}$ (resp. $Q_{s_1+1},\ldots,Q_s, P_1, P_2$) 
are mapped to $P_1$ (resp. $P_2$).
\end{enumerate}
Finally, suppose that $d=n=2$ and $s=1$. Then we may assume that $\beta(Q_1)=P_1$ and 
$\beta(P_1)=\beta(P_2)=P_2$. Then $m_2=2$ as well by Lemma \ref{Lemma 2.2.10}. So, if 
$\{m_1,m_2\}\ne \{2,2\}$, then $d=1$ and $\beta$ is an automorphism.
\QED

As a consequence of Lemma \ref{Lemma 2.2.11}, we can prove the following result, which rectifies 
Theorem 6.1 in \cite{GM2}.

\begin{thm}\label{Theorem 2.2.12}
Let $X$ be a $\Q$-homology plane with an $\A^1$-fibration $\rho : X \to B$. Let $m_1A_1, \ldots,m_nA_n$ 
exhaust all multiple fibers of $\rho$. Suppose that either $n \ge 3$ or $n=2$ and $\{m_1,m_2\}\ne 
\{2,2\}$. Then any \'etale endomorphism $\varphi : X \to X$ is an automorphism.
\end{thm}

\Proof
By Lemma \ref{Lemma 2.2.11}, there exists an automorphism $\beta$ of $B$ such that 
$\rho\cdot\varphi=\beta\cdot\rho$. Since $\beta$ is an automorphism, Lemma \ref{Lemma 2.2.10} implies that 
$\beta$ induces a permutation of the finite set $\{P_1, \ldots,P_n\}$. By replacing $\beta$ by its 
suitable iteration $\beta^r$, we may assume that $\beta$ induces the identity on $\{P_1,\ldots,P_n\}$. 
Since $n \ge 2$ and $\beta$ ( or rather an induced automorphism of the smooth compactification $\ol{B}$ 
of $B$) fixes the point at infinity $P_\infty$. Hence $\beta$ is then the identity automorphism. 

Let $K=k(B)$ be the function field of $B$ and let $X_K$ be the generic fiber of $\rho$. 
Then $X_K$ is isomorphic to the affine line over $K$, and $\varphi$ induces an \'etale 
endomorphism $\varphi_K$ of $X_K$. Since $\varphi_K$ is then finite, $\varphi_K$ is an automorphism. 
Hence $\varphi$ is birational. Then Zariski's Main Theorem implies that $\varphi$ is an open 
immersion. Note that $\Pic(X)_\Q=0$ and $\Gamma(\SO_X)^*=\C^*$. Suppose that $X\ne \varphi(X)$. 
Then $X-\varphi(X)$ has pure codimension one. Since $\Pic(X)_\Q=0$, there exists a regular 
function $h$ on $X$ such that the zero locus $(h)_0$ of $h$ is supported by $X-\varphi(X)$. 
Then $\varphi^*(h)$ is a non-constant invertible function on $X$, which contradicts the property 
$\Gamma(\SO_X)^*=\C^*$. So, $\varphi$ is an automorphism.
\QED

\section{Normalization by \'etale endomorphisms}

Given an \'etale endomorphism $\varphi : X_1 \to X_2$ with $X_1 \cong X_2 \cong X$, consider the 
normalization $\nu : \wt{X}_2 \to X_2$ of $X_2$ in the function field of $X_1$. If $X$ is normal, $X_1$ 
is a Zariski open set of $\wt{X}_2$ with an open immersion $\iota : X_1 \hookrightarrow \wt{X}_2$ such 
that $\varphi=\nu\circ\iota$. We are interested in the complement $\wt{X}_2-X_1$. First of all, we 
shall generalize Lemma \ref{Lemma 2.2.2}.

\begin{lem}\label{Lemma 2.3.1}
Let $X$ and $Y$ be smooth affine surfaces with $\A^1_*$-fibrations $\rho_X : X \to C$ and 
$\rho_Y : Y \to C$, where $C$ is a smooth algebraic curve. Let $\varphi : Y \to X$ be an \'etale morphism 
such that $\rho_Y=\rho_X\circ\varphi$ and $\codim_X(X-\varphi(Y))\ge 2$. Let $ \nu : Z \to X$ be the 
normalization of $X$ in the function field of $Y$. Then the following assertions hold:
\begin{enumerate}
\item[{\rm (1)}]
$Z$ is a smooth affine surface containing $Y$ as a Zariski open set.
\item[{\rm (2)}]
The normalization morphism $\nu : Z \to X$ is an \'etale Galois covering with a cyclic group $G$ of order 
$n$ as the Galois group, where $n:= \deg\varphi$.
\item[{\rm (3)}]
Let $S=\Gamma+\Delta$ be a singular fiber of $\rho_X$, where $\Gamma$ is $\emptyset$ or isomorphic to 
$\A^1_*$ or $\A^1\cup \A^1$ with two $\A^1$'s meeting in one point transversally and $\Delta$ is a 
disjoint union of the affine lines. If $\Gamma\ne \emptyset$ then $\varphi$ is finite over $\Gamma$, 
{\em i.e.}, $\varphi^*(\Gamma)$ is $G$-invariant, and $\varphi$ is totally decomposable over $\Delta$, 
{\em i.e.}, the stabilizer group of each connected component of $\Delta$ is trivial.
\end{enumerate}
\end{lem}
\Proof
If the $\A^1_*$-fibrations $\rho_X$ and $\rho_Y$ are both untwisted, the proof is exactly the same as 
for Lemma \ref{Lemma 2.2.2}. Let $K$ be the function field of $C$ and let $X_K, Y_K$ be the generic 
fibers of $\rho_X, \rho_Y$, respectively. Let $\ol{X}_K, \ol{Y}_K$ be the smooth completions of 
$X_K, Y_K$, respectively. Then $\ol{X}_K$ and $\ol{Y}_K$ are isomorphic to the projective line 
defined over $K$, and $\rho_X$ (resp. $\rho_Y$) is untwisted or twisted according as $\ol{X}_K$ 
(resp. $\ol{Y}_K$) consists of two $K$-rational points or one non $K$-rational point. Since 
$\varphi$ induces a finite $K$-morphism $\ol{\varphi}_K : \ol{Y}_K \to \ol{X}_K$, $\rho_X$ is untwisted 
if so is $\rho_Y$. 

Consider the case where $\rho_Y$ is twisted. Then there exists a double covering $C' \to C$ such that 
$\rho_{Y'} : Y' \to C'$ is an untwisted $\A^1_*$-fibration, where $Y'$ is the normalization of the fiber 
product $Y\times_CC'$ and $\rho_{Y'}$ is a composite of the normalization morphism $Y' \to Y\times_CC'$ 
and the projection $Y\times_CC' \to C'$. Considering all possible singular fibers of $\rho_Y$, we can 
show that $Y'$ is a smooth affine surface. Similarly, we consider the normalization $X'$ of the fiber 
product $X\times_CC'$ and the induced $\A^1_*$-fibration $\rho_{X'} : X' \to C'$. The \'etale morphism 
$\varphi : Y \to X$ induces an \'etale morphism $\varphi' : Y' \to X'$ such that $\rho_{Y'}=\rho_{X'}
\circ\varphi'$. Furthermore, $\deg\varphi'=\deg\varphi$ and $\codim_{X'}(X'-\varphi'(Y')) \ge 2$. Let 
$Z'$ be the normalization of $X'$ in the function field of $Y'$. It is readily verified that $Z'$ is the 
normalization of the fiber product of $Z\times_CC'$. More precisely, $Z'$ has an involution 
$i : Z' \to Z'$ induced by the involution of the double covering $C' \to C$, and $Z$ is the quotient 
variety with respect to $i$. As shown in the untwisted case, $Z'$ is smooth and the normalization 
morphism $\nu' : Z' \to X'$ is an \'etale Galois covering with a cyclic group $G$ of order $n$ as the 
Galois group, where $n=\deg\varphi'$. Since the involution $i$ commutes with the Galois group action, 
the assertions for $\varphi : Y \to X$ follow from the corresponding assertions for 
$\varphi' : Y' \to X'$.
\QED

Let $\rho : X \to C$ be an $\A^1_*$-fibration on a smooth affine surface $X$. Let $S=\Gamma+\Delta$ be 
a singular fiber of $\rho$. We call $S$ a singular fiber {\em of the first kind}
\index{singular fiber!of the first kind} ({\em of the second kind}
\index{singular fiber!of the second kind}, {\em of the third kind}\index{singular fiber!of the third kind}, 
respectively) if $\Gamma_{\red}\cong \A^1_*$ ($\Gamma_{\red}\cong \A^1\cup \A^1$, 
$\Gamma_{\red} = \emptyset$, respectively). 

\begin{cor}\label{Corollary 2.3.2}
With the same notations and assumptions as in Lemma \ref{Lemma 2.3.1}, let $\varphi : Y \to X$ be an 
\'etale endomorphism  satisfying $\rho_Y=\rho_X\circ\varphi$. Then $\varphi$ is an isomorphism if 
$\rho_X$ has a singular fiber of the second kind.
\end{cor}
\Proof
Suppose that $\rho_X$ has a singular fiber of the second kind $S=\Gamma+\Delta$. Write  
$\Gamma=a_1F_1+a_2F_2$ with $F_1\cong F_2\cong \A^1$. Then $F_1$ and $F_2$ meet each other in a 
single point $P$ transversally. Since $\varphi^*(\Gamma)$ is $G$-invariant, the fiber 
$T:=\rho_Y^{-1}(\rho_X(S))$ is a singular fiber of the second kind, and $\Gamma'=b_1G_1+b_2G_2$ 
with $G_1\cong G_2\cong \A^1$ and $G_1, G_2$ meeting in a single point $Q$ transversally if we write 
$T=\Gamma'+\Delta'$. Furthermore, the group $G$ acts on $G_1\cup G_2$ freely because the morphism 
$\varphi$ induces an \'etale finite morphism from $G_1\cup G_2$ onto $F_1\cup F_2$. Nevertheless, it 
is clear that the point $Q$ is fixed under the $G$-action. This implies that the order $n$ of $G$ 
equals one. Since $n=\deg\varphi$, it follows that $\varphi$ is a birational morphism. Then $Z$ is 
isomorphic to $X$. Hence $\varphi : Y \to X$ is an open immersion. Meanwhile, since $X$ and $Y$ are 
affine, the complement $Z-Y$ is of pure codimension one. Since $\codim_X(X-\varphi(Y))\ge 2$, it 
follows that $\varphi$ is an isomorphism.
\QED

\begin{lem}\label{Lemma 2.3.3}
Let $X$ be a smooth affine surface with $\lkd(X)=1$ and let $\varphi : X \to X$ be an \'etale 
endomorphism, which we write $\varphi : X_1 \to X_2$. Let $\nu : \wt{X}_2 \to X_2$ be the normalization 
of $X_2$ in the function field of $X_1$. Then the following assertions hold:
\begin{enumerate}
\item[{\rm (1)}]
$\wt{X}_2$ is a smooth affine surface with $\lkd(\wt{X}_2)=1$ and $X_1$ is an open set of $\wt{X}_2$. 
Furthermore, $\nu : \wt{X}_2 \to X_2$ is an \'etale Galois covering with a cyclic group $G$ of order 
$n:= \deg\varphi$ as the Galois group.
\item[{\rm (2)}]
The $\A^1_*$-fibration $\rho_1 : X_1 \to C_1$ extends to an $\A^1_*$-fibration $\wt{\rho}_1 : 
\wt{X}_2 \to C_1$ such that $\beta\circ\wt{\rho}_1=\rho_2\circ\nu$, where $\beta : C_1 \to C_2$ is an 
automorphism.
\item[{\rm (3)}]
$\wt{X}_2-X_1$ is a disjoint union of irreducible curves isomorphic to the affine line. The number 
$N$ of irreducible components of $\wt{X}_2-X_1$ is zero or given by the following formula:
\[
N = \sum_{i=1}^r(nd_i-d'_i),
\]
where $r$ is the number of singular fibers of $\rho_2$ (and hence of $\rho_1$) and  $d_i$ (resp. $d'_i$) 
is the number of irreducible components in $\rho_2^*(P_i)$ (resp. $\rho_1^*(\beta^{-1}(P_i))$) isomorphic 
to $\A^1$ with $\{P_1, \ldots, P_r\}$ exhausting all points of $C_2$ such that $\rho_2^*(P_i)$ is a 
singular fiber.
\end{enumerate}
\end{lem}
\Proof
Consider smooth completions $X_1 \hookrightarrow V_1$ and $X_2 \hookrightarrow V_2$ such that 
$D_1:=V_1-X_1$ and $D_2:=V_2-X_2$ are divisors with simple normal crossings and that the endomorphism 
$\varphi$ extends to a morphism $\Psi : V_1 \to V_2$. By the logarithmic ramification formula, we have, 
for every $m > 0$,
\[
|m(D_1+K_{V_1})|=|m\Psi^*(D_2+K_{V_2})|+mR_\Phi
\]
with an effective divisor $R_\Psi$. Let $\Phi_1:= \Phi_{|m(D_1+K_{V_1})|}$ and 
$\Phi_2:=\Phi_{|m(D_2+K_{V_2})|}$ be the morphisms defined by the respective linear 
systems. Since $\lkd(X)=1$, $\Phi_1$ and $\Phi_2$ define respectively the $\A^1_*$-fibrations 
$\rho_1 : X_1 \to C_1$ and $\rho_2 : X_2 \to C_2$ provided $m \gg 0$, where $C_1$ and $C_2$ are 
isomorphic to one and the same smooth curve $C$. Furthermore there exists an automorphism 
$\beta : C_1 \to C_2$ such that $\beta\circ\rho_1=\rho_2\circ\varphi$. Now take $Y$ and $X$ 
in Lemma 3.1 to be $X_1$ and $(C_1,\beta)\times_{C_2}X_2$. We denote by the same symbol 
$\varphi$ the induced \'etale morphism $X_1 \to (C_1,\beta)\times_{C_2}X_2$. Then $Z$ is 
isomorphic to $\wt{X}_2$. The assertions (1) and (2) then follow from Lemma \ref{Lemma 2.3.1} and 
its proof. 

Let $S=\Gamma+\Delta$ be a singular fiber of $\rho_2$ lying over a point $P \in C_2$. Then 
$\wt{\rho}_1^*(Q)=\nu^*(\Gamma)+\nu^*(\Delta)$ is a singular fiber of $\wt{\rho}_1$ such that 
$\nu^*(\Gamma)=\varphi^*(\Gamma)$ is $G$-invariant and $\nu^*(\Delta)$ is totally decomposable, where 
$\beta(Q)=P$. If we write $\rho_1^*(Q)=\Gamma'+\Delta'$ then $\Gamma'=\nu^*(\Gamma)$ and 
$\nu^*(\Delta)$ is the $G$-translate of $\Delta'$. Hence we conclude by a fiberwise argument 
that $\wt{X}_2-X_1$ consists of connected components each of which is isomorphic to $\A^1$. 
Now, by the above argument, we have $N=\sum_{i=1}^r(nd_i-d'_i)$.
\QED

Apart from the generalized Jacobian conjecture, we shall make the following remark (cf. \cite{M3}).

\begin{lem}\label{Lemma 2.3.4}
Let $X$ be a normal affine surface with an $\A^1_*$-fibration $\rho : X \to C$. Let $S$ be a singular 
fiber, i.e., a fiber which is scheme-theoretically not isomorphic to $\A^1_*$. Then $S$ is written as 
$S=\Gamma+\Delta$, where
\begin{enumerate}
\item[{\rm (1)}]
$\Gamma=\emptyset$, $\Gamma$ is an irreducible curve with possibly one cyclic quotient singularity and 
the normalization of $\Gamma$ is isomorphic to $\A^1_*$, or $\Gamma=C_1+C_2$ with the normalizations 
$\wt{C}_1, \wt{C}_2$ being isomorphic to $\A^1$ and $C_1, C_2$ meeting in one point which might be 
a cyclic quotient singular point of $X$;
\item[{\rm (2)}]
$\Delta=\emptyset$ or $\Delta=\coprod_i D_i$ with the normalization of $D_i$ being isomorphic to $\A^1$;
\item[{\rm (3)}]
A connected component $D_i$ of $\Delta$ might have at most one singular point which is a cyclic quotient 
singularity.
\item[{\rm (4)}]
$X$ has at worst cyclic quotient singularities.
\end{enumerate}
\end{lem} 
\Proof
Let $\pi : \wh{X} \to X$ be a resolution of singularity. Then $\wh{\rho}:=  \rho\circ\pi : \wh{X} \to C$ 
is an $\A^1_*$-fibration. Furthermore, $\wh{X}$ has a smooth completion $\wh{X} \hookrightarrow W$ such 
that $W-\wh{X}$ is a divisor of simple normal crossings and the $\A^1_*$-fibration $\wh{\rho}$ extends 
to a $\BP^1$-fibration $p : W \to \ol{C}$, where $\ol{C}$ is a smooth completion of $C$. Then 
$\ol{S}:=p^{-1}(\rho(S))$ is a degenerate fiber of $p$. We may assume that $\ol{S}-S$ does not contain 
$(-1)$ components. Let $T$ be a maximal connected union of irreducible components of $\ol{S}-S$ such that 
$T\cap (W-(\wh{X}\cup \ol{S})=\emptyset$. Since $X$ is affine, $T$ meets at least one irreducible 
component of $S$. If $T$ meets one irreducible component of the proper transform of $S$, then the dual 
graph of $T$ is a linear chain and $T$ contracts to a cyclic quotient singular point. In fact, if the 
dual graph of $T$ has a branching point then the successive contractions of $(-1)$ curves in $\ol{S}$ 
will produce either three components meeting in one point or two components meeting in the same point 
of a cross-section (or a double point of $2$-section). This is a contradiction. A similar argument shows 
that if $T$ meets two irreducible components of $S$ then $T$ again contracts to a cyclic 
quotient singular point. Since $\rho$ is an $\A^1_*$-fibration, $T$ cannot meet three or more irreducible 
components of $S$. These considerations imply the stated assertions. 
\QED

\section{Affine pseudo-coverings}

Let $X$ be a smooth affine variety and let $f : Y \to X$ be a morphism of algebraic varieties. 
We say that $f$ is {\em almost surjective}\index{morphism!almost surjective} if 
$\codim_X(X-f(Y)) \ge 2$ and that $Y$ is an {\em affine pseudo-covering}\index{affine pseudo-covering} 
of $X$ if $Y$ is affine, $f$ is \'etale and almost surjective. Note that 
an affine pseudo-covering is not necessarily an ordinary finite covering. If the 
field extension $k(Y)/k(X)$ is a Galois extension, the affine pseudo-covering is called {\em Galois}
\index{affine pseudo-covering!Galois}. We are interested in determining all affine pseudo-coverings 
of the affine plane. We recall some known results \cite{KM}.

\begin{lem}\label{Lemma 2.4.1}
Let $f : Y \to X$ be an \'etale morphism of smooth affine varieties. Suppose that 
$\Pic(X)=0$ and $\Gamma(Y,\SO_Y)^*=\C^*$. Then $Y$ is an affine pseudo-covering of $X$.
\end{lem}
\Proof
It suffices to show that $f$ is almost surjective. Let $A$ (resp. $B$) be the coordinate 
ring of $X$ (resp. $Y$). Since $f$ is a dominant morphism, we may and shall identify $A$ 
with a subalgebra of $B$ via $f^*$. Then $A$ is factorial and $B^*=k^*$ by the hypothesis. 
Let $Z$ be an irreducible subvariety of $X$ of codimension one. Since $A$ is factorial, 
$Z$ is defined by an element $a \in A$. Let $\gp=aA$ be the defining ideal of $Z$. 
Since $B^*=k^*$, the ideal $aB$ is not equal to $B$. Let
\[
aB=\gq_1\cap \cdots \cap \gq_r
\]
be the primary decomposition of the ideal $aB$. Since $B$ is normal, each primary ideal 
$\gq_i\ (1 \le i \le r)$ has height one. Let $V_i$ be the subvariety in $Y$ defined by 
the radical $\GP_i:=\sqrt{\gq_i}$. Let $W_i$ be the subvariety defined by a prime ideal 
$\GP_i\cap A$. Note that $W_i \subseteq Z$ because $\gp \subseteq \GP_i\cap A$. 
It suffices to show that $W_i=Z$ for every $1 \le i \le r$. Suppose to the contrary 
that $W_i$ is a proper subvariety of $Z$. Then $f|_{V_i} : V_i \to W_i$ is a dominant morphism 
with $\dim V_i > \dim W_i$. Hence a general fiber of $f|_{V_i}$ has positive dimension. 
This is, however, impossible because $f$ is \'etale. 
\QED

This lemma implies that if $X$ is a smooth affine variety with $\Pic(X)=0$ and 
$\Gamma(X,\SO_X)^*=\C^*$, then an \'etale endomorphism $\varphi : X \to X$ is almost surjective. 
This is, in particular, the case for $X=\A^n$. Let $f : Y \to X$ be a Galois affine pseudo-covering 
of smooth affine varieties with Galois group $G:= {\rm Gal}(\C(Y)/\C(X))$ and let $\wt{Y}$ be the 
normalization of $X$ in $\C(Y)$. Then $Y$ is a Zariski open set of $\wt{Y}$ and the group $G$ acts 
freely on $\wt{Y}$ so that $X=\wt{Y}/G$. The purity of branch loci implies that the normalization 
morphism $\wt{f} : \wt{Y} \to X$ is a finite \'etale morphism and that $\wt{Y}$ is thereby smooth.

We shall consider various examples of affine pseudo-coverings mostly in the case where $X$ 
is a smooth affine surface. Let $X:=X(d,r)$ be an affine pseudo-plane. Then the universal covering 
$\wt{X}:=\wt{X}(d,r)$ of $X$ has a Galois group $H(d)\cong\Z/d\Z$. By Lemma \ref{Lemma 1.4.20}, 
$\wt{X}$ contains an open set $U_\omega$ which is isomorphic to $\A^2$ and mapped surjectively onto 
$X$ by the covering mapping, where $\omega \in H(d)$. Hence $\A^2$ is a Galois affine pseudo-covering 
of $X$. Slightly generalizing this result, we can prove the following result.

\begin{lem}\label{Lemma 2.4.2}
Let $X$ be a $\Q$-homology plane with $\lkd(X)=-\infty$. Hence there exists an 
$\A^1$-fibration $\rho : X \to C$, where $C \cong \A^1$. Suppose that $\rho$ has a 
unique multiple fiber $dF$, where $F \cong \A^1$ and $d \ge 2$. Then $\A^2$ is a Galois affine 
pseudo-covering of $X$ with Galois group $H(d)$.
\end{lem}
\Proof
The proof is just a repetition of the previous argument used in the construction of the universal 
covering of an affine pseudo-plane. Let $P=\rho(F)$ and let $\mu : \wt{C} \to C$ be a $d$-ple 
cyclic covering which ramifies totally over $P$ and the point at infinity $P_\infty$ of $C$. Let 
$\wt{Y}$ be the normalization of the fiber product $X\times_C\wt{C}$ and $\wt{f} : \wt{Y} \to X$ 
be the composite of the normalization morphism $\wt{Y} \to X\times_C\wt{C}$ and the first 
projection $X\times_C\wt{C} \to X$. Then $\wt{f}$ is a Galois \'etale covering, and the multiple fiber 
$dF$ splits into a disjoint union of $d$ copies of the affine line on $\wt{Y}$. Let $Y$ be an open 
set of $\wt{Y}$ obtained by removing $d-1$ copies of $\A^1$ from $\wt{f}^{-1}(F)$. Then $Y$ is 
isomorphic to $\A^2$ and mapped surjectively onto $X$ by $\wt{f}$.
\QED

We consider if the converse to the above remark holds. Let $X$ be a smooth affine surface and 
let $f : \A^2 \to X$ be a Galois affine pseudo-covering with Galois group  $G$. Let $Y:= \A^2$ 
and let $\wt{f} : \wt{Y} \to X$ be the normalization of $X$ in the function field of $\A^2$. 
Then $\wt{f}$ is a finite \'etale Galois covering with Galois group $G$ and $\wt{Y}$ contains $Y$ 
as a Zariski open set. Furthermore, $\wt{Y}-Y$ is the disjoint union $\coprod_{i=1}^r C_i$ of the 
curves isomorphic to $\A^1$  \footnote{In fact, consider any $\A^1$-fibration $\rho : Y:=\A^2 \to C$. 
Then $\rho$ extends to an $\A^1$-fibration $\wt{\rho} : \wt{Y} \to \wt{C}$, where 
$\rho=\wt{\rho}\mid_Y$ and $\wt{C}$ contains $C$ as an open set. Then it is readily shown that 
$\wt{Y}-Y$ consists of the fiber components of $\wt{\rho}$ which are necessarily isomorphic to $\A^1$.}.

\begin{lem}\label{Lemma 2.4.3}
Let the notations and assumptions be the same as above. We assume, furthermore, that the following 
hypothesis {\em (H)} holds:
\[
g(C_i)\cap h(C_j)=\emptyset \quad\mbox{whenever $g(C_i)\ne h(C_j)$ for 
$1 \le i, j \le r$ and $g, h \in G$.}
\]
Then there exists an $\A^1$-fibration $\wt{\rho} : \wt{Y} \to \wt{C}$ such that the following 
conditions are satisfied.
\begin{enumerate}
\item[{\rm (1)}]
The group $G$ preserves the $\A^1$-fibration $\wt{\rho}$. Namely, for any fiber component $A$ 
of $\wt{\rho}$, the translate $g(A)$ is a fiber component of $\wt{\rho}$ for $g \in G$.
\item[{\rm (2)}]
$g(C_i)$ is a fiber component of $\wt{\rho}$ for $1 \le i \le r$ and $g \in G$.
\item[{\rm (3)}]
The $\A^1$-fibration $\wt{\rho}$ restricts to an $\A^1$-fibration $\rho : Y \to C$, where $C$ is a 
Zariski open set of $\wt{C}$ and isomorphic to $\A^1$ and where $\rho$ is surjective.
\item[{\rm (4)}]
The curve $\wt{C}$ is isomorphic to either $\BP^1$ or $\A^1$.
\item[{\rm (5)}]
There exists a group homomorphism $\alpha : G \to \Aut\wt{C}$ which is injective. Let $\wt{K}$ be the 
function field of $\wt{C}$. Then $G$ is the Galois group of the field extension $\wt{K}/K$, where 
$K=\wt{K}^G$.
\end{enumerate}
\end{lem}
\Proof
Since $f : Y \to X$ is almost surjective, there exists a translate $g(C_i)$ for some $g \in G$ such 
that $g(C_i)\cap Y \ne \emptyset$. Then $g(C_i) \subset Y$. In fact, if not, $g(C_i)\cap C_j\ne 
\emptyset$ for some $j$. This contradicts the hypothesis. Since $g(C_i)$ is isomorphic to $\A^1$, 
Theorem of Abhyankar-Moh-Suzuki implies that there exists a trivial $\A^1$-bundle structure $\rho : 
Y \to C$ with $g(C_i)$ as a fiber, where $C \cong \A^1$. Note that if $h(C_j)\cap Y \ne \emptyset$ 
for $h \in G$, then $h(C_j)$ is a fiber of $\rho$. \footnote{We may choose coordinates $\{x,y\}$ on $Y$ 
so that the curve $g(C_i)$ is defined by $y=0$. Suppose that $h(C_j)$ is defined by $f=0$. Then $f$ 
as well as $f-\alpha$ is irreducible for all $\alpha \in k$ by the theorem of Abhyankar-Moh-Suzuki. 
Suppose $h(C_j)\cap g(C_i)=\emptyset$. Then we have $f(x,0)=c \in k$. Namely, $f-c=yg(x,y)$ with 
$g(x,y) \in k[x,y]$. Since $f-c$ is irreducible, we must have $g(x,y) \in k^*$.}  
Let $F$ be a fiber which is not of the form 
$h(C_j)$ with $1 \le j \le r$ and $h \in G$. Then, for any $g \in G$, the translate $g(F)$ meets 
none of the $h(C_j)$. In fact, if $g(F)\cap h(C_j)\ne \emptyset$, then $F\cap g^{-1}h(C_j) \ne 
\emptyset$, which implies that $F=g^{-1}h(C_j)$ because both are the fibers of $\rho$. This 
is a contradiction to the choice of $F$. It is now clear that $\rho : Y \to C$ extends to an 
$\A^1$-fibration on $\wt{\rho} : \wt{Y} \to \wt{C}$, where the parameter space $\wt{C}$ for $\wt{\rho}$ 
contains $C$ as a Zariski open set. 

We shall show the assertion (5). Since $G$ preserves the $\A^1$-fibration $\wt{\rho} : \wt{Y} \to 
\wt{C}$, there exists a natural group homomorphism $\alpha : G \to \Aut\wt{C}$. We shall show that 
$\alpha$ is injective. Since $\rho : Y \to C$ is a trivial $\A^1$-bundle, the generic fiber $Y_{\wt{K}}$ 
of $\rho$ is the affine line $\A^1_{\wt{K}}:=\Spec \wt{K}[x]$, where $\wt{K}$ is the function field 
of $\wt{C}$. Then we have 
\[
g(x)=a(g)x+b(g) \quad \mbox{with}\quad a(g) \in \wt{K}^*, b(g) \in \wt{K} \quad \mbox{for}\quad g \in G,
\]
where 
\[
a(gh)=a(h)^ga(g) \quad \mbox{and}\quad b(gh)= a(h)^gb(g)+b(h)^g \quad \mbox{for}\quad g, h \in G.
\]
Suppose that $\Ker\alpha \ne (0)$. Then it follows by induction that 
\[
b(g^m)=\left(a(g)^{m-1}+\cdots+a(g)+1\right)b(g) \quad \mbox{for $g \in \Ker\alpha$},
\]
where $b(1)=0$. If $a(g)\ne 1$, then the point given by $x=-b(g)/(a(g)-1)$ is left fixed by the action 
of $g$, but this is not possible because $G$ acts on $\wt{Y}$ freely. If $a(g)=1$ and $m$ is 
the order of $g$, then $b(g)=0$ because $a(g)^{m-1}+\cdots+a(g)+1=m\ne 0$ and $b(g^m)=0$. Then $g$ 
acts on $\wt{Y}$ trivially, and this is a contradiction. Hence $\Ker\alpha=(1)$ and $\alpha$ 
is injective.
\QED

\begin{thm}\label{Theorem 2.4.4}
Let $X$ be a smooth affine surface and let $f : \A^2 \to X$ be a Galois affine 
pseudo-covering with Galois group $G$. Assume that the hypothesis {\em (H)} in Lemma \ref{Lemma 2.4.3} 
is satisfied. Then the following assertions hold.
\begin{enumerate}
\item[{\rm (1)}]
There exist an $\A^1$-fibration $\rho : X \to C$ and a finite $G$-morphism of smooth curves 
$\pi : \wt{C} \to C$ such that the given morphism $f : \A^2 \to X$ factors as 
$f=\wt{f}\circ\iota$, where $\wt{f} : \wt{Y} \to X$ is 
a composite of the normalization morphism $\wt{Y} \to X\times_C\wt{C}$ with the normalization 
$\wt{Y}$ of $X\times_C\wt{C}$ in the function field of $\A^2$ and the canonical 
projection $X\times_C\wt{C} \to X$ and where $\iota : \A^2 \to \wt{Y}$ is an open immersion. 
\item[{\rm (2)}]
$C$ is isomorphic to either $\A^1$ or $\BP^1$.
\item[{\rm (3)}]
If $C$ is isomorphic to $\A^1$ and every fiber of $\rho$ is irreducible, then $X$ is a 
$\Q$-homology plane and the Galois group is a cyclic group. Furthermore, the $\A^1$-fibration 
$\rho : X \to C$ has a single multiple fiber.
\item[{\rm (4)}]
If $C$ is isomorphic to $\BP^1$ and every fiber of $\rho$ is irreducible, then the 
$\A^1$-fibration $\rho : X \to C$ has at most three multiple fibers. If there are 
three multiple fibers $m_iF_i$ with $m_i > 1$ and $F_i \cong \A^1$, then $\{m_1,m_2,m_3\}$ 
is $\{2,2,n\}, \{2,3,3\}, \{2,3,4\}$ or $\{2,3,5\}$ up to permutations.
\end{enumerate}
\end{thm}
\Proof
Our proof consists of several steps. We may assume that the Galois group $G$ is non-trivial. 
If $G$ is trivial, the morphism $f$ is an isomorphism.
\svskip

\noindent
(I)\ Let $Y:=\A^2$ and let $\wt{Y}$ be the normalization of $X$ in the function field of $Y$. 
We follow the proof of Lemma \ref{Lemma 2.4.3}. By the construction of the $\A^1$-fibration 
$\rho : Y \to C$, we find a $G$-stable open set $\wt{W}$ of $Y$ such that $Y-\wt{W}$ 
is a union of finitely many fibers of the $\A^1$-bundle and $\wt{W}= \wt{U} \times \A^1$. 
Note that $\wt{U}$ is an open set of the affine line. Hence we may write 
$\wt{U}= \Spec k[t,u(t)^{-1}]$, where $u(t) \in k[t]$ and where we may assume that 
$t \mid u(t)$. Then $\wt{W}= \Spec A$ with $A=k[t,u(t)^{-1},x]$ and $A^*$ is generated by 
$k^*$ and distinct prime factors of $u(t)$. Since $A$ is a $G$-algebra over $k$, the 
$G$-action preserves $A^*$. Hence $G$ acts on the function field $\wt{K}:=k(t)$. Let $K=\wt{K}^G$ and 
$A_{\wt{K}}=A\otimes_{k[t]}k(t)=\wt{K}[x]$. Then $G$ is the Galois group of the field extension 
$\wt{K}/K$ and $G$ acts on $A_{\wt{K}}$. Hence, as in the proof of Lemma \ref{Lemma 2.4.3}, we may write
\[
g(x)=a(g)x+b(g) \quad \mbox{with}\quad a(g), b(g) \in K \quad \mbox{for}\quad g \in G,
\]
where 
\[
a(gh)=a(h)^ga(g) \quad \mbox{and}\quad b(gh)= a(h)^gb(g)+b(h)^g \quad \mbox{for}\quad g, h \in G.
\]
By Theorem 90 of Hilbert, we find $c \in \wt{K}^*$ such that $a(g)=c(c^g)^{-1}$ for $g \in G$. 
By replacing $x$ by $cx$, we may then assume that $a(g)=1$. Then $b(gh)=b(g)+b(h)^g$ for $g, h \in G$. 
Hence $b(g)=d-d^g$ with $d=\left(\sum_{h\in G}b(h)\right)/|G|$. By replacing $x$ by $x+d$, 
we may assume that $g(x)=x$ for $g \in G$. Namely, a suitable choice of the fiber coordinate $x$ 
on the generic fiber (hence for all fibers over a sufficiently small open set of $\wt{U}$), 
we may assume that the group $G$ acts trivially on the fibers. We may assume that the curve 
$\wt{U}$ is chosen in such a way that $G$ acts trivially on the factor $\A^1$ in the product 
$\wt{W}=\wt{U}\times \A^1$.
\svskip

\noindent
(II)\ Let $U$ be the quotient $\wt{U}/G$. Then $W:=U\times \A^1$ is an open set of $X$, and 
$\wt{W}=W\times_U\wt{U}$. The morphism $f|_{\wt{W}} : \wt{W} \to W$ is induced by the 
quotient morphism $q : \wt{U} \to U$. Let $\wt{B}$ and $B$ be the smooth completions of 
$\wt{U}$ and $U$, respectively. Then $\wt{B}$ and $B$ are isomorphic to $\BP^1$. The morphism 
$q$ extends to a morphism $\pi : \wt{B} \to B$. The $G$-action on $\wt{U}$ extends naturally 
to $\wt{B}$ and $\pi$ is the quotient morphism $\wt{B} \to \wt{B}/G$, where $B=\wt{B}/G$. 
On the other hand, the projection $p_1 : W \to U$ extends to a morphism $\rho : X \to C$, 
where $C$ is an open set of the curve $B$. Since $\rho\circ f : \A^2 \to C$ is a dominant 
morphism, $C$ is isomorphic to $\A^1$ or $\BP^1$. Let $\wt{C}$ be the inverse image 
$\pi^{-1}(C)$ in the curve $\wt{B}$. Then $\pi : \wt{C} \to C$ is a finite (possibly ramified) 
$G$-morphism with $C=\wt{C}/G$. By the above construction, it follows that $f : \A^2 \to X$ 
factors as $\A^2 \st{\sigma} X\times_{C}\wt{C} \st{p_X} X$, where $\sigma$ is a birational 
morphism and the canonical projection $p_X$ is a finite morphism. Hence $\wt{Y}$ is the 
normalization of $X\times_C\wt{C}$ in the function field of $Y=\A^2$, and $f : Y \to X$ 
factors as indicated in the statement, where $\iota$ is an open immersion by Zariski Main 
Theorem. Note by the above construction that $\wt{Y}$ has an $\A^1$-fibration 
$\wt{\rho} : \wt{Y} \to \wt{C}$ such that the restriction 
$\wt{\rho}\mid_Y : Y \to \wt{\rho}(Y)$ is an $\A^1$-bundle over $\A^1$. 
Thus the first and second assertions are verified.
\svskip

\noindent
(III)\ Suppose that $C$ is isomorphic to $\A^1$ and every fiber of $\rho$ is irreducible. 
Then $X$ is a $\Q$-homology plane because $\Pic(X)$ is a finite torsion group. Since $G$ is a finite 
subgroup of $\Aut\A^1$, $G$ is a cyclic group. We shall show that the $\A^1$-fibration $\rho : 
X \to C$ has a single multiple fiber. If there are no multiple fibers, then $\rho : X \to C$ is a 
trivial $\A^1$-bundle over $\A^1$. Hence $X$ is isomorphic to $\A^2$. Then $\wt{f} : \wt{Y} \to 
X$ is an isomorphism. This case was excluded at the beginning. Let $m_iF_i\ (1 \le i \le s)$ exhaust 
all singular fibers of $\rho$. Let $P_i=\rho(F_i)$ and let $P_\infty=B-C$, where $B=\wt{B}/G 
\cong \BP^1$. By the construction of $\wt{Y}$ in the assertion (2) (see also Lemma 
\ref{Lemma 2.4.5} below), we know that 
\begin{enumerate}
\item[(i)]
The Galois covering $\pi : \wt{B} \to B$ ramifies at the points $\wt{\pi}^{-1}(P_i)$ with 
ramification index $N/m_i$ for $1 \le i \le s$ and totally ramifies over the point $P_\infty$;
\item[(ii)]
The curve $\wt{C}$ is isomorphic to $\A^1$ and $\wt{\rho} : \wt{Y} \to \wt{C}$ has $m_i$ connected 
components over every point of $\wt{\pi}^{-1}(P_i)$ for $1 \le i \le s$, each of which is isomorphic 
to $\A^1$;
\item[(iii)]
In order to obtain $Y$, we remove $m_i-1$ components from the fiber of $\wt{\rho}$ over every point 
of $\wt{\pi}^{-1}(P_i)$ for $1 \le i \le s$.
\end{enumerate}
Let $N$ be the order of $G$. By the Riemann-Hurwitz theorem, we have
\[
-2=-2N+\sum_{i=1}^s\frac{N}{m_i}(m_i-1)+(N-1),
\]
from which we obtain 
\[
\frac{1}{m_1}+\cdots+\frac{1}{m_s}=(s-1)+\frac{1}{N}.
\]
This implies that $s=1$ and $m_1=N$.
\svskip

\noindent
(IV)\  We shall consider the assertion (4). Suppose first that there are no multiple fibers of $\rho$. 
Hence $\rho : X \to C$ is an $\A^1$-bundle over $C \cong \BP^1$. Hence $\wt{Y} \cong X\times_C\wt{C}$ 
and $\wt{Y} \to \wt{C}$ is an $\A^1$-bundle over $\wt{C} \cong \BP^1$ as well. Since $\wt{\rho}\mid_Y$ 
is an $\A^1$-bundle over $\A^1$, only one fiber of $\wt{\rho}$ is removed to obtain $Y$. 
On the other hand, since $\pi : \wt{C} \to C$ is a finite $G$-morphism with non-trivial $G$, 
the morphism $f : Y \to X$ is not an \'etale morphism. This is a contradiction. So, $\rho$ has at 
least one multiple fiber. Let $m_iF_i\ (1 \le i \le s)$ exhaust all the multiple fibers of $\rho$, 
where $C \cong \BP^1$ and let $P_i=\rho(F_i)$. Suppose $s > 2$. Let $g : D \to C$ be a 
Galois covering which ramifies over the points $P_i$ with ramification index $m_i$. 
Such a covering exists by \cite{BN} and \cite{Fox}. Let $Z$ be the normalization of 
$X\times_CD$ and let $h : Z \to X$ be a composite of the normalization morphism 
$Z \to X \times_CD$ and the projection $X\times_CD \to X$. Then $h$ is an \'etale finite 
morphism. Since $\wt{f} : \wt{Y} \to X$ is the universal covering morphism, there exists 
a finite \'etale morphism $\tau : \wt{Y} \to Z$ such that $\wt{f}=h\circ\tau$. Hence 
the curve $D$ is a rational curve. This implies that 
\[
2N-2=\sum_{i=1}^s\frac{N}{m_i}(m_i-1).
\]
It then follows that $s=3$ and $\{m_1,m_2,m_3\}$ is $\{2,2,n\}, \{2,3,3\}, \{2,3,4\}$ or 
$\{2,3,5\}$ up to permutations. Hence the assertion (4) holds. 
\QED

Let $X$ be a smooth affine variety. We define the {\em Makar-Limanov invariant}
\index{Makar-Limanov invariant} $\ML(X)$ to be the intersection of the invariant subrings of the 
coordinate ring of $X$ under all possible $G_a$-actions 
on $X$. Let $X$ be a $\Q$-homology plane having the properties in the assertion (3) of Theorem 
\ref{Theorem 2.4.4}. The Makar-Limanov invariant $\ML(X)$ is either a polynomial ring $k[x]$ or the 
constant field $k$. By virtue of Theorem \ref{Theorem 1.4.10}, the affine pseudo-plane $X:=X(d,r)$ 
has the Makar-Limanov invariant $k[x]$ provided $r \ge 2$. On the other hand, we call a smooth affine 
surface $X$ a {\em Platonic $\A^1$-fiber space}\index{Platonic $\A^1$-fiber space} if $X$ has an 
$\A^1$-fibration having the properties in the assertion (4) of Theorem \ref{Theorem 2.4.4}. Let $X$ be 
a smooth affine surface with $\A^2$ as a Galois affine pseudo-covering. Lemma \ref{Lemma 2.4.5} below 
implies that the Picard rank $\dim\Pic(X)_\Q$ can be arbitrarily big.

\begin{lem}\label{Lemma 2.4.5}
Let $\rho : X \to C$ be an $\A^1$-fibration with $C \cong \A^1$ with one multiple 
irreducible fiber and $r$ reducible fibers with respective multiplicity one. Write 
$\rho^*(P_0)=mF_0$ with $m > 1$ and $\rho^*(P_i)=\sum_{j=1}^{n_i}F_{ij}$ with connected 
components $F_{ij}$ for $1 \le i \le r$ and $n_i > 1$, where $F_0$ and $F_{ij}$ are 
isomorphic to $\A^1$. Suppose $m \ge n_i$ for all $i$. Let $\pi : \wt{C} \to C$ be 
a cyclic Galois covering of degree $m$ ramifying totally 
over $P_0$ and the point at infinity $P_\infty$. Let $\wt{Y}$ be the normalization of 
$X\times_C\wt{C}$ and let $\wt{\rho} : \wt{Y} \to \wt{C}$ be the composite of the 
normalization morphism and the canonical projection, where $\wt{C}\cong \A^1$. 
Let $Q_0=\pi^{-1}(P_0)$ and let $\pi^{-1}(P_i)= \{Q_i^{(1)},\ldots,Q_i^{(m)}\}$. 
Then $\wt{\rho}^{-1}(Q_0)$ consists of $m$ affine lines with multiplicity 
one and $\wt{\rho}^{-1}(Q_i^{(\ell)})$ has also $n_i$ disjoint affine lines 
$\sum_{j=1}^{n_i}\wt{F}_{ij}^{(\ell)}$ for $1 \le i \le r$ and $1 \le \ell \le m$. 
Leave only one component among the $m$ affine lines in $\wt{\rho}^{-1}(Q_0)$ and leave 
one component, say $G_i^{(\ell)}$, in $\sum_{j=1}^{n_i}\wt{F}_{ij}^{(\ell)}$ which is chosen 
so that the images of $G_i^{(\ell)}\ (1 \le \ell \le m)$ by $\wt{f} : \wt{Y} \to X$ cover 
$\sum_{j=1}^{n_i}F_{ij}$ for each $1 \le i \le r$. Let $Y$ be the open set of $\wt{Y}$ with 
only one irreducible component left in $\wt{\rho}^{-1}(Q_0)$ and 
$\wt{\rho}^{-1}(Q_i^{(\ell)})$ for $1 \le i \le r$ and $1 \le \ell \le m$. Then $Y$ is 
isomorphic to $\A^2$ and the restriction $f=\wt{f}\mid_Y : Y \to X$ is a Galois \'etale 
pseudo-covering. It is easy to compute $\dim\Pic(X)_\Q= \sum_{i=1}^r(n_i-1)$. Hence 
$\dim\Pic(X)_\Q$ can be arbitrarily big.
\end{lem}
\Proof
The proof should be clear from the construction.
\QED

It is almost clear that a Galois affine pseudo-covering of a topologically simply-connected smooth 
algebraic variety is trivial, that is to say, the covering morphism is an isomorphism.
\footnote{Let $f : Y \to X$ be a Galois affine pseudo-covering, where $\pi_1(X)=(1)$. Let $\wt{f} : 
\wt{Y} \to X$ be the induced Galois \'etale finite covering, where $Y$ is an open set of $\wt{Y}$ and 
$f$ is the restriction of $\wt{f}$ onto $Y$. Then $\wt{f}$ is an isomorphism because $\pi_1(X)=(1)$. 
Hence $f : Y \to X$ is an open immersion such that $\codim_X(X-f(Y))\ge 2$. Since $Y$ is affine, it 
follows that $X=f(Y)$.} There are, 
nonetheless, non-trivial affine pseudo-coverings of the affine space. We just give few examples in 
the case of $\A^1$. It is obvious that we can produce various examples of affine pseudo-coverings of 
$\A^n$ by taking direct products.
\svskip

\noindent
{\bf Example 2.5.1}\label{Example 2.4.1}\ \ Let $C$ be a smooth cubic curve in $\BP^2$. 
Let $P$ be a flex of $C$ and 
let $\ell_P$ be the tangent line of $C$ at $P$. Choose a point $Q$ on $\ell$ such that $Q\ne P$ and 
that the other tangent lines of $8$ other flexes do not meet $\ell_P$ on $Q$. Consider the 
projection of $\BP^2$ to $\BP^1$ with center $Q$. Since every line $\ell$ through $Q$ other than $\ell_P$ 
meets the curve $C$ either in three distinct points or in two distinct points with intersection 
multiplicities $1$ and $2$, we obtain an open set $Y$ of $C$ by omitting the point $P$ and the points 
where the line through $Q$ meets $C$ with intersection multiplicity $2$ and a surjective \'etale morphism 
$f : Y \to \A^1$ by restricting the projection onto $Y$. This morphism has degree $3$ and is non-Galois.
\svskip

\noindent
{\bf Example 2.5.2}\label{Example 2.4.2}\ \ By the same construction as above, we can assume that $Y$ 
is rational. In fact, start with a cuspidal cubic $X_0^2X_2=X_1^3$ and take the point $Q:=(0,0,1)$ 
as the center of the projection. The lines through $Q$ are given by $\{a(X_0-X_1)=X_2 \mid a \in k\cup 
(\infty)\}$. Let $Y$ be the curve $C$ with the cusp, the point $(0,0,1)$ and one more point removed off 
and let $f$ be the restriction of the projection from $Q$ onto $Y$. Then $f : Y \to \A^1$ is a surjective 
\'etale morphism.
\svskip

The following result shows that the surfaces appearing in Theorem \ref{Theorem 2.4.4} cannot be 
affine pseudo-coverinigs of $\A^2$.

\begin{thm}\label{Theorem 2.4.6}
Let $X$ be either a $\Q$-homology plane with $\ML(X)=k$ or a Platonic $\A^1$-fiber space possibly except 
for the case where the multiplicity sequence is $\{2,2,m\}$ with $m > 5$. Then there are no \'etale 
morphisms from $X$ to $\A^2$.
\end{thm}
\Proof
We first treat the case where $X$ is a $\Q$-homology plane with $\ML(X)=k$. Let $f : X \to \A^2$ be 
an \'etale morphism which we do not have to assume to be almost surjective. Let $p : S \to X$ be 
the universal covering, where $S$ is an affine hypersurface $xy=z^n-1$ in $\A^3$. We refer to \cite{KM2} 
for relevant results. Let $\varphi=f\cdot p$. Then 
$\varphi : S \to \A^2$ is an \'etale morphism. We show that such $\varphi$ does not exist. Note that 
the canonical divisor $K_S$ is trivial and $\Gamma(\SO_S)^*=k^*$. By a straightforward computation, 
one can show that a differential $2$-form $\omega=(1/z^{n-1})dx\wedge dy$ is an everywhere nonzero 
regular form. On the other hand, since $\varphi$ is \'etale, there exist elements $f,g$ of 
$\Gamma(\SO_S)$ such that $\varphi^*(df\wedge dg)$ is an everywhere nonzero regular $2$-form. Hence 
we have $\omega=c\varphi^*(df\wedge dg)$. Meanwhile, there exists a free $\Z/n\Z$-action given by 
$(x,y,z) \mapsto (\zeta x,\zeta^{-1}y, \zeta^dz)$, where $\zeta$ is a primitive $n$-th root of unity 
with $n=|G|$ and $0 < d < n$. Then $\zeta(\omega)= \zeta^{-d(n-1)}\omega$ and $\varphi^*(df\wedge dg)$ 
is invariant under this $\Z/n\Z$-action because $\varphi$ splits via $X$. So, $d(n-1)$ is divisible 
by $n$, which is not the case.

Next consider the case where $X$ is a Platonic $\A^1$-fiber space. We consider, in general, the case 
where $X$ has an $\A^1$-fibration $\rho : X \to \BP^1$ whose fibers are all irreducible. Let $m_1F_1, 
\ldots, m_rF_r$ exhaust all multiple fibers with $m_i > 1$ and $F_i \cong \A^1$. As in the first case,
the assertion holds provided $K_X \not\cong \SO_X$. We consider a smooth compactification $V$ of $X$ 
such that $V$ is obtained by a sequence $\sigma$ of blowing-ups from a minimal ruled surface $\Sigma_n$ 
with $n \ge 0$ and the closures $\ol{F}_i$ of $F_i$ are unique $(-1)$ components of the singular fibers 
$\Phi_i$. Let $T$ be the image under $\sigma$ of the unique boundary component which is transversal to 
the $\BP^1$-fibration on $V$ which extends the given $\A^1$-fibration $\rho$. Let $M$ be a minimal 
section of $\Sigma_n$ and write $T \sim M+a\ell$ with a fiber $\ell$. Let $k_i$ be the coefficient of 
$\ol{F}_i$ in $K_V$ which is written as 
\[
K_V \sim -2M-(n+2)\ell+\sum_{i=1}^r(k_i\ol{F}_i+\cdots).
\]
Since $T \sim M+a\ell$, we have 
$$
K_X \sim \left\{(2a-n-2)+\sum_{i=1}^r\frac{k_i}{m_i}\right\}\ell, \eqno{(1)}
$$
where $\ell$ signifies abusedly the restriction of $\ell$ on $X$. The dual graph of $\Phi_i$ 
together with the proper transform $T'$ of $T$ looks like:

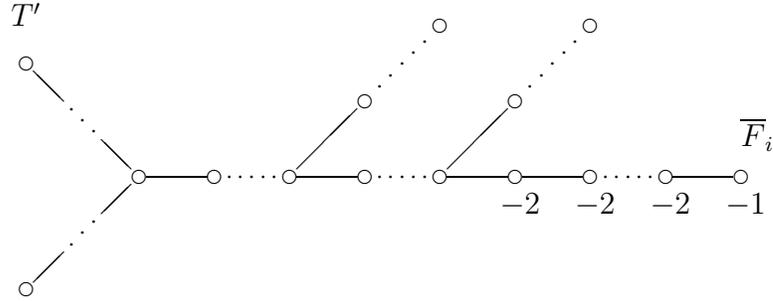
\raisebox{-50mm}{
\begin{picture}(115,50)(0,0)
\unitlength=1mm
\put(10,40){\circle{1.8}}
\put(11,39){\line(1,-1){4}}
\multiput(16,34)(1.5,-1.5){3}{\circle*{0.2}}
\put(20,30){\line(1,-1){4}}
\put(25,25){\circle{1.8}}
\put(10,10){\circle{1.8}}
\put(11,11){\line(1,1){4}}
\multiput(16,16)(1.5,1.5){3}{\circle*{0.2}}
\put(20,20){\line(1,1){4}}
\put(26,25){\line(1,0){8}}
\put(35,25){\circle{1.8}}
\multiput(37,25)(1.5,0){5}{\circle*{0.2}}
\put(45,25){\circle{1.8}}
\put(46,26){\line(1,1){8}}
\put(55,35){\circle{1.8}}
\multiput(57,37)(1.5,1.5){5}{\circle*{0.2}}
\put(65,45){\circle{1.8}}
\put(46,25){\line(1,0){8}}
\put(55,25){\circle{1.8}}
\multiput(57,25)(1.5,0){5}{\circle*{0.2}}
\put(65,25){\circle{1.8}}
\put(66,26){\line(1,1){8}}
\put(75,35){\circle{1.8}}
\multiput(77,37)(1.5,1.5){5}{\circle*{0.2}}
\put(85,45){\circle{1.8}}
\put(66,25){\line(1,0){8}}
\put(75,25){\circle{1.8}}
\put(76,25){\line(1,0){8}}
\put(85,25){\circle{1.8}}
\multiput(87,25)(1.5,0){5}{\circle*{0.2}}
\put(95,25){\circle{1.8}}
\put(96,25){\line(1,0){8}}
\put(105,25){\circle{1.8}}
\put(8,45){$T'$}
\put(73,20){$-2$}
\put(83,20){$-2$}
\put(93,20){$-2$}
\put(103,20){$-1$}
\put(105,29){$\ol{F}_i$}
\end{picture}}

Let $h_i$ be the number of blowing-ups to be needed to produce this dual graph from the single 
fiber $\ell_i$ of $\Sigma_n$. The center of the first blowing-up is called the {\em starting point}
\index{blowing-up!starting point}. We may assume that the starting point on the fiber $\ell$ does not 
lie on the section $T$. We omit the suffix $i$ to denote $m_i, k_i$ and $h_i$. If $m$ is small, 
we can express $k$ in terms of $h$. 
Suppose $m=2$. Then $k=h$. Suppose $m=3$. Then $k=h+1$ if $(-3)\cap T'= \emptyset$ and $k=h$ if 
$(-3)\cap T'\ne \emptyset$, where $T'$ is the proper transform of $T$ on $V$ and $(-3)\cap T'=\emptyset$ 
(resp $\ne \emptyset$) means that the adjacent component of $T'$ in the graph is not (resp. is) a 
$(-3)$-component. Suppose $m=4$. Then $k=h+s-1$ with $s \ge 3$ and $h \ge s+2$ if the dual graph of 
$\Phi$ has more than one branching vertices, where the vertices are connected to more than two adjacent 
vertices, which is necessarily as in Figure 1 below with the proper transform $\ell'$ of $\ell$.

\raisebox{-40mm}{
\begin{picture}(125,40)(0,-10)
\unitlength=1mm
\put(10,30){\circle{1.8}}
\put(10,10){\circle{1.8}}
\put(20,20){\circle{1.8}}
\put(30,20){\circle{1.8}}
\put(60,20){\circle{1.8}}
\put(70,20){\circle{1.8}}
\put(80,30){\circle{1.8}}
\put(80,10){\circle{1.8}}
\put(110,10){\circle{1.8}}
\put(120,10){\circle{1.8}}
\put(11,29){\line(1,-1){8}}
\put(11,11){\line(1,1){8}}
\put(21,20){\line(1,0){8}}
\put(31,20){\line(1,0){9}}
\multiput(41,20)(1.5,0){5}{\circle*{0.2}}
\put(50,20){\line(1,0){9}}
\put(61,20){\line(1,0){8}}
\put(71,21){\line(1,1){8}}
\put(71,19){\line(1,-1){8}}
\put(81,10){\line(1,0){9}}
\multiput(91,10)(1.5,0){5}{\circle*{0.2}}
\put(100,10){\line(1,0){9}}
\put(111,10){\line(1,0){8}}
\put(5,35){$-2$}
\put(13,30){$\ell'$}
\put(5,5){$-2$}
\put(18,15){$-2$}
\put(28,15){$-2$}
\put(58,15){$-3$}
\put(67,15){$-2$}
\put(78,35){$-2$}
\put(78,15){$-2$}
\put(108,15){$-2$}
\put(118,15){$-1$}
\put(118,5){$E_h$}
\put(20,10){$\underbrace{\hspace{50mm}}_{s-2}$}
\put(55,-3){Figure 1}
\end{picture}}
\svskip

\noindent
If the dual graph of $\Phi$ has only one branching vertex and $(-4)\cap T' =\emptyset$ 
(resp. $\not= \emptyset$), then $k=h+2$ (resp. $k=h$). Suppose $m=5$. Then $k=h, k=h+2$ or $k=h+3$, 
resp. if $(-5)\cap T'\ne \emptyset, (-3)\cap T'\ne\emptyset$, or $(-2)\cap T' \ne \emptyset$, resp. 
There are two cases when $(-2)\cap T' \ne \emptyset$. Namely, the proper transform $\ell'$ meets 
either a $(-2)$-component of $(-3)$-component. On the other hand, it is rather easy to show by 
induction on the number of blowing-ups that $k > m$. 

We shall prove the assertion of Theorem \ref{Theorem 2.4.6} in the second case. Recall that 
$T \sim M+a\ell$, where $a \ge n$ or $a=0$. Suppose $T\ne M$. Since $X$ is a Platonic $\A^1$-fiber space, 
we have $r=3$. Hence $I > 0$ in the relation $(1)$, where 
\[
I=\left\{(2a-n-2)+\sum_{i=1}^3\frac{k_i}{m_i}\right\}.
\]
In fact, $k_i/m_i > 1$ and $I > (2a-n-2)+3 \ge (2n-n-2)+3 =n+1 \ge 1$. So, we assume that $T =M$. 
In this case, we show that the following three conditions lead to a contradiction.
\begin{enumerate}
\item[(1)]
For every $1 \le i \le 3$, $k_i$ is divisible by $h_i$ and $k_i > m_i$.
\item[(2)]
$\sum_{i=1}^3\frac{k_i}{m_1}=n+2$.
\item[(3)]
The intersection matrix of $T'+\sum_{i=1}^3(\Phi_i-m_i\ol{F}_i)$ has a positive eigenvalue.
\end{enumerate}
The conditions (1) and (2) are equivalent to saying that $K_X$ is linearly equivalent to $0$. The 
condition (3) is equivalent to saying that $X$ is an affine surface. Let $\Phi$ represent one of 
the three singular fibers $\Phi_i$ with $1 \le i \le 3$. Then $\Phi$ is supported by $\ell'$ and 
the irreducible exceptional curves $E_1, \ldots, E_h$, where $E_h=\ol{F}$ is a unique $(-1)$ 
component. Note that $T=M$ and $M$ is untouched under the sequence $\sigma$ of blowing-ups. We 
shall determine the coefficient $\alpha$ in the expression $A:=
M+\alpha\ell'+\sum_{i=1}^{h-1}\beta_iE_i$, which is subject to the relation 
$\is{A}{\ell'}=\is{A}{E_i}=0$ for $1 \le i \le h-1$. A straightforward computation gives the following 
result.

\begin{lem}\label{Lemma 2.4.7}
We have the following table.
\[{\large
\begin{array}{|c|c|c|} \hline
m  &  k  & \alpha \\ \hline
2  &  h  & \frac{1}{4}h \\ \hline
3  &  h+1  &  \frac{h+3}{9} \\ \hline
3  &  h  &  \frac{h}{9} \\ \hline
4  &  h+s-1  &  \frac{(8s-3)(h-s)-(12s-5)}{16\{2(h-s)-3\}} \\ \hline
4  &  h+2  &  \frac{h+8}{16} \\ \hline
4  &  h  &  \frac{h}{16} \\ \hline
5  &  h+3  &  \frac{h+15}{25} \\ \hline
5  &  h+3  &  \frac{h+11}{25} \\ \hline
5  &  h+2  &  \frac{h+6}{25} \\ \hline
5  &  h  &  \frac{h}{25}  \\ \hline
\end{array}}
\]
\end{lem}

Let $\alpha_i$ be the coefficient of $\ell_i'$ as computed as above for the fiber $\Phi_i$ for 
$1 \le i \le 3$. Then the condition (3) is equivalent to $-n+\alpha_1+\alpha_2+\alpha_3 > 0$. We just 
indicate an outline of the proof in the case where $\{m_1,m_2,m_3\}=\{2,3,5\}, k_1=h_1, k_2=h_2$ 
and $k_3=h_3$. The conditions (1), (2) and (3) impose 
\begin{eqnarray*}
h_1 \ge 4, h_2 \ge 6, h_3 \ge 10, \\ 
\frac{h_1}{2}+\frac{h_2}{3}+\frac{h_3}{5}=n+2, \\
\frac{h_1}{4}+\frac{h_2}{9}+\frac{h_3}{25} > n. 
\end{eqnarray*}
Then we have 
\[
2> \frac{h_1}{4}+\frac{2}{9}h_2+\frac{4}{25}h_3 >3.
\]
This is a contradiction.
\QED

We excluded the case $\{m_1,m_2,m_3\}=\{2,2,m\}$ with $m > 5$. This is because we are not able to 
prove an inequality $\alpha m >k$ with the above notations. So, we make the following conjecture.

\begin{conj}\label{Conjecture 2.4.8}
Let $mF$ be a multiple fiber of a smooth $\A^1$-fiber space with 
irreducible fibers. Let $k$ be the coefficient of $F$ in the canonical divisor $K_X$. Define a rational 
number $\alpha$ by completing the fiber $F$ in a singular fiber of a $\BP^1$-fibration. Then we have 
$\alpha m > k$.
\end{conj}

\begin{remark}
{\em Theorem \ref{Theorem 2.4.6} asserts, in fact, the following result.}
\svskip

Let $X$ be either a $\Q$-homology plane with $\ML(X)=k$ or a Platonic $\A^1$-fiber space possibly 
except for the case where the multiplicity sequence is $\{2,2,m\}$ with $m > 5$. Then there are no 
\'etale morphisms from $X$ to an affine surface with trivial canonical divisor.
\end{remark}

We shall next consider the case of affine pseudo-planes. 

\begin{thm}\label{Theorem 2.4.10}
Let $X$ be an affine pseudo-plane of type $(d,n,r)$. If $r \ne 2$ then there are no \'etale 
morphism from $X$ to $\A^2$.
\end{thm}
\Proof
Let $\rho : X \to A$ be the given $\A^1$-fibration with a multiple fiber $dF_0$, where $A \cong \A^1$.
A simple computation using the smooth compactification $(V,D)$ in Definition \ref{Definition 1.4.5} 
shows that $\Pic(X)$ is isomorphic to $\Z/d\Z$ with a generator $F_0$ and the canonical divisor 
$K_X$ is linearly equivalent to $(r-2)F_0$. Hence $K_X \not\sim 0$ unless $r=2$. If there is an 
\'etale morphism $\varphi : X \to \A^2$, then $K_X$ is the inverse image of the canonical divisor of 
$\A^2$. Hence $K_X \sim 0$. So, we obtain the stated assertion if $r \ne 2$.
\QED

We note that the differential module $\Omega^1_X$ is the direct sum of an invertible sheaf and $\SO_X$ 
by virtue of a result due to Murthy-Swan \cite{MuSw}. So, $\Omega^1_X$ is trivial if and only if 
$K_X \sim 0$. 

\newpage

\printindex
\end{document}